\documentclass[12pt]{amsart}
\title[Affine links]{A note on affine links}
\author{Sadek AL HARBAT}
\address{LAMFA, Universit\'e de Picardie - Jules Verne} 
\email{sadikharbat@math.univ-paris-diderot.fr}

\usepackage{etex} 
\usepackage[latin1]{inputenc}	
\usepackage[T1]{fontenc}
\usepackage[english]{babel}
\usepackage{amsmath,amssymb}
\usepackage{wasysym}
\usepackage{stmaryrd}
\usepackage[table]{xcolor}
\usepackage{color}
\usepackage{dsfont}\let\mathbb\mathds
\usepackage[all]{xy}
\usepackage{graphicx}
\usepackage{mathenv}
\usepackage{lastpage}
\usepackage{fancyhdr}	
\usepackage{subfigure}
\usepackage{geometry}
\usepackage[version=3]{mhchem}
\usepackage[force,almostfull]{textcomp}
\usepackage{array}
\usepackage{lipsum}		
\usepackage{eso-pic}
\usepackage{multirow}
\usepackage{fancybox}
\usepackage{cite}
\usepackage{amsthm}
\usepackage{listings}
\usepackage{lscape}
\usepackage{multicol}		
\usepackage{setspace}
\usepackage{times}
\usepackage{eurosym}
\usepackage{latexsym}
\usepackage{ulem}
\usepackage{lmodern}
\usepackage{colortbl}
\usepackage{rotating}
\usepackage{type1cm}
\usepackage{lettrine}
\usepackage{ragged2e}
\usepackage{mathrsfs}
\usepackage{scrtime}
\usepackage{enumitem}

\usepackage{longtable}
\usepackage{url}
\usepackage{nomentbl}
\usepackage{nomencl}
\usepackage{hyperref}
\hypersetup{
pdfpagemode=UseOutlines,      
pdfstartview=Fit,             
pdffitwindow=true,            
pdfpagelayout=TwoColumnsRight,
pdftoolbar=true,              
pdfmenubar=true,              
bookmarksopen=false,          
bookmarksnumbered=true,       
colorlinks=true,              
pdfauthor={Ton nom},          
pdftitle={Titre PDF},         
pdfcreator=PDFLaTeX,          %
pdfproducer=PDFLaTeX,         %
linkcolor=blue,               
urlcolor=blue,                
anchorcolor=black,            
citecolor=blue,               
frenchlinks=true,             
pdfborder={0 0 0}             
}
\usepackage{pgf, tikz, tikz-cd}
\usetikzlibrary{matrix,positioning}
\usetikzlibrary{arrows,decorations.markings}

\definecolor{aliceblue}{rgb}{0.94,0.97,1.00}
\definecolor{antiquewhite}{rgb}{0.98,0.92,0.84}
\definecolor{antiquewhite1}{rgb}{1.00,0.94,0.86}
\definecolor{antiquewhite2}{rgb}{0.93,0.87,0.80}
\definecolor{antiquewhite3}{rgb}{0.80,0.75,0.69}
\definecolor{antiquewhite4}{rgb}{0.55,0.51,0.47}
\definecolor{aquamarine}{rgb}{0.50,1.00,0.83}
\definecolor{aquamarine1}{rgb}{0.50,1.00,0.83}
\definecolor{aquamarine2}{rgb}{0.46,0.93,0.78}
\definecolor{aquamarine3}{rgb}{0.40,0.80,0.67}
\definecolor{aquamarine4}{rgb}{0.27,0.55,0.45}
\definecolor{azure}{rgb}{0.94,1.00,1.00}
\definecolor{azure1}{rgb}{0.94,1.00,1.00}
\definecolor{azure2}{rgb}{0.88,0.93,0.93}
\definecolor{azure3}{rgb}{0.76,0.80,0.80}
\definecolor{azure4}{rgb}{0.51,0.55,0.55}
\definecolor{beige}{rgb}{0.96,0.96,0.86}
\definecolor{bisque}{rgb}{1.00,0.89,0.77}
\definecolor{bisque1}{rgb}{1.00,0.89,0.77}
\definecolor{bisque2}{rgb}{0.93,0.84,0.72}
\definecolor{bisque3}{rgb}{0.80,0.72,0.62}
\definecolor{bisque4}{rgb}{0.55,0.49,0.42}
\definecolor{black}{rgb}{0.00,0.00,0.00}
\definecolor{blanchedalmond}{rgb}{1.00,0.92,0.80}
\definecolor{blue}{rgb}{0.00,0.00,1.00}
\definecolor{blue1}{rgb}{0.00,0.00,1.00}
\definecolor{blue2}{rgb}{0.00,0.00,0.93}
\definecolor{blue3}{rgb}{0.00,0.00,0.80}
\definecolor{blue4}{rgb}{0.00,0.00,0.55}
\definecolor{blueviolet}{rgb}{0.54,0.17,0.89}
\definecolor{brown}{rgb}{0.65,0.16,0.16}
\definecolor{brown1}{rgb}{1.00,0.25,0.25}
\definecolor{brown2}{rgb}{0.93,0.23,0.23}
\definecolor{brown3}{rgb}{0.80,0.20,0.20}
\definecolor{brown4}{rgb}{0.55,0.14,0.14}
\definecolor{burlywood}{rgb}{0.87,0.72,0.53}
\definecolor{burlywood1}{rgb}{1.00,0.83,0.61}
\definecolor{burlywood2}{rgb}{0.93,0.77,0.57}
\definecolor{burlywood3}{rgb}{0.80,0.67,0.49}
\definecolor{burlywood4}{rgb}{0.55,0.45,0.33}
\definecolor{cadetblue}{rgb}{0.37,0.62,0.63}
\definecolor{cadetblue1}{rgb}{0.60,0.96,1.00}
\definecolor{cadetblue2}{rgb}{0.56,0.90,0.93}
\definecolor{cadetblue3}{rgb}{0.48,0.77,0.80}
\definecolor{cadetblue4}{rgb}{0.33,0.53,0.55}
\definecolor{chartreuse}{rgb}{0.50,1.00,0.00}
\definecolor{chartreuse1}{rgb}{0.50,1.00,0.00}
\definecolor{chartreuse2}{rgb}{0.46,0.93,0.00}
\definecolor{chartreuse3}{rgb}{0.40,0.80,0.00}
\definecolor{chartreuse4}{rgb}{0.27,0.55,0.00}
\definecolor{chocolate}{rgb}{0.82,0.41,0.12}
\definecolor{chocolate1}{rgb}{1.00,0.50,0.14}
\definecolor{chocolate2}{rgb}{0.93,0.46,0.13}
\definecolor{chocolate3}{rgb}{0.80,0.40,0.11}
\definecolor{chocolate4}{rgb}{0.55,0.27,0.07}
\definecolor{coral}{rgb}{1.00,0.50,0.31}
\definecolor{coral1}{rgb}{1.00,0.45,0.34}
\definecolor{coral2}{rgb}{0.93,0.42,0.31}
\definecolor{coral3}{rgb}{0.80,0.36,0.27}
\definecolor{coral4}{rgb}{0.55,0.24,0.18}
\definecolor{cornsilk}{rgb}{1.00,0.97,0.86}
\definecolor{cornflowerblue}{rgb}{0.39,0.58,0.93}
\definecolor{cornsilk1}{rgb}{1.00,0.97,0.86}
\definecolor{cornsilk2}{rgb}{0.93,0.91,0.80}
\definecolor{cornsilk3}{rgb}{0.80,0.78,0.69}
\definecolor{cornsilk4}{rgb}{0.55,0.53,0.47}
\definecolor{cyan}{rgb}{0.00,1.00,1.00}
\definecolor{cyan1}{rgb}{0.00,1.00,1.00}
\definecolor{cyan2}{rgb}{0.00,0.93,0.93}
\definecolor{cyan3}{rgb}{0.00,0.80,0.80}
\definecolor{cyan4}{rgb}{0.00,0.55,0.55}
\definecolor{darkblue}{rgb}{0.00,0.00,0.55}
\definecolor{darkcyan}{rgb}{0.00,0.55,0.55}
\definecolor{darkgoldenrod}{rgb}{0.72,0.53,0.04}
\definecolor{DarkGoldenrod1}{rgb}{1.00,0.73,0.06}
\definecolor{DarkGoldenrod2}{rgb}{0.93,0.68,0.05}
\definecolor{DarkGoldenrod3}{rgb}{0.80,0.58,0.05}
\definecolor{DarkGoldenrod4}{rgb}{0.55,0.40,0.03}
\definecolor{darkgray}{rgb}{0.66,0.66,0.66}
\definecolor{darkgreen}{rgb}{0.00,0.39,0.00}
\definecolor{darkgrey}{rgb}{0.66,0.66,0.66}
\definecolor{darkkhaki}{rgb}{0.74,0.72,0.42}
\definecolor{darkmagenta}{rgb}{0.55,0.00,0.55}
\definecolor{darkOlivegreen}{rgb}{0.33,0.42,0.18}
\definecolor{darkOlivegreen1}{rgb}{0.79,1.00,0.44}
\definecolor{darkOlivegreen2}{rgb}{0.74,0.93,0.41}
\definecolor{darkOlivegreen3}{rgb}{0.64,0.80,0.35}
\definecolor{darkOlivegreen4}{rgb}{0.43,0.55,0.24}
\definecolor{darkolive}{rgb}{0.33,0.42,0.18}
\definecolor{DarkOrange1}{rgb}{1.00,0.50,0.00}
\definecolor{DarkOrange2}{rgb}{0.93,0.46,0.00}
\definecolor{DarkOrange3}{rgb}{0.80,0.40,0.00}
\definecolor{DarkOrange4}{rgb}{0.55,0.27,0.00}
\definecolor{darkorange}{rgb}{1.00,0.55,0.00}
\definecolor{darkorchid}{rgb}{0.60,0.20,0.80}
\definecolor{darkOrchid1}{rgb}{0.75,0.24,1.00}
\definecolor{darkOrchid2}{rgb}{0.70,0.23,0.93}
\definecolor{darkOrchid3}{rgb}{0.60,0.20,0.80}
\definecolor{darkOrchid4}{rgb}{0.41,0.13,0.55}
\definecolor{darkred}{rgb}{0.55,0.00,0.00}
\definecolor{darksalmon}{rgb}{0.91,0.59,0.48}
\definecolor{darkseagreen}{rgb}{0.56,0.74,0.56}
\definecolor{darkseagreen1}{rgb}{0.76,1.00,0.76}
\definecolor{darkseagreen2}{rgb}{0.71,0.93,0.71}
\definecolor{darkseagreen3}{rgb}{0.61,0.80,0.61}
\definecolor{darkseagreen4}{rgb}{0.41,0.55,0.41}
\definecolor{darksea}{rgb}{0.56,0.74,0.56}
\definecolor{darkslategray}{rgb}{0.18,0.31,0.31}
\definecolor{darkslateblue}{rgb}{0.28,0.24,0.55}
\definecolor{darkslategray1}{rgb}{0.59,1.00,1.00}
\definecolor{darkslategray2}{rgb}{0.55,0.93,0.93}
\definecolor{darkslategray3}{rgb}{0.47,0.80,0.80}
\definecolor{darkslategray4}{rgb}{0.32,0.55,0.55}
\definecolor{darkslate}{rgb}{0.18,0.31,0.31}
\definecolor{darkslate1}{rgb}{0.28,0.24,0.55}
\definecolor{darkturquoise}{rgb}{0.00,0.81,0.82}
\definecolor{darkviolet}{rgb}{0.58,0.00,0.83}
\definecolor{deeppink}{rgb}{1.00,0.08,0.58}
\definecolor{deepPink1}{rgb}{1.00,0.08,0.58}
\definecolor{deepPink2}{rgb}{0.93,0.07,0.54}
\definecolor{deepPink3}{rgb}{0.80,0.06,0.46}
\definecolor{deepPink4}{rgb}{0.55,0.04,0.31}
\definecolor{deepskyblue}{rgb}{0.00,0.75,1.00}
\definecolor{deepskyblue1}{rgb}{0.00,0.75,1.00}
\definecolor{deepskyblue2}{rgb}{0.00,0.70,0.93}
\definecolor{deepskyblue3}{rgb}{0.00,0.60,0.80}
\definecolor{deepskyblue4}{rgb}{0.00,0.41,0.55}
\definecolor{deepsky}{rgb}{0.00,0.75,1.00}
\definecolor{dimgray}{rgb}{0.41,0.41,0.41}
\definecolor{dodgerblue}{rgb}{0.12,0.56,1.00}
\definecolor{dodgerblue1}{rgb}{0.12,0.56,1.00}
\definecolor{dodgerblue2}{rgb}{0.11,0.53,0.93}
\definecolor{dodgerblue3}{rgb}{0.09,0.45,0.80}
\definecolor{dodgerblue4}{rgb}{0.06,0.31,0.55}
\definecolor{firebrick}{rgb}{0.70,0.13,0.13}
\definecolor{firebrick1}{rgb}{1.00,0.19,0.19}
\definecolor{firebrick2}{rgb}{0.93,0.17,0.17}
\definecolor{firebrick3}{rgb}{0.80,0.15,0.15}
\definecolor{firebrick4}{rgb}{0.55,0.10,0.10}
\definecolor{floralwhite}{rgb}{1.00,0.98,0.94}
\definecolor{forestgreen}{rgb}{0.13,0.55,0.13}
\definecolor{gainsboro}{rgb}{0.86,0.86,0.86}
\definecolor{ghostwhite}{rgb}{0.97,0.97,1.00}
\definecolor{gold1}{rgb}{1.00,0.84,0.00}
\definecolor{gold2}{rgb}{0.93,0.79,0.00}
\definecolor{gold3}{rgb}{0.80,0.68,0.00}
\definecolor{gold4}{rgb}{0.55,0.46,0.00}
\definecolor{goldenrod}{rgb}{0.85,0.65,0.13}
\definecolor{goldenrod1}{rgb}{1.00,0.76,0.15}
\definecolor{goldenrod2}{rgb}{0.93,0.71,0.13}
\definecolor{goldenrod3}{rgb}{0.80,0.61,0.11}
\definecolor{goldenrod4}{rgb}{0.55,0.41,0.08}
\definecolor{gold}{rgb}{1.00,0.84,0.00}
\definecolor{gray}{rgb}{0.75,0.75,0.75}
\definecolor{gray0}{rgb}{0.00,0.00,0.00}
\definecolor{gray1}{rgb}{0.01,0.01,0.01}
\definecolor{gray2}{rgb}{0.02,0.02,0.02}
\definecolor{gray3}{rgb}{0.03,0.03,0.03}
\definecolor{gray4}{rgb}{0.04,0.04,0.04}
\definecolor{gray5}{rgb}{0.05,0.05,0.05}
\definecolor{gray6}{rgb}{0.06,0.06,0.06}
\definecolor{gray7}{rgb}{0.07,0.07,0.07}
\definecolor{gray8}{rgb}{0.08,0.08,0.08}
\definecolor{gray9}{rgb}{0.09,0.09,0.09}
\definecolor{gray10}{rgb}{0.10,0.10,0.10}
\definecolor{gray11}{rgb}{0.11,0.11,0.11}
\definecolor{gray12}{rgb}{0.12,0.12,0.12}
\definecolor{gray13}{rgb}{0.13,0.13,0.13}
\definecolor{gray14}{rgb}{0.14,0.14,0.14}
\definecolor{gray15}{rgb}{0.15,0.15,0.15}
\definecolor{gray16}{rgb}{0.16,0.16,0.16}
\definecolor{gray17}{rgb}{0.17,0.17,0.17}
\definecolor{gray18}{rgb}{0.18,0.18,0.18}
\definecolor{gray19}{rgb}{0.19,0.19,0.19}
\definecolor{gray20}{rgb}{0.20,0.20,0.20}
\definecolor{gray21}{rgb}{0.21,0.21,0.21}
\definecolor{gray22}{rgb}{0.22,0.22,0.22}
\definecolor{gray23}{rgb}{0.23,0.23,0.23}
\definecolor{gray24}{rgb}{0.24,0.24,0.24}
\definecolor{gray25}{rgb}{0.25,0.25,0.25}
\definecolor{gray26}{rgb}{0.26,0.26,0.26}
\definecolor{gray27}{rgb}{0.27,0.27,0.27}
\definecolor{gray28}{rgb}{0.28,0.28,0.28}
\definecolor{gray29}{rgb}{0.29,0.29,0.29}
\definecolor{gray30}{rgb}{0.30,0.30,0.30}
\definecolor{gray31}{rgb}{0.31,0.31,0.31}
\definecolor{gray32}{rgb}{0.32,0.32,0.32}
\definecolor{gray33}{rgb}{0.33,0.33,0.33}
\definecolor{gray34}{rgb}{0.34,0.34,0.34}
\definecolor{gray35}{rgb}{0.35,0.35,0.35}
\definecolor{gray36}{rgb}{0.36,0.36,0.36}
\definecolor{gray37}{rgb}{0.37,0.37,0.37}
\definecolor{gray38}{rgb}{0.38,0.38,0.38}
\definecolor{gray39}{rgb}{0.39,0.39,0.39}
\definecolor{gray40}{rgb}{0.40,0.40,0.40}
\definecolor{gray41}{rgb}{0.41,0.41,0.41}
\definecolor{gray42}{rgb}{0.42,0.42,0.42}
\definecolor{gray43}{rgb}{0.43,0.43,0.43}
\definecolor{gray44}{rgb}{0.44,0.44,0.44}
\definecolor{gray45}{rgb}{0.45,0.45,0.45}
\definecolor{gray46}{rgb}{0.46,0.46,0.46}
\definecolor{gray47}{rgb}{0.47,0.47,0.47}
\definecolor{gray48}{rgb}{0.48,0.48,0.48}
\definecolor{gray49}{rgb}{0.49,0.49,0.49}
\definecolor{gray50}{rgb}{0.50,0.50,0.50}
\definecolor{gray51}{rgb}{0.51,0.51,0.51}
\definecolor{gray52}{rgb}{0.52,0.52,0.52}
\definecolor{gray53}{rgb}{0.53,0.53,0.53}
\definecolor{gray54}{rgb}{0.54,0.54,0.54}
\definecolor{gray55}{rgb}{0.55,0.55,0.55}
\definecolor{gray56}{rgb}{0.56,0.56,0.56}
\definecolor{gray57}{rgb}{0.57,0.57,0.57}
\definecolor{gray58}{rgb}{0.58,0.58,0.58}
\definecolor{gray59}{rgb}{0.59,0.59,0.59}
\definecolor{gray60}{rgb}{0.60,0.60,0.60}
\definecolor{gray61}{rgb}{0.61,0.61,0.61}
\definecolor{gray62}{rgb}{0.62,0.62,0.62}
\definecolor{gray63}{rgb}{0.63,0.63,0.63}
\definecolor{gray64}{rgb}{0.64,0.64,0.64}
\definecolor{gray65}{rgb}{0.65,0.65,0.65}
\definecolor{gray66}{rgb}{0.66,0.66,0.66}
\definecolor{gray67}{rgb}{0.67,0.67,0.67}
\definecolor{gray68}{rgb}{0.68,0.68,0.68}
\definecolor{gray69}{rgb}{0.69,0.69,0.69}
\definecolor{gray70}{rgb}{0.70,0.70,0.70}
\definecolor{gray71}{rgb}{0.71,0.71,0.71}
\definecolor{gray72}{rgb}{0.72,0.72,0.72}
\definecolor{gray73}{rgb}{0.73,0.73,0.73}
\definecolor{gray74}{rgb}{0.74,0.74,0.74}
\definecolor{gray75}{rgb}{0.75,0.75,0.75}
\definecolor{gray76}{rgb}{0.76,0.76,0.76}
\definecolor{gray77}{rgb}{0.77,0.77,0.77}
\definecolor{gray78}{rgb}{0.78,0.78,0.78}
\definecolor{gray79}{rgb}{0.79,0.79,0.79}
\definecolor{gray80}{rgb}{0.80,0.80,0.80}
\definecolor{gray81}{rgb}{0.81,0.81,0.81}
\definecolor{gray82}{rgb}{0.82,0.82,0.82}
\definecolor{gray83}{rgb}{0.83,0.83,0.83}
\definecolor{gray84}{rgb}{0.84,0.84,0.84}
\definecolor{gray85}{rgb}{0.85,0.85,0.85}
\definecolor{gray86}{rgb}{0.86,0.86,0.86}
\definecolor{gray87}{rgb}{0.87,0.87,0.87}
\definecolor{gray88}{rgb}{0.88,0.88,0.88}
\definecolor{gray89}{rgb}{0.89,0.89,0.89}
\definecolor{gray90}{rgb}{0.90,0.90,0.90}
\definecolor{gray91}{rgb}{0.91,0.91,0.91}
\definecolor{gray92}{rgb}{0.92,0.92,0.92}
\definecolor{gray93}{rgb}{0.93,0.93,0.93}
\definecolor{gray94}{rgb}{0.94,0.94,0.94}
\definecolor{gray95}{rgb}{0.95,0.95,0.95}
\definecolor{gray96}{rgb}{0.96,0.96,0.96}
\definecolor{gray97}{rgb}{0.97,0.97,0.97}
\definecolor{gray98}{rgb}{0.98,0.98,0.98}
\definecolor{gray99}{rgb}{0.99,0.99,0.99}
\definecolor{gray100}{rgb}{1.00,1.00,1.00}
\definecolor{green}{rgb}{0.00,1.00,0.00}
\definecolor{green1}{rgb}{0.00,1.00,0.00}
\definecolor{green2}{rgb}{0.00,0.93,0.00}
\definecolor{green3}{rgb}{0.00,0.80,0.00}
\definecolor{green4}{rgb}{0.00,0.55,0.00}
\definecolor{greenyellow}{rgb}{0.68,1.00,0.18}
\definecolor{grey}{rgb}{0.75,0.75,0.75}
\definecolor{grey0}{rgb}{0.00,0.00,0.00}
\definecolor{grey1}{rgb}{0.01,0.01,0.01}
\definecolor{grey2}{rgb}{0.02,0.02,0.02}
\definecolor{grey3}{rgb}{0.03,0.03,0.03}
\definecolor{grey4}{rgb}{0.04,0.04,0.04}
\definecolor{grey5}{rgb}{0.05,0.05,0.05}
\definecolor{grey6}{rgb}{0.06,0.06,0.06}
\definecolor{grey7}{rgb}{0.07,0.07,0.07}
\definecolor{grey8}{rgb}{0.08,0.08,0.08}
\definecolor{grey9}{rgb}{0.09,0.09,0.09}
\definecolor{grey10}{rgb}{0.10,0.10,0.10}
\definecolor{grey11}{rgb}{0.11,0.11,0.11}
\definecolor{grey12}{rgb}{0.12,0.12,0.12}
\definecolor{grey13}{rgb}{0.13,0.13,0.13}
\definecolor{grey14}{rgb}{0.14,0.14,0.14}
\definecolor{grey15}{rgb}{0.15,0.15,0.15}
\definecolor{grey16}{rgb}{0.16,0.16,0.16}
\definecolor{grey17}{rgb}{0.17,0.17,0.17}
\definecolor{grey18}{rgb}{0.18,0.18,0.18}
\definecolor{grey19}{rgb}{0.19,0.19,0.19}
\definecolor{grey20}{rgb}{0.20,0.20,0.20}
\definecolor{grey21}{rgb}{0.21,0.21,0.21}
\definecolor{grey22}{rgb}{0.22,0.22,0.22}
\definecolor{grey23}{rgb}{0.23,0.23,0.23}
\definecolor{grey24}{rgb}{0.24,0.24,0.24}
\definecolor{grey25}{rgb}{0.25,0.25,0.25}
\definecolor{grey26}{rgb}{0.26,0.26,0.26}
\definecolor{grey27}{rgb}{0.27,0.27,0.27}
\definecolor{grey28}{rgb}{0.28,0.28,0.28}
\definecolor{grey29}{rgb}{0.29,0.29,0.29}
\definecolor{grey30}{rgb}{0.30,0.30,0.30}
\definecolor{grey31}{rgb}{0.31,0.31,0.31}
\definecolor{grey32}{rgb}{0.32,0.32,0.32}
\definecolor{grey33}{rgb}{0.33,0.33,0.33}
\definecolor{grey34}{rgb}{0.34,0.34,0.34}
\definecolor{grey35}{rgb}{0.35,0.35,0.35}
\definecolor{grey36}{rgb}{0.36,0.36,0.36}
\definecolor{grey37}{rgb}{0.37,0.37,0.37}
\definecolor{grey38}{rgb}{0.38,0.38,0.38}
\definecolor{grey39}{rgb}{0.39,0.39,0.39}
\definecolor{grey40}{rgb}{0.40,0.40,0.40}
\definecolor{grey41}{rgb}{0.41,0.41,0.41}
\definecolor{grey42}{rgb}{0.42,0.42,0.42}
\definecolor{grey43}{rgb}{0.43,0.43,0.43}
\definecolor{grey44}{rgb}{0.44,0.44,0.44}
\definecolor{grey45}{rgb}{0.45,0.45,0.45}
\definecolor{grey46}{rgb}{0.46,0.46,0.46}
\definecolor{grey47}{rgb}{0.47,0.47,0.47}
\definecolor{grey48}{rgb}{0.48,0.48,0.48}
\definecolor{grey49}{rgb}{0.49,0.49,0.49}
\definecolor{grey50}{rgb}{0.50,0.50,0.50}
\definecolor{grey51}{rgb}{0.51,0.51,0.51}
\definecolor{grey52}{rgb}{0.52,0.52,0.52}
\definecolor{grey53}{rgb}{0.53,0.53,0.53}
\definecolor{grey54}{rgb}{0.54,0.54,0.54}
\definecolor{grey55}{rgb}{0.55,0.55,0.55}
\definecolor{grey56}{rgb}{0.56,0.56,0.56}
\definecolor{grey57}{rgb}{0.57,0.57,0.57}
\definecolor{grey58}{rgb}{0.58,0.58,0.58}
\definecolor{grey59}{rgb}{0.59,0.59,0.59}
\definecolor{grey60}{rgb}{0.60,0.60,0.60}
\definecolor{grey61}{rgb}{0.61,0.61,0.61}
\definecolor{grey62}{rgb}{0.62,0.62,0.62}
\definecolor{grey63}{rgb}{0.63,0.63,0.63}
\definecolor{grey64}{rgb}{0.64,0.64,0.64}
\definecolor{grey65}{rgb}{0.65,0.65,0.65}
\definecolor{grey66}{rgb}{0.66,0.66,0.66}
\definecolor{grey67}{rgb}{0.67,0.67,0.67}
\definecolor{grey68}{rgb}{0.68,0.68,0.68}
\definecolor{grey69}{rgb}{0.69,0.69,0.69}
\definecolor{grey70}{rgb}{0.70,0.70,0.70}
\definecolor{grey71}{rgb}{0.71,0.71,0.71}
\definecolor{grey72}{rgb}{0.72,0.72,0.72}
\definecolor{grey73}{rgb}{0.73,0.73,0.73}
\definecolor{grey74}{rgb}{0.74,0.74,0.74}
\definecolor{grey75}{rgb}{0.75,0.75,0.75}
\definecolor{grey76}{rgb}{0.76,0.76,0.76}
\definecolor{grey77}{rgb}{0.77,0.77,0.77}
\definecolor{grey78}{rgb}{0.78,0.78,0.78}
\definecolor{grey79}{rgb}{0.79,0.79,0.79}
\definecolor{grey80}{rgb}{0.80,0.80,0.80}
\definecolor{grey81}{rgb}{0.81,0.81,0.81}
\definecolor{grey82}{rgb}{0.82,0.82,0.82}
\definecolor{grey83}{rgb}{0.83,0.83,0.83}
\definecolor{grey84}{rgb}{0.84,0.84,0.84}
\definecolor{grey85}{rgb}{0.85,0.85,0.85}
\definecolor{grey86}{rgb}{0.86,0.86,0.86}
\definecolor{grey87}{rgb}{0.87,0.87,0.87}
\definecolor{grey88}{rgb}{0.88,0.88,0.88}
\definecolor{grey89}{rgb}{0.89,0.89,0.89}
\definecolor{grey90}{rgb}{0.90,0.90,0.90}
\definecolor{grey91}{rgb}{0.91,0.91,0.91}
\definecolor{grey92}{rgb}{0.92,0.92,0.92}
\definecolor{grey93}{rgb}{0.93,0.93,0.93}
\definecolor{grey94}{rgb}{0.94,0.94,0.94}
\definecolor{grey95}{rgb}{0.95,0.95,0.95}
\definecolor{grey96}{rgb}{0.96,0.96,0.96}
\definecolor{grey97}{rgb}{0.97,0.97,0.97}
\definecolor{grey98}{rgb}{0.98,0.98,0.98}
\definecolor{grey99}{rgb}{0.99,0.99,0.99}
\definecolor{grey100}{rgb}{1.00,1.00,1.00}
\definecolor{honeydew}{rgb}{0.94,1.00,0.94}
\definecolor{honeydew1}{rgb}{0.94,1.00,0.94}
\definecolor{honeydew2}{rgb}{0.88,0.93,0.88}
\definecolor{honeydew3}{rgb}{0.76,0.80,0.76}
\definecolor{honeydew4}{rgb}{0.51,0.55,0.51}
\definecolor{hotpink}{rgb}{1.00,0.41,0.71}
\definecolor{hotPink1}{rgb}{1.00,0.43,0.71}
\definecolor{hotPink2}{rgb}{0.93,0.42,0.65}
\definecolor{hotPink3}{rgb}{0.80,0.38,0.56}
\definecolor{hotPink4}{rgb}{0.55,0.23,0.38}
\definecolor{indianred}{rgb}{0.80,0.36,0.36}
\definecolor{indianred1}{rgb}{1.00,0.42,0.42}
\definecolor{indianred2}{rgb}{0.93,0.39,0.39}
\definecolor{indianred3}{rgb}{0.80,0.33,0.33}
\definecolor{indianred4}{rgb}{0.55,0.23,0.23}
\definecolor{ivory}{rgb}{1.00,1.00,0.94}
\definecolor{ivory1}{rgb}{1.00,1.00,0.94}
\definecolor{ivory2}{rgb}{0.93,0.93,0.88}
\definecolor{ivory3}{rgb}{0.80,0.80,0.76}
\definecolor{ivory4}{rgb}{0.55,0.55,0.51}
\definecolor{khaki}{rgb}{0.94,0.90,0.55}
\definecolor{khaki1}{rgb}{1.00,0.96,0.56}
\definecolor{khaki2}{rgb}{0.93,0.90,0.52}
\definecolor{khaki3}{rgb}{0.80,0.78,0.45}
\definecolor{khaki4}{rgb}{0.55,0.53,0.31}
\definecolor{lavenderblush}{rgb}{1.00,0.94,0.96}
\definecolor{lavenderblush1}{rgb}{1.00,0.94,0.96}
\definecolor{lavenderblush2}{rgb}{0.93,0.88,0.90}
\definecolor{lavenderblush3}{rgb}{0.80,0.76,0.77}
\definecolor{lavenderblush4}{rgb}{0.55,0.51,0.53}
\definecolor{lavender}{rgb}{0.90,0.90,0.98}
\definecolor{lawngreen}{rgb}{0.49,0.99,0.00}
\definecolor{lemonchiffon}{rgb}{1.00,0.98,0.80}
\definecolor{lemonchiffon1}{rgb}{1.00,0.98,0.80}
\definecolor{lemonchiffon2}{rgb}{0.93,0.91,0.75}
\definecolor{lemonchiffon3}{rgb}{0.80,0.79,0.65}
\definecolor{lemonchiffon4}{rgb}{0.55,0.54,0.44}
\definecolor{lightblue}{rgb}{0.68,0.85,0.90}
\definecolor{lightblue1}{rgb}{0.75,0.94,1.00}
\definecolor{lightblue2}{rgb}{0.70,0.87,0.93}
\definecolor{lightblue3}{rgb}{0.60,0.75,0.80}
\definecolor{lightblue4}{rgb}{0.41,0.51,0.55}
\definecolor{lightcoral}{rgb}{0.94,0.50,0.50}
\definecolor{lightcyan}{rgb}{0.88,1.00,1.00}
\definecolor{lightcyan1}{rgb}{0.88,1.00,1.00}
\definecolor{lightcyan2}{rgb}{0.82,0.93,0.93}
\definecolor{lightcyan3}{rgb}{0.71,0.80,0.80}
\definecolor{lightcyan4}{rgb}{0.48,0.55,0.55}
\definecolor{lightgoldenrod}{rgb}{0.93,0.87,0.51}
\definecolor{lightgoldenrod0}{rgb}{0.98,0.98,0.82}
\definecolor{lightgoldenrod1}{rgb}{1.00,0.93,0.55}
\definecolor{lightgoldenrod2}{rgb}{0.93,0.86,0.51}
\definecolor{lightgoldenrod3}{rgb}{0.80,0.75,0.44}
\definecolor{lightgoldenrod4}{rgb}{0.55,0.51,0.30}
\definecolor{lightgoldenrodYellow}{rgb}{0.98,0.98,0.82}
\definecolor{lightgray}{rgb}{0.83,0.83,0.83}
\definecolor{lightgreen}{rgb}{0.56,0.93,0.56}
\definecolor{lightgrey}{rgb}{0.83,0.83,0.83}
\definecolor{lightpink}{rgb}{1.00,0.71,0.76}
\definecolor{lightpink1}{rgb}{1.00,0.68,0.73}
\definecolor{lightpink2}{rgb}{0.93,0.64,0.68}
\definecolor{lightpink3}{rgb}{0.80,0.55,0.58}
\definecolor{lightpink4}{rgb}{0.55,0.37,0.40}
\definecolor{lightsalmon}{rgb}{1.00,0.63,0.48}
\definecolor{lightsalmon1}{rgb}{1.00,0.63,0.48}
\definecolor{lightsalmon2}{rgb}{0.93,0.58,0.45}
\definecolor{lightsalmon3}{rgb}{0.80,0.51,0.38}
\definecolor{lightsalmon4}{rgb}{0.55,0.34,0.26}
\definecolor{lightseagreen}{rgb}{0.13,0.70,0.67}
\definecolor{lightsea}{rgb}{0.13,0.70,0.67}
\definecolor{lightsky}{rgb}{0.53,0.81,0.98}
\definecolor{lightSkyblue}{rgb}{0.53,0.81,0.98}
\definecolor{lightSkyblue1}{rgb}{0.69,0.89,1.00}
\definecolor{lightSkyblue2}{rgb}{0.64,0.83,0.93}
\definecolor{lightSkyblue3}{rgb}{0.55,0.71,0.80}
\definecolor{lightSkyblue4}{rgb}{0.38,0.48,0.55}
\definecolor{lightslateblue}{rgb}{0.52,0.44,1.00}
\definecolor{lightslategray}{rgb}{0.47,0.53,0.60}
\definecolor{lightslate}{rgb}{0.47,0.53,0.60}
\definecolor{lightslate1}{rgb}{0.52,0.44,1.00}
\definecolor{lightsteelblue}{rgb}{0.69,0.77,0.87}
\definecolor{lightsteelblue1}{rgb}{0.79,0.88,1.00}
\definecolor{lightsteelblue2}{rgb}{0.74,0.82,0.93}
\definecolor{lightsteelblue3}{rgb}{0.64,0.71,0.80}
\definecolor{lightsteelblue4}{rgb}{0.43,0.48,0.55}
\definecolor{lightyellow}{rgb}{1.00,1.00,0.88}
\definecolor{lightsteel}{rgb}{0.69,0.77,0.87}
\definecolor{lightyellow}{rgb}{1.00,1.00,0.88}
\definecolor{lightyellow1}{rgb}{1.00,1.00,0.88}
\definecolor{lightyellow2}{rgb}{0.93,0.93,0.82}
\definecolor{lightyellow3}{rgb}{0.80,0.80,0.71}
\definecolor{Lightyellow4}{rgb}{0.55,0.55,0.48}
\definecolor{limegreen}{rgb}{0.20,0.80,0.20}
\definecolor{linen}{rgb}{0.98,0.94,0.90}
\definecolor{magenta}{rgb}{1.00,0.00,1.00}
\definecolor{magenta1}{rgb}{1.00,0.00,1.00}
\definecolor{magenta2}{rgb}{0.93,0.00,0.93}
\definecolor{magenta3}{rgb}{0.80,0.00,0.80}
\definecolor{magenta4}{rgb}{0.55,0.00,0.55}
\definecolor{maroon}{rgb}{0.69,0.19,0.38}
\definecolor{maroon1}{rgb}{1.00,0.20,0.70}
\definecolor{maroon2}{rgb}{0.93,0.19,0.65}
\definecolor{maroon3}{rgb}{0.80,0.16,0.56}
\definecolor{maroon4}{rgb}{0.55,0.11,0.38}
\definecolor{mediumaquamarine}{rgb}{0.40,0.80,0.67}
\definecolor{mediumblue}{rgb}{0.00,0.00,0.80}
\definecolor{mediumorchid1}{rgb}{0.88,0.40,1.00}
\definecolor{mediumorchid2}{rgb}{0.82,0.37,0.93}
\definecolor{mediumorchid3}{rgb}{0.71,0.32,0.80}
\definecolor{mediumorchid4}{rgb}{0.48,0.22,0.55}
\definecolor{mediumorchid}{rgb}{0.73,0.33,0.83}
\definecolor{mediumpurple}{rgb}{0.58,0.44,0.86}
\definecolor{mediumpurple1}{rgb}{0.67,0.51,1.00}
\definecolor{mediumpurple2}{rgb}{0.62,0.47,0.93}
\definecolor{mediumpurple3}{rgb}{0.54,0.41,0.80}
\definecolor{mediumpurple4}{rgb}{0.36,0.28,0.55}
\definecolor{medium}{rgb}{0.45,0.45,0.45}
\definecolor{mediumseagreen}{rgb}{0.24,0.70,0.44}
\definecolor{mediumsea}{rgb}{0.24,0.70,0.44}
\definecolor{mediumslateblue}{rgb}{0.48,0.41,0.93}
\definecolor{mediumslate}{rgb}{0.48,0.41,0.93}
\definecolor{mediumspringgreen}{rgb}{0.00,0.98,0.60}
\definecolor{mediumspring}{rgb}{0.00,0.98,0.60}
\definecolor{mediumturquoise}{rgb}{0.28,0.82,0.80}
\definecolor{mediumvioletred}{rgb}{0.78,0.08,0.52}
\definecolor{mediumviolet}{rgb}{0.78,0.08,0.52}
\definecolor{midnightblue}{rgb}{0.10,0.10,0.44}
\definecolor{mintcream}{rgb}{0.96,1.00,0.98}
\definecolor{mistyrose}{rgb}{1.00,0.89,0.88}
\definecolor{mistyrose1}{rgb}{1.00,0.89,0.88}
\definecolor{mistyrose2}{rgb}{0.93,0.84,0.82}
\definecolor{mistyrose3}{rgb}{0.80,0.72,0.71}
\definecolor{mistyrose4}{rgb}{0.55,0.49,0.48}
\definecolor{moccasin}{rgb}{1.00,0.89,0.71}
\definecolor{navajowhite}{rgb}{1.00,0.87,0.68}
\definecolor{navajowhite1}{rgb}{1.00,0.87,0.68}
\definecolor{navajowhite2}{rgb}{0.93,0.81,0.63}
\definecolor{navajowhite3}{rgb}{0.80,0.70,0.55}
\definecolor{navajowhite4}{rgb}{0.55,0.47,0.37}
\definecolor{navyblue}{rgb}{0.00,0.00,0.50}
\definecolor{navy}{rgb}{0.00,0.00,0.50}
\definecolor{oldlace}{rgb}{0.99,0.96,0.90}
\definecolor{olivedrab}{rgb}{0.42,0.56,0.14}
\definecolor{olivedrab1}{rgb}{0.75,1.00,0.24}
\definecolor{olivedrab2}{rgb}{0.70,0.93,0.23}
\definecolor{olivedrab3}{rgb}{0.60,0.80,0.20}
\definecolor{olivedrab4}{rgb}{0.41,0.55,0.13}
\definecolor{orange1}{rgb}{1.00,0.65,0.00}
\definecolor{orange2}{rgb}{0.93,0.60,0.00}
\definecolor{orange3}{rgb}{0.80,0.52,0.00}
\definecolor{orange4}{rgb}{0.55,0.35,0.00}
\definecolor{orangered}{rgb}{1.00,0.27,0.00}
\definecolor{orangered1}{rgb}{1.00,0.27,0.00}
\definecolor{orangered2}{rgb}{0.93,0.25,0.00}
\definecolor{orangered3}{rgb}{0.80,0.22,0.00}
\definecolor{orangered4}{rgb}{0.55,0.15,0.00}
\definecolor{orange}{rgb}{1.00,0.65,0.00}
\definecolor{orchid}{rgb}{0.85,0.44,0.84}
\definecolor{orchid1}{rgb}{1.00,0.51,0.98}
\definecolor{orchid2}{rgb}{0.93,0.48,0.91}
\definecolor{orchid3}{rgb}{0.80,0.41,0.79}
\definecolor{orchid4}{rgb}{0.55,0.28,0.54}
\definecolor{palegoldenrod}{rgb}{0.93,0.91,0.67}
\definecolor{palegreen}{rgb}{0.60,0.98,0.60}
\definecolor{palegreen1}{rgb}{0.60,1.00,0.60}
\definecolor{palegreen2}{rgb}{0.56,0.93,0.56}
\definecolor{palegreen3}{rgb}{0.49,0.80,0.49}
\definecolor{palegreen4}{rgb}{0.33,0.55,0.33}
\definecolor{paleturquoise1}{rgb}{0.73,1.00,1.00}
\definecolor{paleturquoise2}{rgb}{0.68,0.93,0.93}
\definecolor{paleturquoise3}{rgb}{0.59,0.80,0.80}
\definecolor{paleturquoise4}{rgb}{0.40,0.55,0.55}
\definecolor{paleturquoise}{rgb}{0.69,0.93,0.93}
\definecolor{palevioletred}{rgb}{0.86,0.44,0.58}
\definecolor{palevioletred1}{rgb}{1.00,0.51,0.67}
\definecolor{palevioletred2}{rgb}{0.93,0.47,0.62}
\definecolor{palevioletred3}{rgb}{0.80,0.41,0.54}
\definecolor{palevioletred4}{rgb}{0.55,0.28,0.36}
\definecolor{paleviolet}{rgb}{0.86,0.44,0.58}
\definecolor{papayawhip}{rgb}{1.00,0.94,0.84}
\definecolor{peachPuff1}{rgb}{1.00,0.85,0.73}
\definecolor{peachPuff2}{rgb}{0.93,0.80,0.68}
\definecolor{peachPuff3}{rgb}{0.80,0.69,0.58}
\definecolor{peachPuff4}{rgb}{0.55,0.47,0.40}
\definecolor{peachpuff}{rgb}{1.00,0.85,0.73}
\definecolor{peru}{rgb}{0.80,0.52,0.25}
\definecolor{pink}{rgb}{1.00,0.75,0.80}
\definecolor{pink1}{rgb}{1.00,0.71,0.77}
\definecolor{pink2}{rgb}{0.93,0.66,0.72}
\definecolor{pink3}{rgb}{0.80,0.57,0.62}
\definecolor{pink4}{rgb}{0.55,0.39,0.42}
\definecolor{plum}{rgb}{0.87,0.63,0.87}
\definecolor{plum1}{rgb}{1.00,0.73,1.00}
\definecolor{plum2}{rgb}{0.93,0.68,0.93}
\definecolor{plum3}{rgb}{0.80,0.59,0.80}
\definecolor{plum4}{rgb}{0.55,0.40,0.55}
\definecolor{powderblue}{rgb}{0.69,0.88,0.90}
\definecolor{purple}{rgb}{0.63,0.13,0.94}
\definecolor{purple1}{rgb}{0.61,0.19,1.00}
\definecolor{purple2}{rgb}{0.57,0.17,0.93}
\definecolor{purple3}{rgb}{0.49,0.15,0.80}
\definecolor{purple4}{rgb}{0.33,0.10,0.55}
\definecolor{red}{rgb}{1.00,0.00,0.00}
\definecolor{red1}{rgb}{1.00,0.00,0.00}
\definecolor{red2}{rgb}{0.93,0.00,0.00}
\definecolor{red3}{rgb}{0.80,0.00,0.00}
\definecolor{red4}{rgb}{0.55,0.00,0.00}
\definecolor{rosybrown}{rgb}{0.74,0.56,0.56}
\definecolor{rosybrown1}{rgb}{1.00,0.76,0.76}
\definecolor{rosybrown2}{rgb}{0.93,0.71,0.71}
\definecolor{rosybrown3}{rgb}{0.80,0.61,0.61}
\definecolor{rosybrown4}{rgb}{0.55,0.41,0.41}
\definecolor{royalblue}{rgb}{0.25,0.41,0.88}
\definecolor{royalblue1}{rgb}{0.28,0.46,1.00}
\definecolor{royalblue2}{rgb}{0.26,0.43,0.93}
\definecolor{royalblue3}{rgb}{0.23,0.37,0.80}
\definecolor{royalblue4}{rgb}{0.15,0.25,0.55}
\definecolor{saddlebrown}{rgb}{0.55,0.27,0.07}
\definecolor{salmon}{rgb}{0.98,0.50,0.45}
\definecolor{salmon1}{rgb}{1.00,0.55,0.41}
\definecolor{salmon2}{rgb}{0.93,0.51,0.38}
\definecolor{salmon3}{rgb}{0.80,0.44,0.33}
\definecolor{salmon4}{rgb}{0.55,0.30,0.22}
\definecolor{sandybrown}{rgb}{0.96,0.64,0.38}
\definecolor{seagreen}{rgb}{0.18,0.55,0.34}
\definecolor{seagreen1}{rgb}{0.33,1.00,0.62}
\definecolor{seagreen2}{rgb}{0.31,0.93,0.58}
\definecolor{seagreen3}{rgb}{0.26,0.80,0.50}
\definecolor{seagreen4}{rgb}{0.18,0.55,0.34}
\definecolor{seashell}{rgb}{1.00,0.96,0.93}
\definecolor{seashell1}{rgb}{1.00,0.96,0.93}
\definecolor{seashell2}{rgb}{0.93,0.90,0.87}
\definecolor{seashell3}{rgb}{0.80,0.77,0.75}
\definecolor{seashell4}{rgb}{0.55,0.53,0.51}
\definecolor{sienna}{rgb}{0.63,0.32,0.18}
\definecolor{sienna1}{rgb}{1.00,0.51,0.28}
\definecolor{sienna2}{rgb}{0.93,0.47,0.26}
\definecolor{sienna3}{rgb}{0.80,0.41,0.22}
\definecolor{sienna4}{rgb}{0.55,0.28,0.15}
\definecolor{skyblue1}{rgb}{0.53,0.81,1.00}
\definecolor{skyblue2}{rgb}{0.49,0.75,0.93}
\definecolor{skyblue3}{rgb}{0.42,0.65,0.80}
\definecolor{skyblue4}{rgb}{0.29,0.44,0.55}
\definecolor{skyblue}{rgb}{0.53,0.81,0.92}
\definecolor{slateblue1}{rgb}{0.51,0.44,1.00}
\definecolor{slateblue2}{rgb}{0.48,0.40,0.93}
\definecolor{slateblue3}{rgb}{0.41,0.35,0.80}
\definecolor{slateblue4}{rgb}{0.28,0.24,0.55}
\definecolor{slateblue}{rgb}{0.42,0.35,0.80}
\definecolor{slategray}{rgb}{0.44,0.50,0.56}
\definecolor{slategray1}{rgb}{0.78,0.89,1.00}
\definecolor{slategray2}{rgb}{0.73,0.83,0.93}
\definecolor{slategray3}{rgb}{0.62,0.71,0.80}
\definecolor{slategray4}{rgb}{0.42,0.48,0.55}
\definecolor{snow}{rgb}{1.00,0.98,0.98}
\definecolor{snow1}{rgb}{1.00,0.98,0.98}
\definecolor{snow2}{rgb}{0.93,0.91,0.91}
\definecolor{snow3}{rgb}{0.80,0.79,0.79}
\definecolor{snow4}{rgb}{0.55,0.54,0.54}
\definecolor{springgreen1}{rgb}{0.00,1.00,0.50}
\definecolor{springgreen2}{rgb}{0.00,0.93,0.46}
\definecolor{springgreen3}{rgb}{0.00,0.80,0.40}
\definecolor{springgreen4}{rgb}{0.00,0.55,0.27}
\definecolor{springgreen}{rgb}{0.00,1.00,0.50}
\definecolor{steelblue}{rgb}{0.27,0.51,0.71}
\definecolor{steelblue1}{rgb}{0.39,0.72,1.00}
\definecolor{steelblue2}{rgb}{0.36,0.67,0.93}
\definecolor{steelblue3}{rgb}{0.31,0.58,0.80}
\definecolor{steelblue4}{rgb}{0.21,0.39,0.55}
\definecolor{tan}{rgb}{0.82,0.71,0.55}
\definecolor{tan1}{rgb}{1.00,0.65,0.31}
\definecolor{tan2}{rgb}{0.93,0.60,0.29}
\definecolor{tan3}{rgb}{0.80,0.52,0.25}
\definecolor{tan4}{rgb}{0.55,0.35,0.17}
\definecolor{thistle}{rgb}{0.85,0.75,0.85}
\definecolor{thistle1}{rgb}{1.00,0.88,1.00}
\definecolor{thistle2}{rgb}{0.93,0.82,0.93}
\definecolor{thistle3}{rgb}{0.80,0.71,0.80}
\definecolor{thistle4}{rgb}{0.55,0.48,0.55}
\definecolor{tomato}{rgb}{1.00,0.39,0.28}
\definecolor{tomato1}{rgb}{1.00,0.39,0.28}
\definecolor{tomato2}{rgb}{0.93,0.36,0.26}
\definecolor{tomato3}{rgb}{0.80,0.31,0.22}
\definecolor{tomato4}{rgb}{0.55,0.21,0.15}
\definecolor{turquoise}{rgb}{0.25,0.88,0.82}
\definecolor{turquoise1}{rgb}{0.00,0.96,1.00}
\definecolor{turquoise2}{rgb}{0.00,0.90,0.93}
\definecolor{turquoise3}{rgb}{0.00,0.77,0.80}
\definecolor{turquoise4}{rgb}{0.00,0.53,0.55}
\definecolor{violetred1}{rgb}{1.00,0.24,0.59}
\definecolor{violetred2}{rgb}{0.93,0.23,0.55}
\definecolor{violetred3}{rgb}{0.80,0.20,0.47}
\definecolor{violetred4}{rgb}{0.55,0.13,0.32}
\definecolor{violetred}{rgb}{0.82,0.13,0.56}
\definecolor{violet}{rgb}{0.93,0.51,0.93}
\definecolor{wheat}{rgb}{0.96,0.87,0.70}
\definecolor{wheat1}{rgb}{1.00,0.91,0.73}
\definecolor{wheat2}{rgb}{0.93,0.85,0.68}
\definecolor{wheat3}{rgb}{0.80,0.73,0.59}
\definecolor{wheat4}{rgb}{0.55,0.49,0.40}
\definecolor{white}{rgb}{1.00,1.00,1.00}
\definecolor{whitesmoke}{rgb}{0.96,0.96,0.96}
\definecolor{yellow}{rgb}{1.00,1.00,0.00}
\definecolor{yellow1}{rgb}{1.00,1.00,0.00}
\definecolor{yellow2}{rgb}{0.93,0.93,0.00}
\definecolor{yellow3}{rgb}{0.80,0.80,0.00}
\definecolor{yellow4}{rgb}{0.55,0.55,0.00}
\definecolor{yellowgreen}{rgb}{0.60,0.80,0.20}
\makeatletter
\newlength\@tempdim@x
\newlength\@tempdim@y
\newcommand\AtUpperLeftCorner[3]{%
\begingroup
\@tempdim@x=0cm
\@tempdim@y=\paperheight
\advance\@tempdim@x#1
\advance\@tempdim@y-#2
\put(\LenToUnit{\@tempdim@x},\LenToUnit{\@tempdim@y}){#3}%
\endgroup}
\newcommand\AtUpperRightCorner[3]{%
\begingroup
\@tempdim@x=\paperwidth
\@tempdim@y=\paperheight
\advance\@tempdim@x-#1
\advance\@tempdim@y-#2
\put(\LenToUnit{\@tempdim@x},\LenToUnit{\@tempdim@y}){#3}%
\endgroup}
\newcommand\AtLowerLeftCorner[3]{%
\begingroup
\@tempdim@x=0cm
\@tempdim@y=0cm
\advance\@tempdim@x#1
\advance\@tempdim@y#2
\put(\LenToUnit{\@tempdim@x},\LenToUnit{\@tempdim@y}){#3}%
\endgroup}
\newcommand\AtLowerRightCorner[3]{%
\begingroup
\@tempdim@x=\paperwidth
\@tempdim@y=0cm
\advance\@tempdim@x-#1
\advance\@tempdim@y#2
\put(\LenToUnit{\@tempdim@x},\LenToUnit{\@tempdim@y}){#3}%
\endgroup}
\makeatother
\newtheorem{theoreme}{Theorem}[section]

\newtheorem{definition}[theoreme]{Definition}
\newtheorem{proposition}[theoreme]{Proposition}
\newtheorem{lemme}[theoreme]{Lemma}
\newtheorem{corollaire}[theoreme]{Corollary}

\newtheorem{remarque}[theoreme]{Remark}

\newenvironment{demo}{\begin{proof}}{\end{proof}}

\addto\captionsfrenchb{}
\addto\captionsfrenchb{}

    \newlength{\myarrowsize} 
    \newlength{\myoldlinewidth}

    \pgfarrowsdeclare{myto}{myto}{
        \pgfsetdash{}{0pt} 
        \pgfsetbeveljoin 
        \pgfsetroundcap 
        \setlength{\myarrowsize}{0.5pt}
        \addtolength{\myarrowsize}{.5\pgflinewidth}
        \pgfarrowsleftextend{-4\myarrowsize-.5\pgflinewidth} 
        \pgfarrowsrightextend{.7\pgflinewidth}
    }{
        \setlength{\myarrowsize}{0.4pt} 
        \addtolength{\myarrowsize}{.3\pgflinewidth}  
        \setlength{\myoldlinewidth}{\pgflinewidth}
        \pgfsetroundjoin
        \pgfsetlinewidth{0.0001pt}
        \pgfpathmoveto{\pgfpoint{0.43\myarrowsize}{0}}
        \pgfpatharc{0}{70}{0.14\myarrowsize}
        \pgfpatharc{-80}{-169.5}{4\myarrowsize}
        \pgfpatharc{150}{189}{0.95\myarrowsize and 0.95\myarrowsize}
        \pgfpatharc{0}{-40}{15\myarrowsize}
        \pgfpathmoveto{\pgfpoint{0.43\myarrowsize}{0}}
        \pgfpatharc{0}{-70}{0.14\myarrowsize}
        \pgfpatharc{80}{169.5}{4\myarrowsize}
        \pgfpatharc{-150}{-189}{0.95\myarrowsize and 0.95\myarrowsize}
        \pgfpatharc{0}{40}{15\myarrowsize}
        \pgfpathclose
        \pgfsetstrokeopacity{0.25}
        \pgfusepathqfillstroke
    }

    \pgfarrowsdeclare{myonto}{myonto}{
        \pgfsetdash{}{0pt} 
        \pgfsetbeveljoin 
        \pgfsetroundcap 
        \setlength{\myarrowsize}{0.5pt}
        \addtolength{\myarrowsize}{.5\pgflinewidth}
        \pgfarrowsleftextend{-4\myarrowsize-.5\pgflinewidth} 
        \pgfarrowsrightextend{.7\pgflinewidth}
    }{
        \setlength{\myarrowsize}{0.4pt} 
        \addtolength{\myarrowsize}{.3\pgflinewidth}  
        \setlength{\myoldlinewidth}{\pgflinewidth}
        \pgfsetroundjoin
        \pgfsetlinewidth{0.0001pt}
        \pgfpathmoveto{\pgfpoint{0.43\myarrowsize}{0}}
        \pgfpatharc{0}{70}{0.14\myarrowsize}
        \pgfpatharc{-80}{-169.5}{4\myarrowsize}
        \pgfpatharc{150}{189}{0.95\myarrowsize and 0.95\myarrowsize}
        \pgfpatharc{0}{-40}{15\myarrowsize}
        \pgfpathmoveto{\pgfpoint{-7\myarrowsize}{0}}
				\pgfpatharc{0}{70}{0.14\myarrowsize}
        \pgfpatharc{-80}{-169.5}{4\myarrowsize}
        \pgfpatharc{150}{189}{0.95\myarrowsize and 0.95\myarrowsize}
        \pgfpatharc{0}{-40}{15\myarrowsize}
        \pgfpathmoveto{\pgfpoint{0.43\myarrowsize}{0}}
        \pgfpatharc{0}{-70}{0.14\myarrowsize}
        \pgfpatharc{80}{169.5}{4\myarrowsize}
        \pgfpatharc{-150}{-189}{0.95\myarrowsize and 0.95\myarrowsize}
        \pgfpatharc{0}{40}{15\myarrowsize}
			  \pgfpathmoveto{\pgfpoint{-7\myarrowsize}{0}}
        \pgfpatharc{0}{-70}{0.14\myarrowsize}
        \pgfpatharc{80}{169.5}{4\myarrowsize}
        \pgfpatharc{-150}{-189}{0.95\myarrowsize and 0.95\myarrowsize}
        \pgfpatharc{0}{40}{15\myarrowsize}
        \pgfpathclose
        \pgfsetstrokeopacity{0.25}
        \pgfusepathqfillstroke
        \pgfusepathqfillstroke
    }

    \pgfarrowsdeclare{myhook}{myhook}{
        \setlength{\myarrowsize}{0.6pt}
        \addtolength{\myarrowsize}{.5\pgflinewidth}
        \pgfarrowsleftextend{-4\myarrowsize-.5\pgflinewidth} 
        \pgfarrowsrightextend{.7\pgflinewidth}
    }{
        \setlength{\myarrowsize}{0.6pt} 
        \addtolength{\myarrowsize}{.5\pgflinewidth}  
        \pgfsetdash{}{+0pt}
        \pgfsetroundcap
        \pgfpathmoveto{\pgfqpoint{-2pt}{-6\pgflinewidth}}
        \pgfpathcurveto
            {\pgfqpoint{4\pgflinewidth}{-4.667\pgflinewidth}}
            {\pgfqpoint{4\pgflinewidth}{0pt}}
            {\pgfpointorigin}
        \pgfusepathqstroke
    }

		\pgfarrowsdeclare{my to}{my to}
{
  \pgfarrowsleftextend{-2\pgflinewidth}
  \pgfarrowsrightextend{\pgflinewidth}
}
{
  \pgfsetlinewidth{0.8\pgflinewidth}
  \pgfsetdash{}{0pt}
  \pgfsetroundcap
  \pgfsetroundjoin
  \pgfpathmoveto{\pgfpoint{-5.5\pgflinewidth}{7.5\pgflinewidth}}
  \pgfpathcurveto
  {\pgfpoint{-4.0\pgflinewidth}{0.1\pgflinewidth}}
  {\pgfpoint{0pt}{0.25\pgflinewidth}}
  {\pgfpoint{0.75\pgflinewidth}{0pt}}
  \pgfpathcurveto
  {\pgfpoint{0pt}{-0.25\pgflinewidth}}
  {\pgfpoint{-4.0\pgflinewidth}{-0.1\pgflinewidth}}
  {\pgfpoint{-5.5\pgflinewidth}{-7.5\pgflinewidth}}
  \pgfusepathqstroke
}

\tikzstyle{vecArrow} = [thick, decoration={markings,mark=at position
   1 with {\arrow[semithick]{open triangle 60}}},
   double distance=1.4pt, shorten >= 5.5pt,
   preaction = {decorate},
   postaction = {draw,line width=1.4pt, white,shorten >= 4.5pt}]
\tikzstyle{innerWhite} = [semithick, white,line width=1.4pt, shorten >= 4.5pt]

	\makeatletter
	\newcommand\POSITION[3]{%
	\begingroup
	\@tempdim@x=0cm
	\@tempdim@y=\paperheight
	\advance\@tempdim@x#1
	\advance\@tempdim@y-#2
	\put(\LenToUnit{\@tempdim@x},\LenToUnit{\@tempdim@y}){#3}%
	\endgroup
	}
\makeatother

\addto\captionsenglish{}

\begin{document}	
	
	\begin{abstract}
		We   view the $\tilde{A}$-type affine braid group  as a subgroup of the $B$-type braid group. We  show that   the $\tilde{A}$-type affine braid group   surjects onto    the $A$-type braid group    and we detect the kernel of this surjection using Schreier's Theorem.   We then describe an injection of  the $B$-type braid group    into 
		 the $A$-type braid group   which allows us 
		finally to give a definition of affine links, as  closures of affine braids viewed as A-type braids after composing the above injections, and we prove that the two conditions of Markov are necessary and sufficient to get the same affine closure of any two affine braids.		
		 \end{abstract}

		\maketitle

	\section{Introduction}

 In \cite{Sadek_2013_1} we define a tower of affine Temperley-Lieb algebras of 
 type $\tilde A$, we define Markov conditions for traces relative to this tower, hence the affine notion of a  Markov trace,  and we show 
 that there exists a unique Markov trace on the above tower. These results are obtained in a purely algebraic way, independently of their translation in topological terms, that is, in terms of invariants of links, however important. They were the goal of my work \cite{Sadek_Thesis}, in which the consequences in terms of invariants of links were explained. We 
describe below the braid groups of type $A$, $B$ and $\tilde A$ involved, both in an algebraic and a geometric   way, we propose a definition of an affine link  as the closure of an affine braid viewed as an $A$-type braid under an explicit injection, 
and we show that the two conditions of Markov are necessary and sufficient to get the same affine closure of any two affine braids.

\section[Artin groups and Braid groups]{Artin groups and Braid groups} \label{1_1}

		Definitions and results of paragraphs \ref{1_1} and \ref{1_2} are taken mostly from \cite{Paris_2009}.  \\
		
		 	Let $S$ be a finite set. 
		
		\begin{definition}
			A Coxeter matrix over $S$ is a square matrix $M = (m_{st})_{s,t \in S}$ such that \\
		
			\begin{itemize}[label=$\bullet$, font=\normalsize, font=\color{black}, leftmargin=2cm,parsep=0cm, itemsep=0.25cm, topsep=0cm]
				\item $m_{ss} = 1$,			
				\item $m_{st} = m_{ts} $ for any $s,t$ in $S$,				
				\item $m_{st} $ belongs to $ \left\{ 2,3,4 ... \infty \right\}$.\\			
			\end{itemize}
		\end{definition}	
		
		We present a Coxeter matrix by its Dynkin graph $\Gamma = \Gamma(M)$, which is a graph given by vertices and edges. $\Gamma$ has $S$ as a set of vertices, and for any non-equal two vertices $s,t$ in $S$ we have \\ 
				
		\begin{itemize}[label=$\bullet$, font=\normalsize, font=\color{black}, leftmargin=2cm,parsep=0cm, itemsep=0.25cm,topsep=0cm]
			\item $s,t$ are joined by an edge if  $ m_{s,t} = 3 $, 			
			\item $s,t$ are joined by a doubled edge if  $ m_{s,t} = 4 $,			
			\item $s,t$ are joined by an edge  labeled by $m_{s,t}$ if  $ m_{s,t} > 4  $.\\		
		\end{itemize}

		
		\begin{figure}[ht]
			\centering
		
			\begin{tikzpicture}
			 \begin{scope}[xscale = 1.5]
  \draw (-0.7, 0) -- (0.7, 0);
  \filldraw (-0.7,0) circle (2pt);
  \filldraw (0.7,0) circle (2pt);
  \node at (-0.7,-0.5) {$s$};
  \node at (0.7,-0.5) {$t$};
    \end{scope}
 \begin{scope}[xshift =3cm, xscale= 1.5]
   \draw (-0.7, 0.07) -- (0.7, 0.07);
  \draw (-0.7, -0.07) -- (0.7, -0.07);
  \filldraw (-0.7,0) circle (2pt);
  \filldraw (0.7,0) circle (2pt);
  \node at (-0.7,-0.5) {$s$};
  \node at (0.7,-0.5) {$t$};
   
 \end{scope}
 \begin{scope}[xshift =6cm, xscale= 1.5]
 
  \draw (-0.7, 0.0) -- (0.7, 0.0);
  \filldraw (-0.7,0) circle (2pt);
  \filldraw (0.7,0) circle (2pt);
  \node at (-0.7,-0.5) {$s$};
  \node at (0.7,-0.5) {$t$};
  \node at (0, 0.2) {$m_{s,t}$};
 \end{scope}
 \end{tikzpicture}
			\caption{Edges}
		\end{figure}

		\begin{definition}  \label{1_1_2}
			Let $B_{S}$ be the set $\left\{ \sigma _{s};~ s \in S\right\}$.We call the pair $(B,S )$ an Artin system of type $\Gamma$, where $B=B_{\Gamma}$ is the group given by generators and relations as follows: $S$ is the set of generators with relations $ prod ( \sigma _{s},\sigma _{t} : m_{st}) = prod ( \sigma _{t} \sigma _{s} : m_{st})$, for any non-equal $ s,t$ in $S$ with $m_{s,t} \neq \infty $. \\
		
			We call $ B$ the Artin group of type $ \Gamma$.\\ 
		\end{definition}
		
		In our work we treat many kinds of Artin groups in which $ m_{s,t} \leq 4 $. The   relations appearing in the definition are called "braid relations".\\
		
		Let $\Gamma$ be a Dynkin graph. Let  $(B,S )$ be the related Artin system. Take $V$ to be the real vector space with $\mathds{R}$-basis $\left\{ e_{s}; s \in S\right\}$ which plays the role of the set of simple roots. The root system gives rise to simple reflections hence to a reflection group generated by those simple reflections, say $ W_{S}$. By the natural linear representation of $ W_{S}$ we can realize it as a subgroup in $GL(V)$ (the group of endomorphisms of $V$). Let $R$ be the set of reflections of $ W_{S}$ (the set of conjugates of simple reflections). Take $r$ in $R$, since it is a reflection it fixes a hyper-plane in $V$, say $ H_{r} $. In fact $ W_{S}$ acts freely on the complement of $\cup_{r\in R}H_{r}$ in $V$. By extending the action of $ W_{S}$ up to $V_{\mathds{C}} = \mathds{C} \otimes _{\mathds{R}} V $ we see that $ W_{S}$ acts freely  on the complement of $\cup_{r\in R}\mathds{C} \otimes _{\mathds{R}}H_{r}$ in $V_{\mathds{C}}$. We call this complement $ M_{\Gamma}$. We set $N_{\Gamma} = M_{\Gamma}/ W_{S} $.\\
		
		\begin{remarque}
			Before stating the theorem, we have to notice that the above argument is valid in the case of finite $ W_{S}$. For when the group $ W_{S}$ is infinite we have to replace $V$ by $ U \subset V $ (the so-called Tits cone), and to replace as well $V_{\mathds{C}}$ by $( U+iV ) \subset V_{\mathds{C}}$.  $M_{\Gamma }$ is to be $ ( U+iV )- \cup_{r\in R}\mathds{C} \otimes _{\mathds{R}}H_{r}$. The action of $ W_{S}$ on $M_{\Gamma }$ is free, and as above $N_{\Gamma}$ is  $M_{\Gamma }$  modulo the action of  $ W_{S}$. \\
		\end {remarque}
		
		\begin{definition}
			The braid group of $\Gamma$-type is $\pi_{1}(N_{\Gamma})$, the fundamental group of the space $N_{\Gamma}$.\\
		\end{definition}
		
		\begin{theoreme}
			(Brieskorn-Van der Lek). $\pi_{1}(N_{\Gamma}) \simeq  B $.\\		
		\end{theoreme}
		
		Briefly: the $\Gamma$-type Artin group is given by a presentation, while the  $\Gamma$-type braid group is a fundamental group. Although the Brieskorn-Van der Lek isomorphism is not canonical, we will not make the distinction in this work. We will call each of these  groups a $\Gamma$-type braid group.

	\section[$A$-type braid groups]{$A$-type braid groups} \label{1_2}
		
		An $A$-type braid group with $n$ generators is historically the first braid group. We give its presentation by generators and relations, then a geometrical one (by means of braids with $n+1$ strands). Many interesting basic facts show the reasons for which it has such respectable position in the group theory, in addition to many other branches of mathematics:  for example it has a faithful representation  in $Aut(F_{n+1})$, the group of automorphisms of the free group with $n$ generators;  it has  strong relations with link theory (here comes the well known Alexander theorem); in addition of its Garsiditude, in fact, it is the first group to be called a Garside group.\\ 
		
		\subsection{Presentations}
		
			Let $n \geq 1$ be an integer.\\
		
			\begin{definition}
				The $A$-type braid group $B(A_{n})$ with $n$ generators is the group presented by a generator set $S=\left\{\sigma_{1}, \sigma_{2}, \dots ,\sigma_{n}\right\}$ and the relations\\
		
				\begin{itemize}[label=$\bullet$, font=\normalsize, font=\color{black}, leftmargin=2cm,parsep=0cm, itemsep=0.25cm,topsep=0cm]
					\item $\sigma_{i} \sigma_{j} =\sigma_{j} \sigma_{i} $  where $1\leq i,j\leq n$ and $ \left| i-j\right| \geq 2$,
					\item $\sigma_{i}\sigma_{i+1}\sigma_{i} = \sigma_{i+1}\sigma_{i}\sigma_{i+1}$ where $1\leq i\leq n-1$.\\ 
				\end{itemize}
			\end{definition}
		
			Thus the related Dynkin diagram is	\\
				\vspace{0.5cm}	
			 			
					\begin{figure}[ht]
				\centering
				\begin{tikzpicture}

  \filldraw (0,0) circle (2pt);
  \node at (0,-0.5) {$\sigma_{1}$}; 
   
  \draw (0,0) -- (1.5, 0);

  \filldraw (1.5,0) circle (2pt);
  \node at (1.5,-0.5) {$\sigma_{2}$};

  \draw (1.5,0) -- (3, 0);

  \node at (3.5,0) {$\dots$};

  \draw (4,0) -- (5.5, 0);
  
  \filldraw (5.5,0) circle (2pt);
  \node at (5.5,-0.5) {$\sigma_{n-1}$};
 
  \draw (5.5,0) -- (7, 0);
  
  \filldraw (7,0) circle (2pt);
  \node at (7,-0.5) {$\sigma_{n}$};

               \end{tikzpicture}
			 \caption{$\Gamma_{A}$}
			\end{figure}

		\vspace{0.5cm}		
			 Let $P_{n}, \dots, P_{n+1} $ be distinct points in the plane $\mathds{R}^{2}$. We  define a rough braid on $n+1$ strands to be an $n$-tuple $\beta = (b_{1}, \dots , b_{n+1})$, where  $b_{k}$  is a path  $b_{k} : \left[0,1\right] \rightarrow  \mathds{R}$ such that\\ 
	
			\begin{itemize}[label=$\bullet$, font=\normalsize, font=\color{black}, leftmargin=2cm, parsep=0cm, itemsep=0.25cm,topsep=0cm]
				\item For any $k$ in  $ \left\{ 1, \dots, n+1\right\} $ we have $b_{k}(0) = P_{k}$,
				\item For any $k$ in  $ \left\{ 1, \dots, n+1\right\} $ there exists a permutation $ x = \theta (\beta) \in Sym_{n+1}$ such that  $b_{k}(1) = P_{x(k)} $,
				\item For any non-equal $k$ and $l$ in  $ \left\{ 1, \dots, n+1\right\} $, for all $t\in \left[0,1\right] $ we have $ b_{k}(t) \neq b_{l}(t) $. \\
			\end{itemize}
	
			By definition: two rough braids $\alpha$ and $\beta $ are homotopic  if there exists a continuous family of rough braids $\left\{\gamma_{s}\right\}_{s\in \left[0,1\right]}$ such that $ \gamma_{0} = \alpha $ and $ \gamma_{1} = \beta $ . This is an equivalence relation. \\
	
			\begin{definition}
				A braid on $n+1$ strands  is a homotopy class of rough braids on $n+1$ strands.
			\end{definition}
	
	\vspace{0.5cm}
	
			The well known geometric interpretation of the elements of $B(A_{n})$ viewed as braids in the space  is the following \\ 
	
			 \vspace{0.5cm}

			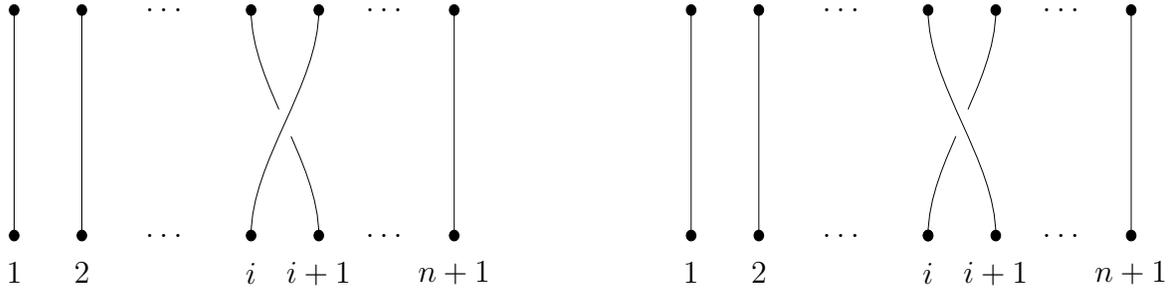
\begin{figure}[ht]
				\centering
				%
\begin{tikzpicture}
\begin{scope}[xscale = 0.9]
  \draw[white,line width = 2pt] (-1,1)-- +(0,-3);

  \filldraw (0,1) circle (2pt);  
  \draw[line width = 0.3pt] (0,1)-- +(0,-3);
  \filldraw (0,-2) circle (2pt);
  \node at (0, -2.5) {$1$};

  \filldraw (1,1) circle (2pt);  
  \draw[line width = 0.3pt] (1,1)-- +(0,-3);
  \filldraw (1,-2) circle (2pt);
  \node at (1, -2.5) {$2$};

	\node at (2.25,1) {$\dots$};
	\node at (2.25,-2) {$\dots$};
	
	\filldraw (3.5,1) circle (2pt); 
	\draw (3.5,1) [line width = 0.3pt] .. controls +(0,-1)  and +(0,1) ..  ++(1,-3);
	\filldraw (4.5,-2) circle (2pt); 
	\node at (3.5, -2.5) {$i$};
	
	\filldraw (4.5,1) circle (2pt); 	
	\fill[white] (4,-0.5) circle (0.2cm);
	\draw (4.5,1) [line width = 0.3pt] .. controls +(0,-1)  and +(0,1) ..  +(-1,-3);
	\filldraw (3.5,-2) circle (2pt); 
	\node at (4.5, -2.5) {$i+1$};
	
	\node at (5.5,1) {$\dots$};
	\node at (5.5,-2) {$\dots$};
		
  \filldraw (6.5,1) circle (2pt);  
  \draw[line width = 0.3pt] (6.5,1)-- +(0,-3);
  \filldraw (6.5,-2) circle (2pt);
  \node at (6.5, -2.5) {$n+1$};


  \draw[white,line width = 2pt] (9,1)-- +(0,-3);

  \filldraw (10,1) circle (2pt); 
  \draw[line width = 0.3pt] (10,1)-- +(0,-3);
	\filldraw (10,-2) circle (2pt); 
  \node at (10, -2.5) {$1$};
	
	\filldraw (11,1) circle (2pt); 
  \draw[line width = 0.3pt] (11,1)-- +(0,-3);
	\filldraw (11,-2) circle (2pt);
  \node at (11, -2.5) {$2$};

	\node at (12.25,1) {$\dots$};
	\node at (12.25,-2) {$\dots$};
	
	\draw (14.5,1) [line width = 0.3pt] .. controls +(0,-1)  and +(0,1) ..  +(-1,-3);
	\filldraw (13.5,1) circle (2pt);
	
	\filldraw (14.5,1) circle (2pt);
	\node at (14.5, -2.5) {$i+1$};
	\filldraw (14.5,-2) circle (2pt);
	\fill[white] (14,-0.5) circle (0.2cm);	
	
	\filldraw (13.5,-2) circle (2pt);
	\draw (13.5,1) [line width = 0.3pt] .. controls +(0,-1)  and +(0,1) ..  ++(1,-3);
	\filldraw (13.5,-2) circle (2pt);
	\node at (13.5, -2.5) {$i$};
	
	\node at (15.5,1) {$\dots$};
	\node at (15.5,-2) {$\dots$};
		
	\filldraw (16.5,1) circle (2pt);	
  \draw[line width = 0.3pt] (16.5,1)-- +(0,-3);
	\filldraw (16.5,-2) circle (2pt);
  \node at (16.5, -2.5) {$n+1$};

 \end{scope}
        \end{tikzpicture}

				\caption{$\sigma_{i} ~~~~~~~~ \& ~~~~~~~~ \sigma^{-1}_{i}$}
			\end{figure}
			 
			 \vspace{0.5cm}
			 			
			\begin{figure}[ht]
				\centering
				%
				\begin{tikzpicture}
				\begin{scope}[xscale = 1.5]

	\filldraw (1,1) circle (2pt);
	\filldraw (2,1) circle (2pt);
	\filldraw (5,1) circle (2pt);
	
	\filldraw (1,-1) circle (2pt);
	\filldraw (2,-1) circle (2pt);
	\filldraw (5,-1) circle (2pt);

	\draw (1,1) -- +(0,-2);
	\draw (2,1) -- +(0,-2);
	\draw (5,1) -- +(0,-2);
	\node at (3.5,1) {$\dots$};
	\node at (3.5,-1) {$\dots$};
	\node at (1,-1.5) {$1$};
	\node at (2,-1.5) {$2$};
	hgcjh
\node at (5,-1.5) {$n+1$};
 \end{scope}
      \end{tikzpicture}

				\caption{Id}
			\end{figure}
			 
			 \vspace{0.5cm}	
			
			We compose two braids in  the way that one would expect, that is for any two braids $X,Y$ the composed braid $XY$ is the braid obtained by putting $X$ at the top and $Y$ at the bottom, welding the bottom end points of $X$ with the upper ones of $Y$ (the $i$-th with the $i$-th, $1\leq i \leq n+1 $) as follows:  

			\vspace{2.5cm}

			\begin{figure}[ht]
				\centering
	\begin{tikzpicture}
\begin{scope}[xscale = 1.25]

	\filldraw (1,1) circle (2pt);
	\filldraw (2,1) circle (2pt);
	\filldraw (5,1) circle (2pt);
	
  \draw[line width = 0.3pt] (1,1) -- +(0,-2);
  \draw[line width = 0.3pt] (2,1) -- +(0,-2);
  \draw[line width = 0.3pt] (5,1) -- +(0,-2);
	
	\filldraw (1,-1) circle (2pt);
	\filldraw (2,-1) circle (2pt);
	\filldraw (5,-1) circle (2pt);

  \filldraw[fill = white, draw = black] (0.7, 0.5) -- (5.3, 0.5) -- (5.3,-0.5) -- (0.7, -0.5) -- cycle;
  \node at (3, 0) {$x$};
	\node at (3.5,1) {$\dots$};
	\node at (3.5,-1) {$\dots$};
	
	\node at (1, -1.5) {$1$};
	\node at (2, -1.5) {$2$};
	\node at (5, -1.5) {$n+1$};
\end{scope}

\begin{scope}[xscale = 1.25, yshift = - 3cm]

	\filldraw (1,1) circle (2pt);
	\filldraw (2,1) circle (2pt);
	\filldraw (5,1) circle (2pt);

  \draw[line width = 0.3pt] (1,1) -- +(0,-2);
  \draw[line width = 0.3pt] (2,1) -- +(0,-2);
  \draw[line width = 0.3pt] (5,1) -- +(0,-2);
	
	\filldraw (1,-1) circle (2pt);
	\filldraw (2,-1) circle (2pt);
	\filldraw (5,-1) circle (2pt);
	
  \filldraw[fill = white, draw = black] (0.7, 0.5) -- (5.3, 0.5) -- (5.3,-0.5) -- (0.7, -0.5) -- cycle;
  \node at (3, 0) {$y$};

	\node at (3.5,1) {$\dots$};
	\node at (3.5,-1) {$\dots$};

	\node at (1, -1.5) {$1$};
	\node at (2, -1.5) {$2$};
	\node at (5, -1.5) {$n+1$};
\end{scope}

\node at (9, -1.5) {$\longrightarrow$};

\begin{scope}[xscale = 1.25, yshift = - 1.5cm, xshift = 8cm]

	\filldraw (1,1) circle (2pt);
	\filldraw (2,1) circle (2pt);
	\filldraw (5,1) circle (2pt);

  \draw[line width = 0.3pt] (1,2) -- +(0,-4);
  \draw[line width = 0.3pt] (2,2) -- +(0,-4);
  \draw[line width = 0.3pt] (5,2) -- +(0,-4);
	
	\filldraw (1,-1) circle (2pt);
	\filldraw (2,-1) circle (2pt);
	\filldraw (5,-1) circle (2pt);
	
	\filldraw (1,2) circle (2pt);
	\filldraw (2,2) circle (2pt);
	\filldraw (5,2) circle (2pt);
		
	\filldraw[fill = white, draw = black] (0.7, 1.5) -- (5.3, 1.5) -- (5.3,0.5) -- (0.7, 0.5) -- cycle;
  \node at (3, 1) {$x$};
  \filldraw[fill = white, draw = black] (0.7, -1.5) -- (5.3, -1.5) -- (5.3,-0.5) -- (0.7, -0.5) -- cycle;
  \node at (3, -1) {$y$};
	\node at (3.5,2) {$\dots$};
	\node at (3.5,-2) {$\dots$};
	\node at (3.5,0) {$\dots$};

	\filldraw (1,-2) circle (2pt);
	\filldraw (2,-2) circle (2pt);
	\filldraw (5,-2) circle (2pt);
	
	\node at (1, -2.5) {$1$};
	\node at (2, -2.5) {$2$};
	\node at (5, -2.5) {$n+1$};
\end{scope}
\end{tikzpicture}
 
		 				\caption{$X,Y \rightarrow XY$}
			\end{figure}
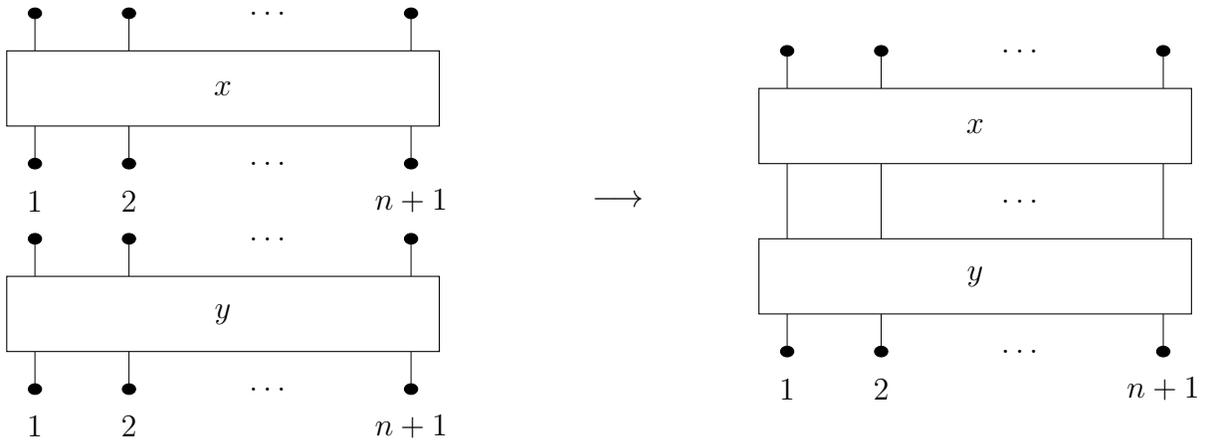
	
	\vspace{1cm}
	
			The natural embedding \\
			
			\begin{eqnarray}
				x_{n-1} : B(A_{n-1}) &\longrightarrow& B(A_{n}) \nonumber\\
				\sigma_{i} &\longmapsto& \sigma_{i} \text{ for } 1\leq i\leq n-1 \nonumber,
			\end{eqnarray}
			
				\vspace{0.5cm}
				
			can be realized geometrically by adding the $(n+1)$-th strand\\
	
	\vspace{1cm}
			
			\begin{figure}[ht]
				\centering
			\begin{tikzpicture}
			\begin{scope}[xscale = 0.85]

				\draw[draw= white, line width = 2pt] (-1,4)-- (-1,-1);

				\draw[draw= white, double = white, line width = 0.1cm] (-0.5, 1.2) .. controls +(-0.3, 0.3) and +(0, -0.3) .. (-1.4, 1.5) .. controls +(0,0.3) and +(-0.3,-0.3) .. (-0.5,1.8); 

				\filldraw[white] (-1,1.35) circle (4pt);
				\draw[draw= white, line width = 2pt] (-1,1.5)-- (-1,0);


			 \filldraw (0,4) circle (2pt);
				\draw[line width = 0.3pt] (0,3) -- (0,4);

				\filldraw (1,4) circle (2pt);
				\draw[line width = 0.3pt] (1,3) -- (1,4);

				\filldraw (5,4) circle (2pt);
				\draw[line width = 0.3pt] (5,3) -- (5,4);

				\draw[line width = 0.3pt] (-0.5,3) -- (5.5,3);

				\draw[line width = 0.3pt] (5.5,3) -- (5.5,0);

				\draw[line width = 0.3pt] (5.5,0) -- (-0.5,0);
				
				\draw[line width = 0.3pt] (-0.5,0) -- (-0.5,3);  
				
				\draw[line width = 0.3pt] (0,0) -- (0,-1);
				\filldraw (0,-1) circle (2pt);
				\node at (0, -1.5) {$1$};

				\draw[line width = 0.3pt] (1,0) -- (1,-1);
				\filldraw (1,-1) circle (2pt);
				\node at (1, -1.5) {2};

				\draw[line width = 0.3pt] (5,0) -- (5,-1);
				\filldraw (5,-1) circle (2pt);
				\node at (5, -1.5) {n};

				\draw[line width =1.5pt][->] (6.25,1.5) -- (7.5,1.5);

\end{scope}

\begin{scope}[xscale = 0.85, xshift = 10cm]

				\draw[draw= white, line width = 2pt] (-1,4)-- (-1,-1);

				\draw[draw= white, double = white, line width = 0.1cm] (-0.5, 1.2) .. controls +(-0.3, 0.3) and +(0, -0.3) .. (-1.4, 1.5) .. controls +(0,0.3) and +(-0.3,-0.3) .. (-0.5,1.8); 

				\filldraw[white] (-1,1.35) circle (4pt);
				\draw[draw= white, line width = 2pt] (-1,1.5)-- (-1,0);


			 \filldraw (0,4) circle (2pt);
				\draw[line width = 0.3pt] (0,3) -- (0,4);

				\filldraw (1,4) circle (2pt);
				\draw[line width = 0.3pt] (1,3) -- (1,4);

				\filldraw (5,4) circle (2pt);
				\draw[line width = 0.3pt] (5,3) -- (5,4);

				\draw[line width = 0.3pt] (-0.5,3) -- (5.5,3);

				\draw[line width = 0.3pt] (5.5,3) -- (5.5,0);

				\draw[line width = 0.3pt] (5.5,0) -- (-0.5,0);
				
				\draw[line width = 0.3pt] (-0.5,0) -- (-0.5,3);  
				
				\draw[line width = 0.3pt] (0,0) -- (0,-1);
				\filldraw (0,-1) circle (2pt);
				\node at (0, -1.5) {$1$};

				\draw[line width = 0.3pt] (1,0) -- (1,-1);
				\filldraw (1,-1) circle (2pt);
				\node at (1, -1.5) {2};

				\draw[line width = 0.3pt] (5,0) -- (5,-1);
				\filldraw (5,-1) circle (2pt);
				\node at (5, -1.5) {n};

				\filldraw (6,4) circle (2pt);
				\draw[line width = 0.3pt] (6,4) -- (6,-1);
				\filldraw (6,-1) circle (2pt);
				\node at (6, -1.5) {n+1};

\end{scope}

			\end{tikzpicture} 
				\caption{$x_{n-1}$}
			\end{figure}

			 \clearpage
			
		\subsection{A faithful representation}
		
		We can realize $B(A_{n})$ as a subgroup of the group of automorphisms of the free group with $n+1$ generators via the "Artin representation", which will be briefly defined in what follows: Let $F_{n+1}$ be the free group with $n+1$ generators $x_{1}, .., x_{n+1}$. Let $Aut(F_{n+1})$ be the group of automorphisms of $F_{n+1}$. For $1\leq k \leq n$, we define $t_{k} $ in $Aut(F_{n+1})$ as follows for $i \neq k,k+1:$ 		
			\begin{eqnarray}
				t_{k}: F_{n}&\longrightarrow& F_{n} \nonumber\\
				x_{i} &\longmapsto& x_{i}, \nonumber\\
				x_{k} &\longmapsto& x^{-1}_{n} x_{k+1}x_{k}, \nonumber\\
				x_{k+1} &\longmapsto& x_{k}. \nonumber
			\end{eqnarray}

			It is easy to show that the map $ \rho: B(A_{n}) \longrightarrow Aut(F_{n+1}) $, which sends $\sigma_{k}$ to $t_{k}$, defines a representation of $ B(A_{n}) $ in $Aut(F_{n+1}) $ called the Artin representation.\\
			
			\begin{theoreme}  
				(Artin) The representation $\rho$ is faithful.
			\end{theoreme}		
			
	\section[$B$-type braid groups]{$B$-type braid groups} \label{1_3}
						
		The $B$-type braid group with $n+1$ generators  $B(B_{n+1})$ plays a role in the theory of low dimension topological spaces, in addition to the fact that it is very useful in investigating the structure of the affine braid group which is the center of interest of this work in general.  Definitions and results in this paragraph  are taken from \cite{Graham_Lehrer_2003}, where we can see more details about the group $B(B_{n+1})$.   
				
		\subsection{Presentations}	
			
			\vspace{0.5cm}				
			\begin{definition}
				The $B$-type braid group with $n+1$ generators $B(B_{n+1})$ is the group presented by a generators set $ \left\{ \sigma_{1}, \sigma_{2}, \dots,  \sigma_{n}, t   \right\}$ and the relations\\
		
				\begin{itemize}[label=$\bullet$, font=\normalsize, font=\color{black}, leftmargin=2cm, parsep=0cm, itemsep=0.25cm, topsep=0cm]
					\item $\sigma_{i} \sigma_{j} =\sigma_{j} \sigma_{i} $  where $1\leq i,j\leq n$ when $ \left| i-j\right| \geq 2$,
					\item $\sigma_{i}\sigma_{i+1}\sigma_{i} = \sigma_{i+1}\sigma_{i}\sigma_{i+1}$ when $1\leq i\leq n$,
					\item $\sigma_{i} t = t \sigma_{i} $ when $2\leq i \leq n$, 
					\item $\sigma_{1} t \sigma_{1} t = t \sigma_{1} t \sigma_{1} $.\\
				\end{itemize}	
		     \end{definition}
		
		The related Dynkin diagram is \\
		
		
		\begin{figure}[ht]
			\centering
				\begin{tikzpicture}
				\begin{scope}[xscale = 1]

 \filldraw (-1.5,0) circle (2pt);
  \node at (-1.5,-0.5) {$t$}; 

  \draw (-1.5,-0.07) -- (0, -0.07);
  \draw (-1.5,0.07) -- (0, 0.07);
  
  \filldraw (0,0) circle (2pt);
  \node at (0,-0.5) {$\sigma_{1}$}; 
   
  \draw (0,0) -- (1.5, 0);

  \filldraw (1.5,0) circle (2pt);
  \node at (1.5,-0.5) {$\sigma_{2}$};

  \draw (1.5,0) -- (3, 0);

  \node at (3.5,0) {$\dots$};

  \draw (4,0) -- (5.5, 0);
  
  \filldraw (5.5,0) circle (2pt);
  \node at (5.5,-0.5) {$\sigma_{n-1}$};
 
  \draw (5.5,0) -- (7, 0);
  
  \filldraw (7,0) circle (2pt);
  \node at (7,-0.5) {$\sigma_{n}$};

\end{scope}
	\end{tikzpicture}
			\caption{$\Gamma_{B_{n+1}}$}
		\end{figure}
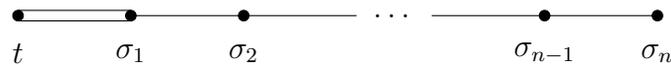
		
		Set $\phi_{n+1} = t\sigma_{1} \dots  \sigma_{n}$. Set $a_{n+1} = \phi_{n+1} \sigma_{n} \phi^{-1}_{n+1} $. We can see directly that $\phi_{n+1} \sigma_{i} \phi^{-1}_{n+1} = \sigma_{i+1}$ for all $ 1\leq i\leq n-1 $, with  $\phi_{n+1} \sigma_{n} \phi^{-1}_{n+1} = a_{n+1}$. What is more, we have the following presentation:\\
		
		\begin{proposition} \label{1_3_2}
			$B(B_{n+1})$ is presented by the set of generators $S'= \left\{ \sigma_{1},  \dots,  \sigma_{n}, a_{n+1}, \phi_{n+1}  \right\}$ together with the relations\\

			\begin{itemize}[label=$\bullet$, font=\normalsize, font=\color{black}, leftmargin=2cm, parsep=0cm, itemsep=0.25cm, topsep=0cm]
				\item $\sigma_{i} \sigma_{j} =\sigma_{j} \sigma_{i} $  where $1\leq i,j\leq n$ when $ \left| i-j\right| \geq 2$,
				\item $\sigma_{i}\sigma_{i+1}\sigma_{i} = \sigma_{i+1}\sigma_{i}\sigma_{i+1}$ when $1\leq i\leq n-1$,
				\item $\sigma_{i} a_{n+1} = a_{n+1} \sigma_{i} $ when $2\leq i \leq n-1$,
				\item $\sigma_{1} a_{n+1} \sigma_{1} = a_{n+1} \sigma_{1} a_{n+1}$,
				\item $\sigma_{n} a_{n+1} \sigma_{n} = a_{n+1} \sigma_{n} a_{n+1}$, 
				\item $\phi_{n+1} \sigma_{i} \phi^{-1}_{n+1} = \sigma_{i+1}$ when $ 1\leq i\leq n-1 $,
				\item $\phi_{n+1} \sigma_{n} \phi^{-1}_{n+1} = a_{n+1}$,
				\item $\phi_{n+1}  a_{n+1} \phi^{-1}_{n+1} = \sigma_{1}$.\\
			\end{itemize}
		\end{proposition} 

		Notice that the relations involving $ \phi_{n+1}$ are not braid relations, i.e., $(B(B_{n+1}),S')$ cannot be viewed as an Artin system.\\
  
  \vspace{0.5cm}
  
		Considering those relations we see directly that $\phi_{n+1}$ is acting by automorphisms on the normal subgroup of $B(B_{n+1})$ generated by $ \left\{ \sigma_{1}, \dots,  \sigma_{n}, a_{n+1} \right\}$. More, this set $ \left\{ \sigma_{1}, \dots,  \sigma_{n}, a_{n+1} \right\}$ with the first five systems of relations forms a presentation by generators and relations of this subgroup which will be the subject of the next section. The element $\phi_{n+1}$ generates a free subgroup of rank 1, which we denote by  $\Phi_{n+1}$. In other terms:	\\ 
			
		\begin{eqnarray}
			B(B_{n+1}) = \left\langle \sigma_{1}, \dots,  \sigma_{n}, a_{n+1} \right\rangle_{B(B_{n+1})} \rtimes \Phi_{n+1}. \nonumber\\\nonumber
		\end{eqnarray}
		
			\vspace{0.5cm}
			
		The geometric presentation of $B(B_{n+1})$ as braids is given in a very similar way of that of $B(A_{n})$ but in $n+2$ strands, where the first strand remains point-wise fixed. The generators are presented as follows: \\
	
		 \clearpage
		
		\begin{figure}[ht]
			\centering
				\begin{tikzpicture}
				\begin{scope}[xscale = 0.9]

  \draw[line width = 2pt] (-1,1)-- +(0,-3);

  \filldraw (0,1) circle (2pt);  
  \draw[line width = 0.3pt] (0,1)-- +(0,-3);
  \filldraw (0,-2) circle (2pt);
  \node at (0, -2.5) {$1$};

  \filldraw (1,1) circle (2pt);  
  \draw[line width = 0.3pt] (1,1)-- +(0,-3);
  \filldraw (1,-2) circle (2pt);
  \node at (1, -2.5) {$2$};

	\node at (2.25,1) {$\dots$};
	\node at (2.25,-2) {$\dots$};
	
	\filldraw (3.5,1) circle (2pt); 
	\draw (3.5,1) [line width = 0.3pt] .. controls +(0,-1)  and +(0,1) ..  ++(1,-3);
	\filldraw (4.5,-2) circle (2pt); 
	\node at (3.5, -2.5) {$i$};
	
	\filldraw (4.5,1) circle (2pt); 	
	\fill[white] (4,-0.5) circle (0.2cm);
	\draw (4.5,1) [line width = 0.3pt] .. controls +(0,-1)  and +(0,1) ..  +(-1,-3);
	\filldraw (3.5,-2) circle (2pt); 
	\node at (4.5, -2.5) {$i+1$};
	
	\node at (5.5,1) {$\dots$};
	\node at (5.5,-2) {$\dots$};
		
  \filldraw (6.5,1) circle (2pt);  
  \draw[line width = 0.3pt] (6.5,1)-- +(0,-3);
  \filldraw (6.5,-2) circle (2pt);
  \node at (6.5, -2.5) {$n+1$};


  \draw[line width = 2pt] (9,1)-- +(0,-3);

  \filldraw (10,1) circle (2pt); 
  \draw[line width = 0.3pt] (10,1)-- +(0,-3);
	\filldraw (10,-2) circle (2pt); 
  \node at (10, -2.5) {$1$};
	
	\filldraw (11,1) circle (2pt); 
  \draw[line width = 0.3pt] (11,1)-- +(0,-3);
	\filldraw (11,-2) circle (2pt);
  \node at (11, -2.5) {$2$};

	\node at (12.25,1) {$\dots$};
	\node at (12.25,-2) {$\dots$};
	
	\draw (14.5,1) [line width = 0.3pt] .. controls +(0,-1)  and +(0,1) ..  +(-1,-3);
	\filldraw (13.5,1) circle (2pt);
	
	\filldraw (14.5,1) circle (2pt);
	\node at (14.5, -2.5) {$i+1$};
	\filldraw (14.5,-2) circle (2pt);
	\fill[white] (14,-0.5) circle (0.2cm);	
	
	\filldraw (13.5,-2) circle (2pt);
	\draw (13.5,1) [line width = 0.3pt] .. controls +(0,-1)  and +(0,1) ..  ++(1,-3);
	\filldraw (13.5,-2) circle (2pt);
	\node at (13.5, -2.5) {$i$};
	
	\node at (15.5,1) {$\dots$};
	\node at (15.5,-2) {$\dots$};
		
	\filldraw (16.5,1) circle (2pt);	
  \draw[line width = 0.3pt] (16.5,1)-- +(0,-3);
	\filldraw (16.5,-2) circle (2pt);
  \node at (16.5, -2.5) {$n+1$};

\end{scope}
					\end{tikzpicture}
			\caption{$\sigma_{i} ~~~~~~~~ \& ~~~~~~~~ \sigma^{-1}_{i}$}
		\end{figure}
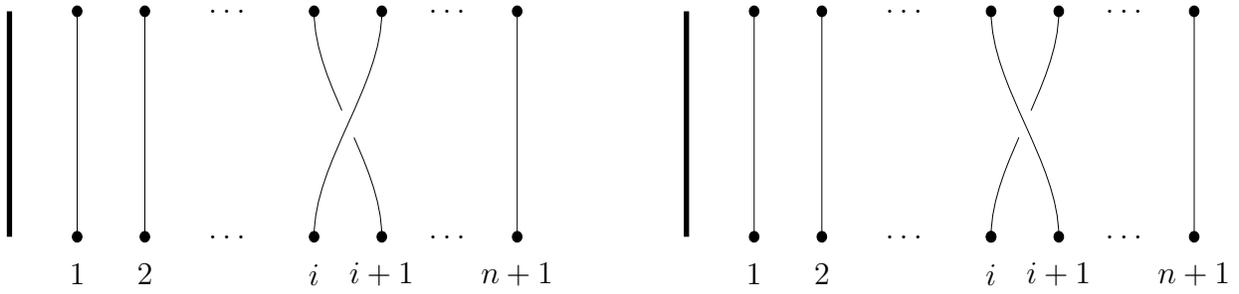
		
		\vspace{0.5cm}
		
		\begin{figure}[ht]
			\centering
			\begin{tikzpicture}
			\begin{scope}[xscale = 0.7]
  
  \filldraw (1.5,2) circle (2pt);
  \draw[line width = 0.3pt] (1.5,2).. controls +(0, -0.5) and +(0, 0.5) .. (-0.5,1).. controls +(0,-0.6) and +(-0.4,1.2) .. (1.3,-0.5).. controls +(0.2,-0.6) and +(0,1) .. (1.5, -2);
  \filldraw (1.5,-2) circle (2pt);  
  \filldraw[white] (0,0.5) circle (4pt);  
	\node at (, -2.5) {$1$};

  \draw[line width = 2pt] (0,2)-- +(0,-4);

  \filldraw (3,2) circle (2pt);  
  \draw[line width = 0.3pt] (3,2)-- +(0,-4);
  \filldraw (3,-2) circle (2pt);
  \node at (3, -2.5) {$2$};
  
  \node at (5,2) {$\dots$};
  \node at (5,-2) {$\dots$}; 
  
  \filldraw (6.5,2) circle (2pt);
  \draw[line width = 0.3pt] (6.5,2)-- +(0,-4);
  \filldraw (6.5,-2) circle (2pt);
  \node at (6.5, -2.5) {$n$};
  
  \filldraw (8,2) circle (2pt);
  \draw[line width = 0.3pt] (8,2)-- +(0,-4); 
  \filldraw (8,-2) circle (2pt);
  \node at (8, -2.5) {$n+1$};

\end{scope}

			\end{tikzpicture}
			\caption{$t$}
		\end{figure}
	
		\vspace{0.75cm}
		
	\section{Affine braids: the group $B(\tilde{A_{n}})$} \label{1_4}       
		\vspace{0.75cm}
		
		The $\tilde{A}$-type affine braid group in $n+1$ generators is the braid group under question in this work. Geometrically, one can see several presentations in the literature, among which we choose one which is compatible with our viewpoint on this group (as a base point for a better understanding of a special kind of links in the space). We show the strong connection between this group and the two groups mentioned above. The arrows in this section are well known, except for the surjection of $B(\tilde{A_{n}})$ onto $B(A_{n})$ which allows us to see $B(\tilde{A_{n}})$ as a semi-direct product of $B(A_{n})$ with a "huge" group, an infinitely generated free group (the semi-direct does not seem of a great use). In the other hand this surjection (and others induced by it) plays an essential role in \cite{Sadek_2013_1}. While concerning $B(\tilde{A_{n}})$ presentations, we define a new presentation (not far from the old one) called the parabolic-like presentation.  \\ 
	     \subsection{Presentations}
			
			\vspace{0.75cm}
			
			\begin{definition}\label{1_4_1}
				The $\tilde{A}$-type braid group with $n+1$ generators $B(\tilde{A}_{n})$ is the group presented by a set of generators $ \left\{ \sigma_{1}, \dots,  \sigma_{n}, a_{n+1}   \right\}$ together with the following defining relations\\

				\begin{itemize}[label=$\bullet$, font=\normalsize, font=\color{black}, leftmargin=2cm, parsep=0cm, itemsep=0.25cm, topsep=0cm]
					\item $\sigma_{i} \sigma_{j} =\sigma_{j} \sigma_{i} $  where $1\leq i,j\leq n$ when $ \left| i-j\right| \geq 2$,
					\item $\sigma_{i}\sigma_{i+1}\sigma_{i} = \sigma_{i+1}\sigma_{i}\sigma_{i+1}$ when $1\leq i\leq n-1$,
					\item $\sigma_{i} a_{n+1} = a_{n+1} \sigma_{i} $ when $2\leq i \leq n-1$,
					\item $\sigma_{1} a_{n+1} \sigma_{1} = a_{n+1} \sigma_{1} a_{n+1}$,
					\item $\sigma_{n} a_{n+1} \sigma_{n} = a_{n+1} \sigma_{n} a_{n+1}$. \\
				\end{itemize}	
			\end{definition}
			
			\begin{remarque}
				In the literature $a_{n+1}$ in this definition is often called $\sigma_{n+1}$, but since we are interested in viewing  $B(\tilde{A_{n-1}})$ as a subgroup of $B(\tilde{A_{n}})$ for $2\leq n $, there would be a confusion between $ a_{n} $ and $\sigma_{n}$ when they are  seen as elements in $B(\tilde{A_{n}})$. Thus we consider the group $B(A_{n})$ generated by  $ \left\{ \sigma_{1}, \dots, \sigma_{n}  \right\}$, then we "affinize" it by adding $a_{n+1}$.\\
			\end{remarque}
			
			\begin{remarque}
				It is not a coincidence that the generators here have the same symbols as the generators of $B(B_{n+1})$ in the second presentation of $B(B_{n+1})$. There is an obvious homomorphism $B(\tilde{A_{n}}) \longrightarrow B(B_{n+1})$, saying that this homomorphism is injective is equivalent to saying that in $B(B_{n+1})$ the subgroup generated by the elements $\sigma_{1}, \dots,  \sigma_{n}, a_{n+1}$ is presented by generators and relations in the following way: it has for generating set $ \left\{ \sigma_{1}, \dots,  \sigma_{n}, a_{n+1}   \right\}$ together with the relations\\
				
				\begin{itemize}[label=$\bullet$, font=\normalsize, font=\color{black}, leftmargin=2cm, parsep=0cm, itemsep=0.25cm, topsep=0cm]
					\item $\sigma_{i} \sigma_{j} =\sigma_{j} \sigma_{i} $  where $1\leq i,j\leq n$ when $ \left| i-j\right| \geq 2$,
					\item $\sigma_{i}\sigma_{i+1}\sigma_{i} = \sigma_{i+1}\sigma_{i}\sigma_{i+1}$ when $1\leq i\leq n-1$,
					\item $\sigma_{i} a_{n+1} = a_{n+1} \sigma_{i} $ when $2\leq i \leq n-1$,
					\item $\sigma_{1} a_{n+1} \sigma_{1} = a_{n+1} \sigma_{1} a_{n+1}$,
					\item $\sigma_{n} a_{n+1} \sigma_{n} = a_{n+1} \sigma_{n} a_{n+1}$. 	\\
				\end{itemize}	

				This is true, since $\phi_{n+1}$ acts on this very group by automorphisms, thus the relations between its generators  do not add any other relation than the  length-respecting  braid relations already existing in the definition. \\
			\end{remarque}
			
			\begin{proposition} \label{pr_1_4_4}
				Let $x$ be in $B(B_{n+1})$. Suppose that $x$ is expressed as a word in the generators: $ \sigma_{1}, \sigma_{2}.. \sigma_{n},t $. Then \\
		
	\centerline{$x\in B( \tilde{A_{n}}) \Longleftrightarrow $ the sum of the exponents  of $t$ in $x$  is zero.}  
	
			\end{proposition}
		\vspace{0.2cm} 
			\begin{demo}
		
				Suppose that $x= u_{1} t^{b_{1}} u_{2} t^{b_{2}} u_{3} .. u_{m} t^{b_{m}} u_{m+1} $, where $b_{i}$ is an  integer and $u_{i} $ is in  $B(A_{n})$ for all $i$. 
	
				We have $\phi_{n+1} = t\sigma_{1} .. \sigma_{n-1} \sigma_{n}$ by definition . Set $ z^{-1} = \sigma_{1} .. \sigma_{n-1} \sigma_{n}$, that gives $ t=\phi_{n+1}z $. We denote the action of $\phi^{r}_{n+1}$ on an element $ e $ in $B(\tilde{A_{n}}) $ by $[e]^{r}$, for any integer $r$. For example $\phi_{n+1}z = [z]^{1}\phi_{n+1}$.			
				\begin{eqnarray}
					\text{Now } t^{b_{i}} &=& \underbrace{\phi_{n+1}z \phi_{n+1}z  .. \phi_{n+1}z}_{b_{i} ~\text{times}} \text{,~which is equal to } \phi^{b_{i}}_{n+1} \prod\limits^{j=b_{i}}_{j=0} [z]^{j-b_{i}}. \nonumber\\\nonumber\\
					\text{Set } Z_{b_{i}} &=& \prod\limits^{j=b_{i}}_{j=0} [z]^{j-b_{i}} \text{,~which is in } B(\tilde{A_{n}}). \text{ Thus: }\nonumber\\\nonumber
				\end{eqnarray}
	
				$x= u_{1} \phi^{b_{1}}Z_{b_{1}} u_{2} \phi^{b_{2}}Z_{b_{2}} u_{3} .. u_{m} \phi^{b_{m}}Z_{b_{m}} u_{m+1} $. By pushing the "$\phi^{b_{i}} $"s to the right (acting on the "$u_{i}$"s as well as on the "$Z_{b_{1}}$"s) we get the following expression of $x$:\\ 
	
				$x =  \lambda \phi^{\big(\sum\limits^{i=m}_{i=1} b_{i}\big)}$ where $\lambda \in \tilde{A_{n}} $. This is the unique decomposition of $x$ in the semi-direct product $ \left\langle \sigma_{1}, \sigma_{2}.. \sigma_{n}, a_{n+1} \right\rangle_{B(B_{n+1})} \rtimes \Phi_{n+1}$ .\\
				
	
				Now $x$ is in $B(\tilde{A_{n}}) $ if and only if $\sum\limits^{i=m}_{i=1} b^{i} = 0 $. The proposition follows.
					
			\end{demo}
			
			
			Inspired by the injection $B(\tilde{A_{n}}) \hookrightarrow B(B_{n+1})$, we explain the geometric presentation we choose in this work. Actually, affine braids with $n+1$ generators could be viewed for example as cylindrical braids within a  cylinder labeled by $ n+1 $ points on each of its circles, the strings of the braids are not allowed to make perfect rounds. We choose to view the affine braids  as $B$-type braids: using the proposition above, the affine braids are $B$-type braids in which the number of positive rounds equals the number of negative rounds (around the fixed strand). The "$\sigma$"s are as presented above, while $a_{n+1} $ and  $\phi_{n+1}$ are the following braids:
			
		
			\begin{figure}[ht]
				\centering
				\begin{tikzpicture}
				\begin{scope}[scale = 0.6]
  \draw[line width = 2pt] (0,2)-- +(0,-4);

  \filldraw[white] (0,1.425) circle (4pt);  
  \filldraw[white] (0,-1.425) circle (4pt);

	\draw[draw=white, double =black, line width = 2pt, double distance= 0.4pt] (1.5,2).. controls +(0, -0.5) and +(0, 0.5) .. (-0.5,1).. controls +(0,-0.6) and +(-0.4,1.2) .. (9.25,0).. controls +(0.3,-0.6)  and +(0,1) .. (9.5, -2);
 	\filldraw (1.5,2) circle (2pt);
	\filldraw (9.5,-2) circle (2pt);
	\filldraw[white] (5.7,0) circle (4pt);

  \draw[draw=white, double =black, line width = 2pt, double distance= 0.4pt] (1.5,-2).. controls +(0, 0.5) and +(0, -0.5) .. (-0.50,-1).. controls +(0,0.6) and +(-0.4,-1.2) .. (9.25,0).. controls +(0.3,0.6)  and +(0,-1) .. (9.5, 2);
 	\filldraw (1.5,-2) circle (2pt);
	\filldraw (9.5,2) circle (2pt);
	
	\filldraw[white] (0,0.8) circle (4pt); 
	\filldraw[white] (0,-0.8) circle (4pt); 
	
	\filldraw[white] (0,0.8) circle (4pt); 
	\filldraw[white] (0,-0.8) circle (4pt);

 	\filldraw[white] (3,0.75) circle (4pt);
	\filldraw[white] (3,-0.75) circle (4pt);
	
 	\filldraw[white] (8,0.45) circle (4pt);
	\filldraw[white] (8,-0.45) circle (4pt);
	
	\draw[draw= white, double = black, line width = 0.3pt, double distance= 2pt ] (0,1)-- +(0,-2);
	\draw[draw= white, double = black, line width = 0.3pt] (3,2) -- +(0,-4);
	
	\node at (5,2) {$\dots$};
	\node at (5,-2) {$\dots$};
	
  \draw[draw= white, double = black, line width = 0.3pt] (8,2) -- +(0,-4);

	\node at (1.5, -2.5) {$1$};
	\node at (3, -2.5) {$2$};
	\node at (8, -2.5) {$n$};
	\node at (9.5, -2.5) {$n+1$};
	
 	\filldraw (3,-2) circle (2pt);
	\filldraw (3,2) circle (2pt);
	
 	\filldraw (8,-2) circle (2pt);
	\filldraw (8,2) circle (2pt);
		
	\node at (4.5, -4) {$(a_{n+1})$};

\end{scope}


\begin{scope}[scale = 0.6, xshift = 14cm, yshift = 7cm]

  \draw[line width = 2pt] (0,-5)-- +(0,-4);
	\filldraw[white] (0,-6.35) circle (4pt); 
			
	\filldraw (1.5,-5) circle (2pt);  
  \draw[line width = 0.4pt] (1.5,-5).. controls +(0, -0.5) and +(0, 0.5) .. (-0.5,-6).. controls +(0,-0.6) and +(-0.4,0.0) .. (9.5,-9);
	\filldraw (9.5,-9) circle (2pt); 
	\node at (9.5, -9.5) {$n+1$};
	
	\draw[draw= white, double = black, line width = 2pt, double distance= 0.4pt ] (3,-5)-- +(-1.5,-4);
	\filldraw (3,-5) circle (2pt); 
  \filldraw (1.5,-9) circle (2pt);
  \node at (1.5, -9.5) {$1$};

	\draw[draw= white, double = black, line width = 2pt, double distance= 0.4pt ] (4.5,-5)-- +(-1.5,-4);
  \filldraw (4.5,-5) circle (2pt); 
  \filldraw (3,-9) circle (2pt);
  \node at (3, -9.5) {$2$};
	
	\node at (6.5,-5) {$\dots$};
	\node at (4.5,-9) {$\dots$};
	
	\draw[draw= white, double = black, line width = 2pt, double distance= 0.4pt ] (8,-5)-- +(-1.5,-4);
	\filldraw (8,-5) circle (2pt); 
  \filldraw (6.5,-9) circle (2pt);
  \node at (6.5, -9.5) {$n-1$};

	\draw[draw= white, double = black, line width = 2pt, double distance= 0.4pt ] (9.5,-5)-- +(-1.5,-4);
	\filldraw (9.5,-5) circle (2pt); 
  \filldraw (8,-9) circle (2pt);
  \node at (8, -9.5) {$n$};

	\node at (5, -11) {$(\phi_{n+1})$};

\end{scope}

				\end{tikzpicture}
				\caption{$ $}
			\end{figure}								

			Notice that $\phi_{n+1}$ is not an affine braid.\\
			
			\begin{proposition}\label{1_4_5}
				The following group homomorphism is injective							
				\begin{eqnarray}
					F_{n}: B(\tilde{A_{n-1}}) &\longrightarrow& B(\tilde{A_{n}})  \nonumber\\
					\sigma_{i} &\longmapsto& \sigma_{i}$ ~~~ \text{for} $1\leq i\leq n-1 \nonumber\\
					a_{n} &\longmapsto& \sigma_{n} a_{n+1}\sigma^{-1}_{n} \nonumber
				\end{eqnarray}
			\end{proposition}
			\begin{demo}
			See \ref{pr_1_4_4}: $F_{n}$ is a restriction of the injection $y_{n}$ to $B(\tilde{A_{n}}) $.
			\end{demo}
			We give now a new presentation of  $B(\tilde{A_{n}})$, in which the defining relations are positive, and where  $F_{n}$ is obtained by simply adding one generator to those of $B(\tilde{A_{n-1}})$.  \\
                 
              	By definition $B(\tilde{A_{n}}) $ has $ \left\{ \sigma_{1}, \dots, \sigma_{n}, a_{n+1}   \right\}$ as a set of generators together with the following defining relations:\\
         	
			\begin{itemize}[label=$\bullet$, font=\normalsize, font=\color{black}, leftmargin=2cm, parsep=0cm, itemsep=0.25cm, topsep=0cm]
				 \item[(1')] $\sigma_{i} \sigma_{j} =\sigma_{j} \sigma_{i} $  where $1\leq i,j\leq n$ when $ \left| i-j\right| \geq 2$,
				 \item[(2')] $\sigma_{i}\sigma_{i+1}\sigma_{i} = \sigma_{i+1}\sigma_{i}\sigma_{i+1}$ when $1\leq i\leq n-1$,
				 \item[(3')] $\sigma_{i} a_{n+1} = a_{n+1} \sigma_{i} $ when $2\leq i \leq n-1$, 
				 \item[(4')] $\sigma_{1} a_{n+1} \sigma_{1} = a_{n+1} \sigma_{1} a_{n+1}$,
				 \item[(5')] $\sigma_{n} a_{n+1} \sigma_{n} = a_{n+1} \sigma_{n} a_{n+1}$ for $n\geq 2$.\\
			\end{itemize}
         	
			This presentation is to be called the formal presentation from now on. For the moment $n$ is to be greater than or equal to 3. $B(\tilde{A_{2}})$ is generated by $\sigma_{1}, \sigma_{2}$ and $a_{3}$, we see that $a_{n+1} $ 
			can be seen as follows: 

			\begin{eqnarray}
				a_{n+1} = \sigma_{n}^{-1} \dots  \sigma_{3}^{-1} a_{3} \sigma_{3} \dots \sigma_{n}. \nonumber
			\end{eqnarray}	
			
			Our aim is to show that $B(\tilde{A_{n}})$ can be generated by $ S'= \left\{ \sigma_{1}, \dots,  \sigma_{n}, a_{3}   \right\}$ with defining relations:\\
         
			\begin{itemize}[label=$\bullet$, font=\normalsize, font=\color{black}, leftmargin=2cm, parsep=0cm, itemsep=0.25cm, topsep=0cm]
				 \item[(1)] $\sigma_{i} \sigma_{j} =\sigma_{j} \sigma_{i} $  where $1\leq i,j\leq n$ when $ \left| i-j\right| \geq 2$,
				 \item[(2)] $\sigma_{i}\sigma_{i+1}\sigma_{i} = \sigma_{i+1}\sigma_{i}\sigma_{i+1}$ when $1\leq i\leq n-1$, 
				 \item[(3)] $\sigma_{1}a_{3}\sigma_{1} = a_{3}\sigma_{1}a_{3}$,
				 \item[(4)] $\sigma_{3}a_{3}\sigma_{3} = a_{3}\sigma_{3}a_{3}$,
				 \item[(5)] $\sigma_{i}a_{3} = a_{3}\sigma_{i}$ when  $4\leq i\leq n$,
				 \item[(6)] $\sigma_{3}\sigma_{2}a_{3}\sigma_{3} = \sigma_{2}a_{3}\sigma_{3}\sigma_{2}$.\\
			\end{itemize}  
         
			First we show that the formal presentation gives the new one, now $ \left\{ \sigma_{1}, \dots,  \sigma_{n}, a_{3}   \right\}$ generates $B(\tilde{A_{n}})$ indeed, since			
			\begin{eqnarray}
				a_{n+1} = \sigma_{n}^{-1} \dots \sigma_{3}^{-1} a_{3} \sigma_{3} 
				\dots  \sigma_{n}, ~~~~\text{   which gives }~~~~a_{3} = \sigma_{3} \dots  \sigma_{n} a_{n+1} \sigma_{n}^{-1} \dots \sigma_{3}^{-1}.\nonumber
			\end{eqnarray}
	
			We see that (3) follows directly from the fact that $a_{3}$ and $a_{n+1}$ are conjugate, while (5) could be seen to be valid geometrically or by a direct computation: 			
			\begin{eqnarray}
			\text{for} ~4\leq i\leq n,~ \text{we have}:	~\sigma_{i}a_{3} &=& \sigma_{i} \sigma_{3} .. \sigma_{n} a_{n+1} \sigma_{n}^{-1} .. \sigma_{3}^{-1} = \sigma_{3} .. \sigma_{n} \sigma_{i-1} a_{n+1} \sigma_{n}^{-1} .. \sigma_{3}^{-1}\nonumber\\\nonumber\\
				&=& \sigma_{3} .. \sigma_{n}  a_{n+1} \sigma_{i-1} \sigma_{n}^{-1} .. \sigma_{3}^{-1} = \sigma_{3} .. \sigma_{n}  a_{n+1} \sigma_{n}^{-1} .. \sigma_{3}^{-1} \sigma_{i} \nonumber\\\nonumber\\
				&=& a_{3}\sigma_{i}.\nonumber
			\end{eqnarray}	
			
			Now we treat (4):			
			\begin{eqnarray}	
				\sigma_{3}a_{3}\sigma_{3} &=& \sigma_{3}^{2} .. \sigma_{n} a_{n+1} \sigma_{n}^{-1} .. \sigma_{4}^{-1}.\nonumber
			\end{eqnarray}
			
			\begin{eqnarray}
				\text{Thus,}~a_{3}\sigma_{3}a_{3} &=& \sigma_{3} .. \sigma_{n} a_{n+1} \underbrace{\sigma_{n}^{-1} .. \sigma_{4}^{-1} \sigma_{3} \sigma_{4} .. \sigma_{n}}_{} a_{n+1} \sigma_{n}^{-1} .. \sigma_{3}^{-1}\nonumber\\\nonumber\\
				&=& \sigma_{3} .. \sigma_{n} a_{n+1} \sigma_{3} .. \sigma_{n-1} \sigma_{n} \sigma_{n-1}^{-1} .. \sigma_{3}^{-1} a_{n+1} \sigma_{n}^{-1} .. \sigma_{3}^{-1} \nonumber\\\nonumber\\
				&=& \sigma_{3} .. \sigma_{n}  \sigma_{3} .. \sigma_{n-1} a_{n+1}\sigma_{n} a_{n+1}\sigma_{n-1}^{-1} ..\sigma_{3}^{-1} \sigma_{n}^{-1} .. \sigma_{3}^{-1} \nonumber\\\nonumber\\
				&=& \sigma_{3} .. \sigma_{n}  \sigma_{3} .. \sigma_{n-1} \sigma_{n}a_{n+1} \sigma_{n}\sigma_{n-1}^{-1} ..\sigma_{3}^{-1} \sigma_{n}^{-1} .. \sigma_{3}^{-1}\nonumber\\\nonumber\\
				&=& \sigma_{3} .. \sigma_{n}  \sigma_{3} .. \sigma_{n-1} \sigma_{n}a_{n+1} \sigma_{n} \sigma_{n}^{-1} .. \sigma_{3}^{-1}\sigma_{n}^{-1} ..\sigma_{4}^{-1} \nonumber\\\nonumber\\
				&=& \sigma_{3} .. \sigma_{n}  \underbrace{\sigma_{3} .. \sigma_{n-1} \sigma_{n} \sigma_{n-1}^{-1} .. \sigma_{3}^{-1}}_{} a_{n+1} \sigma_{n}^{-1} ..\sigma_{4}^{-1} \nonumber\\\nonumber\\
				&=& \sigma_{3} .. \sigma_{n}  \sigma_{n}^{-1} .. \sigma_{4}^{-1} \sigma_{3} \sigma_{4} .. \sigma_{n} a_{n+1} \sigma_{n}^{-1} ..\sigma_{4}^{-1} \nonumber\\\nonumber\\
				&=& \sigma_{3}^{2}\sigma_{4} .. \sigma_{n} a_{n+1} \sigma_{n}^{-1} ..\sigma_{4}^{-1}\nonumber\\\nonumber\\
				&=& \sigma_{3}a_{3}\sigma_{3}. \nonumber	\\\nonumber
			\end{eqnarray}

			We see that (6) is equivalent to    			
			\begin{eqnarray}
				\sigma_{2}\sigma_{3}^{-1}a_{3}\sigma_{3} = \sigma_{3}^{-1}a_{3}\sigma_{3}\sigma_{2}~~~~~~ ... (6''). \nonumber
			\end{eqnarray}
         
			But $\sigma_{3}^{-1}a_{3}\sigma_{3} = \sigma_{4} .. \sigma_{n} a_{n+1} \sigma_{n}^{-1} .. \sigma_{4}^{-1}$, and (6) follows.\\
			
			Now we show that the new presentation gives the formal one.\\
         
			We are reduced to show that the new presentation gives (3'), (4') and (5'), i.e., the formal relations which involve $a_{n+1}$ have to be shown using the new relations with: 			
			\begin{eqnarray}
				a_{n+1} = \sigma_{n}^{-1} .. \sigma_{3}^{-1} a_{3} \sigma_{3} .. \sigma_{n}.\nonumber
			\end{eqnarray}
                   
			We start by dealing with (3'): let $3\leq i \leq n-1$. We compute:
			
			\begin{eqnarray}
				\sigma_{i} a_{n+1} &=& \sigma_{i} \sigma_{n}^{-1} .. \sigma_{3}^{-1} a_{3} \sigma_{3} .. \sigma_{n} = \sigma_{n}^{-1} .. \sigma_{3}^{-1} \underbrace{\sigma_{i+1} a_{3}}_{(5)} \sigma_{3} .. \sigma_{n}\nonumber\\\nonumber\\
				&=& \sigma_{n}^{-1} .. \sigma_{3}^{-1} a_{3} \sigma_{i+1}\sigma_{3} .. \sigma_{n} = \sigma_{n}^{-1} .. \sigma_{3}^{-1} a_{3} \sigma_{3} .. \sigma_{n} \sigma_{i} \nonumber\\\nonumber\\
				&=& a_{n+1} \sigma_{i}. \nonumber
			\end{eqnarray}
			
			\vspace{-0.5cm}		
			\begin{eqnarray}			
			\text{Moreover, }\sigma_{2} a_{n+1} &=& \sigma_{2} \sigma_{n}^{-1} .. \sigma_{3}^{-1} a_{3} \sigma_{3} .. \sigma_{n} = \sigma_{n}^{-1} .. \underbrace{\sigma_{2} \sigma_{3}^{-1} a_{3} \sigma_{3}}_{(6'')} .. \sigma_{n} = \sigma_{n}^{-1} .. \sigma_{3}^{-1} a_{3} \sigma_{3} \sigma_{2}  .. \sigma_{n} \nonumber\\
				&=& a_{n+1} \sigma_{2}.  \nonumber\\\nonumber
			\end{eqnarray}			
			
			Hence, $\sigma_{i} a_{n+1} = a_{n+1} \sigma_{i} $ when $2\leq i \leq n-1$. Thus, (3') is proved.\\
          
			Since $a_{n+1}$ and $a_{3}$ are conjugate, (3) gives directly (4').\\
          
			Now we want to show that $\underbrace{\sigma_{n} a_{n+1} \sigma_{n}}_{:=x} = \underbrace{a_{n+1} \sigma_{n} a_{n+1}}_{:=y}.$			
			\begin{eqnarray}
				x &=& \sigma_{n-1}^{-1} .. \sigma_{3}^{-1} a_{3} \sigma_{3} .. \sigma_{n}^{2},\nonumber\\\nonumber\\
				y &=& \sigma_{n}^{-1} .. \sigma_{3}^{-1} a_{3} \sigma_{3} .. \sigma_{n-1}\sigma_{n}\sigma_{n-1}^{-1} .. \sigma_{3}^{-1} a_{3} \sigma_{3} .. \sigma_{n} = \sigma_{n}^{-1} .. \sigma_{3}^{-1} a_{3} \sigma_{n}^{-1} .. \sigma_{4}^{-1}\sigma_{3}\sigma_{4} .. \sigma_{n} a_{3} \sigma_{3} .. \sigma_{n}\nonumber\\\nonumber\\
				&=& \sigma_{n}^{-1} .. \sigma_{3}^{-1}  \sigma_{n}^{-1} .. \sigma_{4}^{-1} \underbrace{a_{3}\sigma_{3}a_{3}}_{(4)}\sigma_{4} .. \sigma_{n} \sigma_{3} .. \sigma_{n} = \sigma_{n}^{-1} .. \sigma_{3}^{-1}  \sigma_{n}^{-1} .. \sigma_{4}^{-1}\sigma_{3}a_{3}\sigma_{3}\sigma_{4} .. \sigma_{n} \sigma_{3} .. \sigma_{n}\nonumber\\\nonumber\\
				&=& \sigma_{n-1}^{-1} .. \sigma_{3}^{-1}\sigma_{n}^{-1} .. \sigma_{3}^{-1} \sigma_{3}a_{3}\sigma_{3}\sigma_{4} .. \sigma_{n} \sigma_{3} .. \sigma_{n}  = \sigma_{n-1}^{-1} .. \sigma_{3}^{-1}a_{3}\sigma_{n}^{-1} .. \sigma_{4}^{-1}\sigma_{3}\sigma_{4} .. \sigma_{n} \sigma_{3} .. \sigma_{n} \nonumber\\\nonumber\\
				&=&  \sigma_{n-1}^{-1} .. \sigma_{3}^{-1}a_{3}\sigma_{3} .. \sigma_{n-1}\sigma_{n}\sigma_{n-1}^{-1} .. \sigma_{3}^{-1} \sigma_{3} .. \sigma_{n}  = \sigma_{n-1}^{-1} .. \sigma_{3}^{-1}a_{3}\sigma_{3} .. \sigma_{n-1}\sigma_{n}^{2} \nonumber\\\nonumber\\
				&=& x. \nonumber
			\end{eqnarray}
			\vspace{0.25cm}
				
			Finally, $F_{n}(a_{3}) = F_{n}(\sigma_{3} .. \sigma_{n-1} a_{n} \sigma_{n-1}^{-1} .. \sigma_{3}^{-1} ) = \sigma_{3} .. \sigma_{n-1} F_{n}(a_{n}) \sigma_{n-1}^{-1} .. \sigma_{3}^{-1}$. This is equal to $ \sigma_{3} .. \sigma_{n} a_{n+1} \sigma_{n}^{-1} .. \sigma_{3}^{-1}$, thus to $  a_{3}$.\\
			
          	
			Now $F_{n}$ would have the following form with the new presentation:\\		
			\begin{eqnarray}
				F_{n}: B(\tilde{A_{n-1}}) &\longrightarrow& B(\tilde{A_{n}}) \nonumber\\
				\sigma_{i} &\longmapsto& \sigma_{i}$ for $1\leq i\leq n-1\nonumber\\
				a_{3} &\longmapsto& a_{3}\nonumber
			\end{eqnarray}
			
			Notice that we could have started with the group $B(\tilde{A_{1}})$, which is a free group in two letters $\sigma_{1}$ and $a_{3}$, with a change in the sixth relation. On the other hand, it is obvious that $(B(\tilde{A_{n+1}}),S')$ is not an Artin System.\\
			
			 \subsection{$B(A_{n})$ as a quotient of $B(\tilde{A_{n}})$}\label{1_4_2} 
			
			Now we consider the element $e$ in $ B(\tilde{A_{n}})$ given as  			
			\begin{eqnarray}
				e =a^{-1}_{n+1}\sigma^{-1}_{n}\sigma^{-1}_{n-1} ...\sigma^{-1}_{2}\sigma_{1}\sigma_{2} ...\sigma_{n-1}\sigma_{n} = a^{-1}_{n+1}\sigma_{1}\sigma_{2} ...\sigma_{n-1}\sigma_{n}\sigma^{-1}_{n-1} ...\sigma^{-1}_{2}\sigma^{-1}_{1}.\nonumber
			\end{eqnarray}
			
			Let $N_{e}$ be the normal subgroup of $B(\tilde{A_{n}})$ generated by $e$. Consider the quotient $Q = B(\tilde{A_{n}})/N_{e}$.\\

			\begin{lemme}
				In $B(\tilde{A_{n}})$ (and in $B({A_{n}})$ as well), the element $b:= \sigma^{-1}_{n}\sigma^{-1}_{n-1} ...\sigma^{-1}_{2}\sigma_{1}\sigma_{2} ...\sigma_{n-1}\sigma_{n}$ verifies the following relations:\\
								
				\begin{enumerate}[label=\arabic*), font=\normalsize, font=\color{black}, leftmargin=2cm, parsep=0cm, itemsep=0.25cm, topsep=0cm]
					\item $\sigma_{i} b = b \sigma_{i}$ for $ 2 \leq i \leq n-1 $,
					\item $\sigma_{n} b \sigma_{n} = b \sigma_{n} b$,
					\item $\sigma_{1} b \sigma_{1} = b \sigma_{1} b$.\\
				\end{enumerate}
			\end{lemme}
			
			\begin{demo}
				\begin{enumerate}[label=\arabic*), font=\normalsize, font=\color{black}, leftmargin=1cm, parsep=0cm, itemsep=0.25cm, topsep=0cm]
					\item  
					\begin{eqnarray}
						\sigma_{i} b &=& \sigma_{i} \sigma^{-1}_{n}\sigma^{-1}_{n-1} ...\sigma^{-1}_{2}\sigma_{1}\sigma_{2} ...\sigma_{n-1}\sigma_{n} = \sigma^{-1}_{n}...\underbrace{\sigma_{i}\sigma^{-1}_{i+1}\sigma^{-1}_{i}}_{\sigma^{-1}_{i+1}\sigma^{-1}_{i}\sigma_{i+1}}...\sigma_{1}\sigma_{2}  ...\sigma_{n-1}\sigma_{n} \nonumber\\\nonumber\\
						&=& \sigma^{-1}_{n}\sigma^{-1}_{n-1} ...\sigma^{-1}_{2} \sigma_{1} ... \underbrace{\sigma_{i+1}\sigma_{i}\sigma_{i+1}}_{\sigma_{i}\sigma_{i+1}\sigma_{i}} .. \sigma_{1}= b \sigma_{i}. \nonumber
					\end{eqnarray}
					\item [3)] 
					\begin{eqnarray}
						& & \sigma^{-1}_{n}\sigma^{-1}_{n-1} ...\sigma^{-1}_{2}\sigma_{1}\sigma_{2} ...\sigma_{n-1}\sigma_{n} \sigma_{1} \sigma^{-1}_{n}\sigma^{-1}_{n-1} ...\sigma^{-1}_{2}\sigma_{1}\sigma_{2} ...\sigma_{n-1}\sigma_{n}\nonumber\\\nonumber\\
						&=& \sigma^{-1}_{n}\sigma^{-1}_{n-1} ...\sigma^{-1}_{2}\underbrace{\sigma_{1}\sigma_{2}\sigma_{1}}_{\sigma_{2}\sigma_{1} \sigma_{2}} \sigma_{2}^{-1} \sigma_{1}\sigma_{2} ...\sigma_{n-1}\sigma_{n} = \sigma^{2}_{1}\sigma^{-1}_{n}\sigma^{-1}_{n-1} ...\sigma^{-1}_{3}\sigma_{2}\sigma_{3} ...\sigma_{n-1}\sigma_{n} = \sigma_{1} b \sigma_{1}. \nonumber\\\nonumber
					\end{eqnarray}
				\end{enumerate}

				In the same way we deal with (2), hence the proof is done.
     
			\end{demo}
			
			We see directly that, when replacing $b$ by $a_{n+1}$ in this lemma, we get the defining relations of  $B(\tilde{A_{n}})$ in which $a_{n+1}$ is involved. Now the group $Q$ is generated by the set $ \left\{ \sigma_{1}, \dots,  \sigma_{n}, a_{n+1} \right\}$ with the defining relations :\\
			\begin{itemize}[label=$\bullet$, font=\normalsize, font=\color{black}, leftmargin=2cm, parsep=0cm, itemsep=0.25cm, topsep=0cm]
				\item[$\bullet$] $\sigma_{i} \sigma_{j} =\sigma_{j} \sigma_{i} $  where $1\leq i,j\leq n$ when $ \left| i-j\right| \geq 2$,
				\item[$\bullet$] $\sigma_{i}\sigma_{i+1}\sigma_{i} = \sigma_{i+1}\sigma_{i}\sigma_{i+1}$ when $1\leq i\leq n-1$, 
				\item[$\bullet$] $\sigma_{i} a_{n+1} = a_{n+1} \sigma_{i} $ when $2\leq i \leq n-1$,
				\item[$\bullet$] $\sigma_{1} a_{n+1} \sigma_{1} = a_{n+1} \sigma_{1} a_{n+1}$ and $\sigma_{n} a_{n+1} \sigma_{n} = a_{n+1} \sigma_{n} a_{n+1}$. for $n\geq2$,
				\item[$\bullet$] $a_{n+1} = \sigma^{-1}_{n}\sigma^{-1}_{n-1} ...\sigma^{-1}_{2}\sigma_{1}\sigma_{2} ...\sigma_{n-1}\sigma_{n}$.\\
			\end{itemize}
						
			\begin{theoreme}\label{1_4_7}
				The map $f : Q \longrightarrow B(A_{n})$ defined by 				
				\begin{eqnarray}
					~~~~~~~\sigma_{i} &\longmapsto& \sigma_{i}$ for $1\leq i\leq n, \nonumber\\
					~~~~~~~a_{n+1} &\longmapsto& \sigma^{-1}_{n}\sigma^{-1}_{n-1} ...\sigma^{-1}_{2}\sigma_{1}\sigma_{2} ...\sigma_{n-1}\sigma_{n}, \nonumber
				\end{eqnarray}

				is a group isomorphism.\\
			\end{theoreme}
			
			\begin{demo}
      
				By the lemma we see that $f$ is a homomorphism, and it is surjective. Moreover, the following map:	
				\vspace{-0.5cm}
				\begin{eqnarray}
					g:B(A_{n}) &\longrightarrow& Q  ~~\text{given by} \nonumber\\
					\sigma_{i} &\longmapsto& \sigma_{i}$ for $1\leq i\leq n, \nonumber
				\end{eqnarray}
				
				is a group homomorphism, surjective indeed, for $\sigma^{-1}_{n}\sigma^{-1}_{n-1} ...\sigma^{-1}_{2}\sigma_{1}\sigma_{2} ...\sigma_{n-1}\sigma_{n}$ is send to itself, hence to $a_{n+1}$.\\
       
				Obviously: $f\circ g = Id_{B(A_{n})}$ and $g\circ f = Id_{Q}$, so $f$ is indeed an isomorphism with $f^{-1} = g$. 
				
			\end{demo}

			With the notation as above: 			
			\begin{eqnarray}
				B(\tilde{A_{n}}) \cong  B(A_{n}) \ltimes  N_{e} \nonumber
			\end{eqnarray}
			
			Clearly this surjection respects the inclusion $F_{n}$, in other terms the following diagram commutes:  			
		 \begin{figure}[ht]
			 \centering

		\begin{tikzpicture}

			\matrix[matrix of math nodes,row sep=1cm,column sep=1cm]{
			|(A)| B(\tilde{A_{n-1}})   & & & & &  |(B)| B(\tilde{A_{n}})   \\
			                           & & & & &                           \\								
			                           & & & & &                           \\								
			|(C)|  B(A_{n-1})          & & & & &  |(D)| B(A_{n})           \\
				};

				\path (B) edge[-myhook,line width=0.42pt] node[above, xshift=2mm, yshift=0mm, rotate=0] {\footnotesize $F_{n}$}     (A);
				\path (A) edge[-myto,line width=0.42pt]      (B);
				
				\path (D) edge[-myhook,line width=0.42pt] node[above, xshift=1mm, yshift=0mm, rotate=0] {\footnotesize $x_{n}$}    (C);
				\path (C) edge[-myto,line width=0.42pt]      (D);
		
				\path (A) edge[-myonto,line width=0.42pt] node[above, xshift=-5mm, yshift=-5mm, rotate=0] {\footnotesize $\beta_{n-1}$}     (C);

				\path (B) edge[-myonto,line width=0.42pt] node[above, xshift=4mm, yshift=-4.5mm, rotate=0] {\footnotesize $\beta_{n}$}     (D);

		\end{tikzpicture}
			 \end{figure}

	     \subsection{On the structure of $B(\tilde{A_{n}})$: Schreier's theorem}
			
			\vspace{0.5cm}
			
			\begin{lemme}
				Let $G$ be a free group with $S=\left\{ F_{0}, F_{1}, .. ,F_{n} \right\}$ as a set of free generators, let $N$ be the subgroup:
				\begin{eqnarray}
					N= \left\{ F^{\epsilon_{1}}_{i_{1}} .. F^{\epsilon_{j}}_{i_{j}} ..  F^{\epsilon_{m}}_{i_{m}}; \sum \epsilon_{j} =0 \right\} \nonumber
				\end{eqnarray}
				
				That is the normal group of words in which the number of positive letters is equal to the number of negative letters.\\
    
				Then, $N$ is an infinitely generated free group over the letters $F^{j}_{0} F_{i} F^{-(j+1)}_{0}$ where $j$ is an arbitrary integer .\\
			\end{lemme}

			\begin{demo}
    
				By Schreier's theorem we see that $N$ is a free group. Now we apply Schreier's algorithm to show that $F^{j}_{0} F_{i} F^{-(j+1)}_{0}$ could be viewed as free generators, to do so we consider the following surjection:
				\begin{eqnarray}
					G &\longrightarrow& \left\langle F_{0} \right\rangle  ~~\text{given by:}\nonumber\\
					F_{i} &\mapsto& F_{0}$ for $1\leq i\leq n, \nonumber		
				\end{eqnarray}
				
				where $\left\langle F_{0} \right\rangle$ is to be the subgroup of $G$ generated by $F_{0}$, that is the free group with one generator $F_{0}$. It is obvious that $N$ is the kernel of this homomorphism (we could actually choose any of the letters of $G$ to generate a range in order to have $N$ as a kernel), hence
				\begin{eqnarray}
					 G \cong  \left\langle F_{0} \right\rangle \ltimes  N. \nonumber
				\end{eqnarray}

				So we can consider the set $\left\{ F^{j}_{0}; j\in \mathbb{Z} \right\}$ as a full set of representatives of right cosets of $N$ in $G$, thus it could be considered as a "Schreier's system". Now we define the following map:
				\begin{eqnarray}
					\phi: G &\longrightarrow& \left\langle F_{0} \right\rangle ~~\text{given by:} \nonumber\\
					\phi(x) &=& F^{k}_{0}$ when  $x\in N F^{k}_{0}. \nonumber
				\end{eqnarray}

				By Schreier's theorem:  $ gs \phi(gs)^{-1} $ are free generators of $N$, where $g$ runs over the "Schreier's system" and $s$ runs over the set of free generators of $G$ .\\
				
				The element $gs$ is of the form $F^{j}_{0}F_{i}$, while $ \phi(gs)^{-1} $ is of the form:				
				\begin{eqnarray}
					\phi(F^{j}_{0}F_{i})^{-1} = \phi (F^{-1}_{i}F^{-j}_{0}) = \phi(\underbrace{F^{-1}_{i}F^{-j}_{0}F^{j+1}_{0}}_{ \in N} F^{-(j+1)}_{0}) = F^{-(j+1)}_{0}. \nonumber\\\nonumber
				\end{eqnarray}

				Thus $ gs \phi(gs)^{-1} $ is of the form $F^{j}_{0} F_{i}F^{-(j+1)}_{0} $ where $ 0 \leq i \leq n $. The lemma follows.
       
			\end{demo}
		
			\vspace{0.75cm}
			
			Now the following diagram commutes:
			
			\clearpage
			\begin{figure} 
				\centering
				

		\begin{tikzpicture}

			\matrix[matrix of math nodes,row sep=1cm,column sep=1cm]{
			|(A)| B(\tilde{A_{n}}) & & & &    |(B)| B(A_{n})  \\
			                       & & & &                   \\								
			|(C)| B(B_{n+1})       & & & &                    \\
				};
				
				\path (A) edge[-myonto,line width=0.42pt]  node[above, xshift=0mm, yshift=1mm, rotate=0] {\footnotesize $f$}    (B);

				\path (C) edge[-myhook,line width=0.42pt] node[below, xshift=-3mm, yshift=3mm, rotate=0] {\footnotesize $i_{n}$}    (A);
				\path (A) edge[-myto,line width=0.42pt]      (C);
								
				\path (C) edge[-myonto,line width=0.42pt]   node[above, xshift=5mm, , yshift=-4mm, rotate=0] {\footnotesize $\alpha_{n}$}   (B);

	 		\end{tikzpicture}

			\end{figure}
		
			We can see that $ker(f) = ker(\alpha_{n}) \cap i_{n}(B(\tilde{A_{n}})) $. But $N_{e}= ker(f)$, so we are reduced to understanding the structure of $ ker(\alpha_{n})$. Set $G':= ker(\alpha_{n})$.
			From  \cite{Digne_Gomi_2001},  section 5 and proposition 16, we get:   				
				
			\begin{lemme}\label{1_4_10} 
			\cite{Digne_Gomi_2001}With the above notations, $G'$ is a free subgroup of $B(B_{n+1})$ generated by the free generators $F_{i}$, where $ 0 \leq i \leq n, ~F_{i} = \sigma_{i} .. \sigma_{1} t \sigma^{-1}_{1} .. \sigma^{-1}_{i} $ and  $F_{0} = t $.
			\end{lemme}

			An element $x$ is in $N_{e}$ if and only if  it is a word in $G'$ and it is in $B(\tilde{A_{n}})$, but since any element of $B(B_{n+1})$ is in  $B(\tilde{A_{n}})$ if and only if the sum of exponents of $t$  is zero in an (every) expression of it, so we can apply the last lemma, taking $G'$ for $G$. Then $N_{e}$ is to be $N$, so our group $N_{e}$ is an  infinitely generated free group.
			
\subsection{$B(A_{n+1})$, $B(B_{n})$ and $B(\tilde{A_{n}})$, arrows}\label{1_4_4}
 
	$ \   $
				\vspace{0.5cm}
 
  In what follows we show the net of arrows between the three types of braid groups mentioned above. We are interested in investigating which arrows among those do respect the injections between groups of a given type in different number of generators. Roughly speaking: the arrows should be thought of as arrows defined over the "towers of groups", precisely the towers come from the injections between groups of the same types. 
			
			Consider $x_{n}$, the injection  $B(A_{n-1}) \hookrightarrow B(A_{n})$ mentioned in 1.2.1. Geometrically $B(B_{n})$ embeds into $B(B_{n+1})$ by adding the $n+1$-th stand, that is: 

\begin{figure}[h]
				\centering
				\begin{tikzpicture}
				\begin{scope}[xscale = 0.85]

				\draw[line width = 2pt] (-1,4)-- (-1,-1);

				\draw[draw= white, double = black, line width = 0.1cm] (-0.5, 1.2) .. controls +(-0.3, 0.3) and +(0, -0.3) .. (-1.4, 1.5) .. controls +(0,0.3) and +(-0.3,-0.3) .. (-0.5,1.8); 

				\filldraw[white] (-1,1.35) circle (4pt);
				\draw[line width = 2pt] (-1,1.5)-- (-1,0);


			 \filldraw (0,4) circle (2pt);
				\draw[line width = 0.3pt] (0,3) -- (0,4);

				\filldraw (1,4) circle (2pt);
				\draw[line width = 0.3pt] (1,3) -- (1,4);

				\filldraw (5,4) circle (2pt);
				\draw[line width = 0.3pt] (5,3) -- (5,4);

				\draw[line width = 0.3pt] (-0.5,3) -- (5.5,3);

				\draw[line width = 0.3pt] (5.5,3) -- (5.5,0);

				\draw[line width = 0.3pt] (5.5,0) -- (-0.5,0);
				
				\draw[line width = 0.3pt] (-0.5,0) -- (-0.5,3);  
				
				\draw[line width = 0.3pt] (0,0) -- (0,-1);
				\filldraw (0,-1) circle (2pt);
				\node at (0, -1.5) {$1$};

				\draw[line width = 0.3pt] (1,0) -- (1,-1);
				\filldraw (1,-1) circle (2pt);
				\node at (1, -1.5) {2};

				\draw[line width = 0.3pt] (5,0) -- (5,-1);
				\filldraw (5,-1) circle (2pt);
				\node at (5, -1.5) {n};

				\draw[line width =1.5pt][->] (6.25,1.5) -- (7.5,1.5);

\end{scope}

\begin{scope}[xscale = 0.85, xshift = 10cm]

				\draw[line width = 2pt] (-1,4)-- (-1,-1);

				\draw[draw= white, double = black, line width = 0.1cm] (-0.5, 1.2) .. controls +(-0.3, 0.3) and +(0, -0.3) .. (-1.4, 1.5) .. controls +(0,0.3) and +(-0.3,-0.3) .. (-0.5,1.8); 

				\filldraw[white] (-1,1.35) circle (4pt);
				\draw[line width = 2pt] (-1,1.5)-- (-1,0);


			 \filldraw (0,4) circle (2pt);
				\draw[line width = 0.3pt] (0,3) -- (0,4);

				\filldraw (1,4) circle (2pt);
				\draw[line width = 0.3pt] (1,3) -- (1,4);

				\filldraw (5,4) circle (2pt);
				\draw[line width = 0.3pt] (5,3) -- (5,4);

				\draw[line width = 0.3pt] (-0.5,3) -- (5.5,3);

				\draw[line width = 0.3pt] (5.5,3) -- (5.5,0);

				\draw[line width = 0.3pt] (5.5,0) -- (-0.5,0);
				
				\draw[line width = 0.3pt] (-0.5,0) -- (-0.5,3);  
				
				\draw[line width = 0.3pt] (0,0) -- (0,-1);
				\filldraw (0,-1) circle (2pt);
				\node at (0, -1.5) {$1$};

				\draw[line width = 0.3pt] (1,0) -- (1,-1);
				\filldraw (1,-1) circle (2pt);
				\node at (1, -1.5) {2};

				\draw[line width = 0.3pt] (5,0) -- (5,-1);
				\filldraw (5,-1) circle (2pt);
				\node at (5, -1.5) {n};

				\filldraw (6,4) circle (2pt);
				\draw[line width = 0.3pt] (6,4) -- (6,-1);
				\filldraw (6,-1) circle (2pt);
				\node at (6, -1.5) {n+1};

\end{scope}

				\end{tikzpicture}
				\caption{$y_{n}$}
			\end{figure}

			Let $y_{n}$ be the injection $B(B_{n}) \hookrightarrow B(B_{n+1})$.\\
		
			$B(A_{n})$ injects in $B(B_{n+1})$ by sending $\sigma_{i} $ to $\sigma_{i} $ for $1\leq i \leq n-1$, let us call this injection $z_{n}$. Take $T$ to be the normal subgroup in $B(B_{n+1})$ generated by $t$, that is the subgroup generated by $x t x^{-1}$ for $ x \in B(B_{n+1})$. Obviously $B(B_{n+1})/T = B(A_{n})$. in other words we have the following exact sequence $1 \rightarrow T \rightarrow B(B_{n+1}) \rightarrow B(A_{n}) \rightarrow 1 $. Call $ \alpha_{n} $ the surjection  $B(B_{n+1}) \twoheadrightarrow B(A_{n})$. \\ 

			Geometrically $B(A_{n})$ injects into $B(B_{n+1})$ by adding the first (fixed) strand, while $B(B_{n+1})$  surjects onto $B(A_{n})$ by removing the very same strand \\
			
			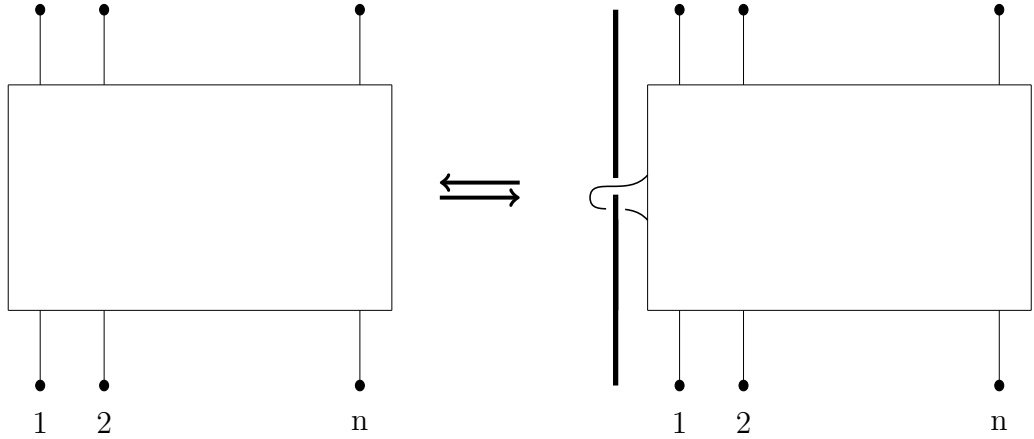
\begin{figure}[h]
				\centering
				\begin{tikzpicture}
				\begin{scope}[xscale = 0.85]

				\draw[draw= white, line width = 2pt] (-1,4)-- (-1,-1);

				\draw[draw= white, double = white, line width = 0.1cm] (-0.5, 1.2) .. controls +(-0.3, 0.3) and +(0, -0.3) .. (-1.4, 1.5) .. controls +(0,0.3) and +(-0.3,-0.3) .. (-0.5,1.8); 

				\filldraw[white] (-1,1.35) circle (4pt);
				\draw[draw= white, line width = 2pt] (-1,1.5)-- (-1,0);


			 \filldraw (0,4) circle (2pt);
				\draw[line width = 0.3pt] (0,3) -- (0,4);

				\filldraw (1,4) circle (2pt);
				\draw[line width = 0.3pt] (1,3) -- (1,4);

				\filldraw (5,4) circle (2pt);
				\draw[line width = 0.3pt] (5,3) -- (5,4);

				\draw[line width = 0.3pt] (-0.5,3) -- (5.5,3);

				\draw[line width = 0.3pt] (5.5,3) -- (5.5,0);

				\draw[line width = 0.3pt] (5.5,0) -- (-0.5,0);
				
				\draw[line width = 0.3pt] (-0.5,0) -- (-0.5,3);  
				
				\draw[line width = 0.3pt] (0,0) -- (0,-1);
				\filldraw (0,-1) circle (2pt);
				\node at (0, -1.5) {$1$};

				\draw[line width = 0.3pt] (1,0) -- (1,-1);
				\filldraw (1,-1) circle (2pt);
				\node at (1, -1.5) {2};

				\draw[line width = 0.3pt] (5,0) -- (5,-1);
				\filldraw (5,-1) circle (2pt);
				\node at (5, -1.5) {n};

				\draw[line width =1.5pt][->] (6.25,1.5) -- (7.5,1.5);
				\draw[line width =1.5pt][<-] (6.25,1.7) -- (7.5,1.7);

\end{scope}

\begin{scope}[xscale = 0.85, xshift = 10cm]

				\draw[line width = 2pt] (-1,4)-- (-1,-1);

				\draw[draw= white, double = black, line width = 0.1cm] (-0.5, 1.2) .. controls +(-0.3, 0.3) and +(0, -0.3) .. (-1.4, 1.5) .. controls +(0,0.3) and +(-0.3,-0.3) .. (-0.5,1.8); 

				\filldraw[white] (-1,1.35) circle (4pt);
				\draw[line width = 2pt] (-1,1.5)-- (-1,0);


			 \filldraw (0,4) circle (2pt);
				\draw[line width = 0.3pt] (0,3) -- (0,4);

				\filldraw (1,4) circle (2pt);
				\draw[line width = 0.3pt] (1,3) -- (1,4);

				\filldraw (5,4) circle (2pt);
				\draw[line width = 0.3pt] (5,3) -- (5,4);

				\draw[line width = 0.3pt] (-0.5,3) -- (5.5,3);

				\draw[line width = 0.3pt] (5.5,3) -- (5.5,0);

				\draw[line width = 0.3pt] (5.5,0) -- (-0.5,0);
				
				\draw[line width = 0.3pt] (-0.5,0) -- (-0.5,3);  
				
				\draw[line width = 0.3pt] (0,0) -- (0,-1);
				\filldraw (0,-1) circle (2pt);
				\node at (0, -1.5) {$1$};

				\draw[line width = 0.3pt] (1,0) -- (1,-1);
				\filldraw (1,-1) circle (2pt);
				\node at (1, -1.5) {2};

				\draw[line width = 0.3pt] (5,0) -- (5,-1);
				\filldraw (5,-1) circle (2pt);
				\node at (5, -1.5) {n};

\end{scope}

				\end{tikzpicture}
				\caption{$B(A_{n})\leftrightarrow B(B_{n+1})$}
			\end{figure}
			
		 
			We get the following commutative diagrams:\\ 

		  	\begin{figure}[ht]
		  		\centering
				

		\begin{tikzpicture}

			\matrix[matrix of math nodes,row sep=1cm,column sep=1cm]{
			|(A)| B(B_{n})         & & &     |(B)| B(B_{n+1}) & &     |(AA)| B(B_{n})         & & &     |(BB)| B(B_{n+1}) \\
			                       & & &                      & &                             & & &                      \\								
			|(C)| B(A_{n-1})       & & &     |(D)| B(A_{n})   & &     |(CC)| B(A_{n-1})       & & &     |(DD)| B(A_{n})   \\
				};
				
				\path (A) edge[-myhook,line width=0.42pt] node[above, xshift=-5mm, yshift=-2mm, rotate=0] {\footnotesize $z_{n-1}$}   (C);
				\path (C) edge[-myto,line width=0.42pt]      (A);
				
				\path (B) edge[-myhook,line width=0.42pt] node[above, xshift=1.5mm, yshift=0mm, rotate=0] {\footnotesize $y_{n}$}     (A);
				\path (A) edge[-myto,line width=0.42pt]      (B);
				
				\path (D) edge[-myhook,line width=0.42pt] node[below, xshift=1.5mm, yshift=0mm, rotate=0] {\footnotesize $x_{n}$}    (C);
				\path (C) edge[-myto,line width=0.42pt]      (D);
								
				\path (B) edge[-myhook,line width=0.42pt] node[above, xshift=3mm, , yshift=-2mm, rotate=0] {\footnotesize $z_{n}$}     (D);
				\path (D) edge[-myto,line width=0.42pt]      (B);

				\node at (-5, -2.5) {{$(1)$}};

				\path (AA) edge[-myonto,line width=0.42pt] node[above, xshift=-5mm, yshift=-2mm, rotate=0] {\footnotesize $\alpha_{n-1}$}    (CC);

				\path (BB) edge[-myhook,line width=0.42pt] node[above, xshift=-1.5mm, yshift=0mm, rotate=0] {\footnotesize $y_{n}$}   (AA);
				\path (AA) edge[-myto,line width=0.42pt]      (BB);
				
				\path (DD) edge[-myhook,line width=0.42pt] node[below, xshift=1.5mm, yshift=0mm, rotate=0] {\footnotesize $x_{n}$}    (CC);
				\path (CC) edge[-myto,line width=0.42pt]      (DD);
								
				\path (BB) edge[-myonto,line width=0.42pt] node[above, xshift=3mm, yshift=-2mm, rotate=0] {\footnotesize $\alpha_{n}$}    (DD);

				\node at (5, -2.5) {{$(2)$}};

		\end{tikzpicture}

			\end{figure}	
										
			\vspace{0.25cm}

			Diagram 1 commutes obviously, while for diagram  2 it is clear that			
			\begin{eqnarray}
				x_{n} \alpha_{n-1} (\sigma_{i}) &=& \alpha_{n}y_{n} (\sigma_{i}) \text{ for } 1 \leq i \leq n-1, \nonumber\\
				x_{n} \alpha_{n-1} (t) &=& 1 = \alpha_{n}y_{n} (t).\nonumber
			\end{eqnarray}
			
			The embedding of $B(\tilde{A_{n-1}})$ into $B(\tilde{A_{n}})$ is not as obvious as the other two embeddings. Consider the injection $y_{n}: B(B_{n}) \hookrightarrow B(B_{n+1})$ (see \cite{Geck_Lambropoulou_1997}), which sends $a_{n}$ to  $$ t\sigma_{1} .. \sigma_{n-2} \sigma_{n-1}\sigma^{-1}_{n-2} .. \sigma^{-1}_{1}t^{-1},$$  which is equal to			
			\begin{eqnarray}			
				t\sigma_{1} .. \sigma_{n-2} \sigma_{n-1} \underbrace{\sigma_{n}\sigma_{n-1}\sigma^{-1}_{n-1}\sigma^{-1}_{n}}_{1} \sigma^{-1}_{n-2} .. \sigma^{-1}_{1}t^{-1} &=& t \sigma_{1} .. \sigma_{n-2} \underbrace{\sigma_{n-1}\sigma_{n}\sigma_{n-1}}_{} \sigma^{-1}_{n-1}\sigma^{-1}_{n} \sigma^{-1}_{n-2} .. \sigma^{-1}_{1}t^{-1} \nonumber\\\nonumber\\
				&=& t \sigma_{1} .. \sigma_{n-2} \sigma_{n}\sigma_{n-1}\sigma_{n} \sigma^{-1}_{n-1} \sigma^{-1}_{n} \sigma^{-1}_{n-2} .. \sigma^{-1}_{1}t^{-1}  \nonumber\\\nonumber\\
				&=& \sigma_{n}\underbrace{t \sigma_{1} ..  \sigma_{n-1}\sigma_{n} \sigma^{-1}_{n-1}  .. \sigma^{-1}_{1}t^{-1}}_{a_{n+1}}\sigma^{-1}_{n}.  \nonumber\\\nonumber
			\end{eqnarray}			
			
			In other terms $y_{n}(a_{n}) = a_{n+1} $ which is in $B(\tilde{A_{n}})$, in other terms the restriction of $y_{n}$ to $B(\tilde{A_{n-1}})$ is equal to $ F_{n}$ as defined in Proposition \ref{1_4_5}, thus, this restriction is injective:

			\begin{eqnarray}
				F_{n}: B(\tilde{A_{n-1}}) &\longrightarrow& B(\tilde{A_{n}}) \nonumber\\        
         			\sigma_{i} &\longmapsto& \sigma_{i} \text{ for } 1\leq i\leq n-1\nonumber\\
        			a_{n} &\longmapsto& \sigma_{n} a_{n+1}\sigma^{-1}_{n}.\nonumber
			\end{eqnarray}	
			
			Set $ I_{n}$ to be the injection $ B(A_{n}) \hookrightarrow B(\tilde{A_{n}})$. Set $ i_{n}$ to be  the injection $ B(\tilde{A_{n}}) \hookrightarrow B(B_{n+1})$. We have the following commutative diagram: \\
		
			\vspace{0.5cm} 
		
		  	\begin{figure}[ht]
		  		\centering
				

		\begin{tikzpicture}

			\matrix[matrix of math nodes,row sep=1cm,column sep=1cm]{
			   |(A)| B(B_{n})    &                          & &                        & |(B)| B(B_{n+1})   \\
			                     &                          & &                        &                    \\								
			                     & |(C)| B(\tilde{A_{n-1}}) & & |(D)| B(\tilde{A_{n}}) &                    \\
				                   &                          & &                        &                    \\								
			   |(E)| B(A_{n-1})  &                          & &                        & |(F)| B(A_{n-1})   \\	
				};
				
				\path (A) edge[-myhook,line width=0.42pt] node[above, xshift=-5mm, yshift=-2mm, rotate=0] {\footnotesize $i_{n-1}$}   (C);
				\path (C) edge[-myto,line width=0.42pt]      (A);
				
				\path (B) edge[-myhook,line width=0.42pt] node[above, xshift=2mm, yshift=0mm, rotate=0] {\footnotesize $y_{n}$}     (A);
				\path (A) edge[-myto,line width=0.42pt]      (B);
				
				\path (D) edge[-myhook,line width=0.42pt] node[above, xshift=1mm, yshift=0mm, rotate=0] {\footnotesize $F_{n}$}    (C);
				\path (C) edge[-myto,line width=0.42pt]      (D);
								
				\path (B) edge[-myhook,line width=0.42pt] node[above, xshift=4mm, , yshift=-2mm, rotate=0] {\footnotesize $i_{n}$}     (D);
				\path (D) edge[-myto,line width=0.42pt]      (B);

				\path (C) edge[-myhook,line width=0.42pt] node[above, xshift=-5mm, yshift=-2mm, rotate=0] {\footnotesize $I_{n-1}$}   (E);
				\path (E) edge[-myto,line width=0.42pt]      (C);

				\path (D) edge[-myhook,line width=0.42pt] node[above, xshift=3mm, yshift=-2mm, rotate=0] {\footnotesize $I_{n}$}   (F);
				\path (F) edge[-myto,line width=0.42pt]      (D);

				\path (F) edge[-myhook,line width=0.42pt] node[above, xshift=2mm, yshift=0mm, rotate=0] {\footnotesize $x_{n}$}   (E);
				\path (E) edge[-myto,line width=0.42pt]      (F);
				
				\path (B) edge[-myhook,line width=0.42pt] node[above, xshift=3mm, yshift=0mm, rotate=0] {\footnotesize $z_{n}$}   (F);
				\path (F) edge[-myto,line width=0.42pt]      (B);
				
				\path (A) edge[-myhook,line width=0.42pt] node[above, xshift=-5mm, yshift=0mm, rotate=0] {\footnotesize $z_{n-1}$}   (E);
				\path (E) edge[-myto,line width=0.42pt]      (A);

		\end{tikzpicture}

			\end{figure}
			
			\vspace{0.5cm}

			Geometrically, we realize in $ B(\tilde{A_{n}}) $ the generators $\sigma_{i} $ of $B(\tilde{A_{n-1}})$ (for $ 1 \leq i \leq n-1$) in the natural way, while concerning $a_{n}$ we see that
								
		  	\begin{figure}[h]
		  		\centering
					\begin{tikzpicture}
					
 \begin{scope}[scale = 0.5]

				\draw[line width = 2pt] (0,2.0)-- +(0,-6);
				\filldraw[white] (0,0.05) circle (4pt);
				\filldraw[white] (0,-2.0) circle (4pt);

				\draw[double = black, line width = 0.1pt, dashed] (-0.5,0.75) -- +(10,0);
				\draw[double = black, line width = 0.1pt, dashed] (-0.5,-2.5) -- +(10,0);
				
				\draw[line width = 0.3pt] (1.5,2.0) -- +(0,-0.5);
				\draw[line width = 0.3pt] (1.5,1.5).. controls +(0, -0.5) and +(0, 0.5) .. (-0.5,0.25).. controls +(0,-0.6) and +(-0.4,1.2) .. (8.5,-1).. controls +(-0,-1)  and +(0,0.5) .. (8.5, -1.75) .. controls +(0.1,-2)  and +(0,2) .. (7, -4);
				\filldraw[white] (0,0.7) circle (4pt); 

				\filldraw (1.5,2.0) circle (2pt);
				\filldraw[white] (8.3,-0.8) circle (4pt);				
				
				\draw[line width = 0.3pt] (1.5,-3) -- +(0,-1);
				\draw[line width = 0.3pt] (1.5,-3) .. controls +(0,0.5) and +(0, -0.5) .. (-0.50,-1.5).. controls +(0,0.7) and +(-0.6,-1.5) .. (8.5,-0.5) .. controls +(0,0.3) and +(0,-0.2) .. (8.5,0.4) .. controls +(-0.3,1)  and +(0,-1) .. (7,2);

				\filldraw[white] (0,-1.3) circle (4pt);

				\filldraw (1.5,-4) circle (2pt);

				\draw[draw= white, double = black, line width = 2pt, double distance= 0.3pt] (3,2.0) -- +(0,-6);

				\draw[draw= white, double = black, line width = 2pt, double distance= 0.3pt] (8.5,2)-- +(0,0) .. controls +(0,-1) and +(0,0.75) .. (7,0.5) .. controls +(0,0) and +(0,0.5) .. (7,-2) .. controls +(0,-1.5) and +(0,1.5) .. (8.5,-4);

				\node at (1.5, -4.5) {$1$};
				\node at (3, -4.5) {$2$};
				\node at (7, -4.5) {$n$};
				\node at (8.5, -4.5) {$n+1$};
				
				\node at (10, 1.5) {$\sigma_{n}$};			
				\node at (10, -3) {$\sigma_{n}^{-1}$};
				\node at (-1.5, -0.5) {$a_{n+1}$};

				\filldraw (3,-4) circle (2pt);
				\filldraw (3,2.0) circle (2pt);
				
				\filldraw (7,-4) circle (2pt);
				\filldraw (7,2.0) circle (2pt);

				\draw[line width =1.5pt][<->] (11,-0.7) -- (12.5,-0.7);
				\filldraw (8.5,2.0) circle (2pt);
				\filldraw (8.5,-4) circle (2pt);
				\node at (5, -5.5) {$(\sigma_{n}a_{n+1}\sigma_{n}^{-1})$};

				\draw[line width = 2pt] (0,1.2)-- +(0,-1);
				\draw[line width = 2pt] (0,-0.7)-- +(0,-1);
\end{scope}

 \begin{scope}[scale = 0.5, xshift = 14cm]

  \draw[line width = 2pt] (0,2.0)-- +(0,-6);

				\filldraw[white] (0,1.4) circle (4pt);  
 				\filldraw[white] (0,-2.65) circle (4pt);

				\draw[line width = 0.3pt] (1.5,2).. controls +(0, -0.5) and +(0, 0.5) .. (-0.5,1).. controls +(0,-0.6) and +(-0.4,1.2) .. (8.25,0).. controls +(0.3,-0.6)  and +(0,1) .. (8.5, -4);
				\filldraw (1.5,2.0) circle (2pt);
				\filldraw (8.5,-4) circle (2pt);
				\filldraw[white] (8.3,0) circle (4pt); 
				\draw[line width = 0.3pt] (1.5,-3.5) -- +(0,-0.5);
				\draw[line width = 0.3pt] (1.5,-3.5).. controls +(0, 0.5) and +(0, -0.5) .. (-0.50,-2.2).. controls +(0,0.6) and +(-0.4,-2.0) .. (8.25,0) .. controls +(0.3,0.6)  and +(0,-1) .. (8.5, 2.0);
				\filldraw (1.5,-4) circle (2pt);
				\filldraw (8.5,2.0) circle (2pt);
				
				\filldraw[white] (0,0.85) circle (4pt); 
				
				\filldraw[white] (0,-2) circle (4pt);

				\filldraw[white] (3,0.75) circle (4pt);
				\filldraw[white] (3,-1.7) circle (4pt);

				\draw[draw= white, double = black, line width = 0.3pt, double distance= 2pt] (0,1)-- +(0,-3.4);
				\draw[line width = 0.3pt] (3,2.0) -- +(0,-6);
				
				\node at (5,2.0) {$\dots$};
				\node at (5,-4) {$\dots$};
				
				\draw[line width = 0.3pt] (10,2.0) -- +(0,-6);

				\node at (1.5, -4.5) {$1$};
				\node at (3, -4.5) {$2$};
				\node at (8.5, -4.5) {$n$};
				\node at (10, -4.5) {$n+1$};

				\filldraw (3,-4) circle (2pt);
				\filldraw (3,2.0) circle (2pt);
				
				\filldraw (10,-4) circle (2pt);
				\filldraw (10,2.0) circle (2pt);

				\node at (5, -5.5) {$(a_{n})$};

\end{scope}

						\end{tikzpicture}
				\caption{$F_{n}$}
			\end{figure}								
			
			Recall that $B(B_{n})$ is generated by $ \left\{ \sigma_{1}, \dots,  \sigma_{n}, a_{n}, \phi_{n}  \right\}$. Now we consider $\phi_{n} $ as an automorphism of $ B(\tilde{A_{n-1}}) $. We call $\phi_{n} $ the Dynkin automorphism of order $n$, since it generates a subgroup of $Aut( B(\tilde{A_{n-1}}) )$ of order $n$. It shifts the generators of the Dynkin diagram one step counter clockwise ($ \sigma_{1}\mapsto  \sigma_{2} \mapsto  \dots  \mapsto  \sigma_{n-1} \mapsto a_{n} \mapsto \sigma_{1} $). In order to simplify we call it the  Dynkin automorphism, referring to it by  $\phi$ when there is no ambiguity.\\ 
						
			On the other hand we see that in $ B(\tilde{A_{n}}) $ for $ 2 \leq i \leq n-1 $: 
			\begin{eqnarray}
				\sigma_{n} \sigma_{n-1}.. \sigma_{1}a_{n+1} \sigma_{i} &=& \sigma_{i-1} \sigma_{n}\sigma_{n-1}.. \sigma_{1}a_{n+1},\nonumber\\
				\sigma_{n} \sigma_{n-1}.. \sigma_{1}a_{n+1} \sigma_{1} &=& a_{n} \sigma_{n} \sigma_{n-1}.. \sigma_{1}a_{n+1}, \nonumber\\
				\sigma_{n} \sigma_{n-1}.. \sigma_{1}a_{n+1} a_{n} &=& \sigma_{n-1} \sigma_{n} \sigma_{n-1}.. \sigma_{1}a_{n+1}. \nonumber
			\end{eqnarray}
					
			The last equality comes from the fact that $ a_{n+1} a_{n}=  \sigma^{-1}_{n} a_{n}\sigma_{n}a_{n} = \sigma^{-1}_{n} \sigma_{n} a_{n} \sigma_{n}$, which is equal to $ a_{n} \sigma_{n} =  \sigma_{n}a_{n+1} $. Hence $ \sigma_{n} \sigma_{n-1}.. \sigma_{1}a_{n+1} $ acts on the elements of $ B(\tilde{A_{n-1}}) $ exactly the way as $\phi^{-1}_{n}$ does in $B(B_{n})$. \\
		
			\begin{definition}
				In $ B(\tilde{A_{n-1}}) $ we call $ \sigma_{n} \sigma_{n-1} \dots \sigma_{1}a_{n+1} $ the dominating element, and we denote it by $D_{n}$. When there is no ambiguity we call it $ D $.\\ 
			\end{definition}

			The following diagram commutes\\
			
		  	\begin{figure}[ht]
		  		\centering
				

		\begin{tikzpicture}

			\matrix[matrix of math nodes,row sep=1cm,column sep=1cm]{
			|(A)| B(B_{n})                 & & & &                              \\
			                               & & & &                              \\
			                               & & & &                              \\
			|(C)| B(\tilde{A_{n-1}})       & & & &    |(D)| B(\tilde{A_{n}})    \\
				};
				
				\path (A) edge[-myhook,line width=0.42pt] node[above, xshift=-5mm, yshift=-2mm, rotate=0] {\footnotesize $i_{n-1}$}   (C);
				\path (C) edge[-myto,line width=0.42pt]      (A);

				\path (D) edge[-myhook,line width=0.42pt] node[below, xshift=1.5mm, yshift=0mm, rotate=0] {\footnotesize $F_{n}$}    (C);
				\path (C) edge[-myto,line width=0.42pt]      (D);

				\path (A) edge[-myto,line width=0.42pt] node[above, xshift=3mm, , yshift=-2mm, rotate=0] {\footnotesize $f_{n}$}          (D);

		\end{tikzpicture}

			\end{figure}
			
			\vspace{-0.5cm}			
			
			\begin{eqnarray}
				\text{where } f_{n}: B(B_{n}) &\longrightarrow& B(\tilde{A_{n}}) \nonumber\\
				\sigma_{i} &\longmapsto& \sigma_{i}$ for $1\leq i\leq n-1\nonumber\\
				a_{n} &\longmapsto& \sigma_{n} a_{n+1}\sigma^{-1}_{n} \nonumber\\				
				\phi_{n} &\longmapsto& D^{-1}_{n}.\nonumber\\\nonumber
			\end{eqnarray}
			
			We consider the group $B(\tilde{A_{n}})$ modulo the action of $\left\langle  \phi_{n+1} \right\rangle_{Aut(B(\tilde{A_{n}}))}$ (the subgroup of $Aut(B(\tilde{A_{n}}))$ of order $n+1$ generated by the Dynkin automorphism). This group is isomorphic to the free group in one letter.\\ 
		
		We have seen that $B(\tilde{A_{n}})$ surjects onto $B(A_{n})$. Call this surjection $\beta_{n}$, we get the following diagram: 
		
			
		  	
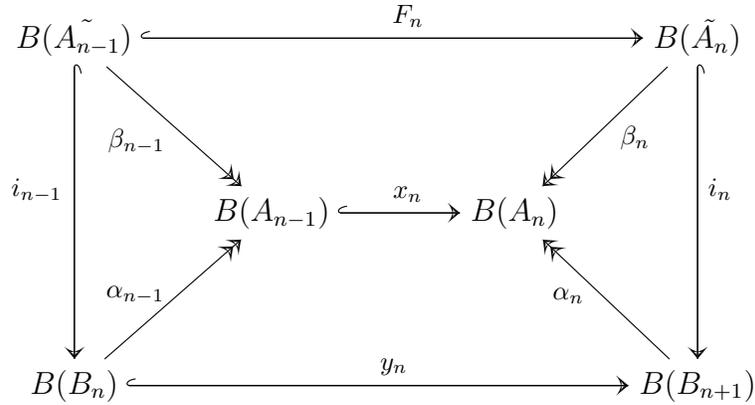
\begin{figure}[ht]
		  		\centering
				

		\begin{tikzpicture}

			\matrix[matrix of math nodes,row sep=0.8cm,column sep=0.8cm]{
			|(A)| B(\tilde{A_{n-1}})   &                   & &                    &    |(B)| B(\tilde{A_{n}})   \\
			                           &                   & &                    &                             \\								
			                           & |(C)|  B(A_{n-1}) & &    |(D)| B(A_{n})  &                             \\
				                      &                   & &                    &                             \\								
			|(E)| B(B_{n})             &                   & &                    &    |(F)| B(B_{n+1})         \\	
				};

				\path (B) edge[-myhook,line width=0.42pt] node[above, xshift=2mm, yshift=0mm, rotate=0] {\footnotesize $F_{n}$}     (A);
				\path (A) edge[-myto,line width=0.42pt]      (B);
				
				\path (D) edge[-myhook,line width=0.42pt] node[above, xshift=1mm, yshift=0mm, rotate=0] {\footnotesize $x_{n}$}    (C);
				\path (C) edge[-myto,line width=0.42pt]      (D);

				\path (F) edge[-myhook,line width=0.42pt] node[above, xshift=2mm, yshift=0mm, rotate=0] {\footnotesize $y_{n}$}   (E);
				\path (E) edge[-myto,line width=0.42pt]      (F);

				\path (A) edge[-myonto,line width=0.42pt] node[above, xshift=-5mm, yshift=-5mm, rotate=0] {\footnotesize $\beta_{n-1}$}     (C);
				
				\path (E) edge[-myonto,line width=0.42pt] node[above, xshift=-5mm, yshift=-2mm, rotate=0] {\footnotesize $\alpha_{n-1}$}     (C);

				\path (E) edge[-myhook,line width=0.42pt] node[above, xshift=-5mm, yshift=0mm, rotate=0] {\footnotesize $i_{n-1}$}   (A);
				\path (A) edge[-myto,line width=0.42pt]      (E);
				
				\path (F) edge[-myhook,line width=0.42pt] node[below, xshift=3mm, yshift=6mm, rotate=0] {\footnotesize $i_{n}$}   (B);
				\path (B) edge[-myto,line width=0.42pt]      (F);
				
				\path (B) edge[-myonto,line width=0.42pt] node[above, xshift=4mm, yshift=-4.5mm, rotate=0] {\footnotesize $\beta_{n}$}     (D);
				
				\path (F) edge[-myonto,line width=0.42pt] node[above, xshift=-5mm, yshift=-2mm, rotate=0] {\footnotesize $\alpha_{n}$}     (D);

		\end{tikzpicture}

			\end{figure}
			
			We have $\alpha_{n-1} i_{n-1}(a_{n-1}) =  \alpha_{n-1}(t\sigma_{1} .. \sigma_{n-2} \sigma_{n-1}\sigma^{-1}_{n-2} .. \sigma^{-1}_{1}t^{-1})  = \sigma_{1} .. \sigma_{n-2} \sigma_{n-1}\sigma^{-1}_{n-2} .. \sigma^{-1}_{1}$, which is equal to $\beta_{n-1} (a_{n-1})$. Hence $\alpha_{n-1} i_{n-1}=\beta_{n-1}$. But we know already that $i_{n} F_{n} =y_{n} i_{n-1}$ and $x_{n} \alpha_{n-1}= \alpha_{n} y_{n}$. Hence $x_{n} \beta_{n-1} = \beta_{n} F_{n}$.\\
		
			Finally we present the arrows shown earlier by the following two diagrams: 	\\
					
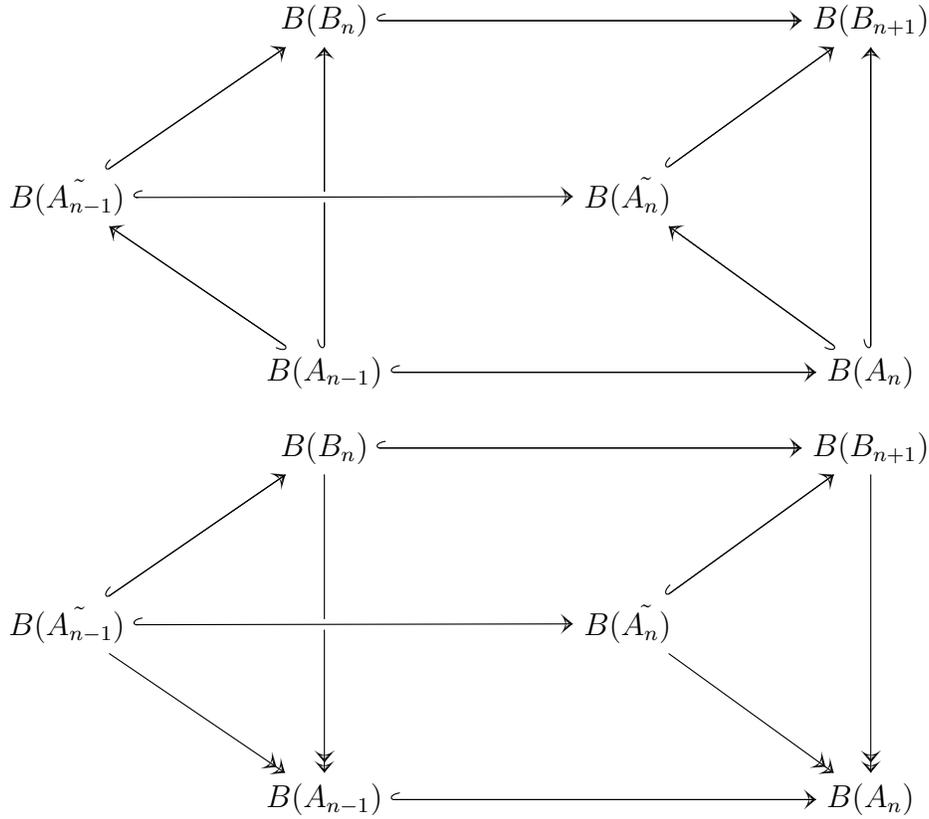
		  	\begin{center}
				

		\begin{tikzpicture}

			\matrix[matrix of math nodes,row sep=0.8cm,column sep=0.8cm]{
				                               & & |(A)| B(B_{n})   & & &                         & & |(B)| B(B_{n+1})     \\
				                               & &                  & & &                         & &                      \\								
				      |(C)| B(\tilde{A_{n-1}}) & &                  & & &  |(D)| B(\tilde{A_{n})} & &                      \\
  			                                    & &                  & & &                         & &                      \\								
				                               & & |(E)| B(A_{n-1}) & & &                         & & |(F)| B(A_{n})       \\
				};

				\path (B) edge[-myhook,line width=0.42pt]                                  (A);
				\path (A) edge[-myto,line width=0.42pt]                                    (B);
				
				\path (B) edge[-myhook,line width=0.42pt,transform canvas={xshift=0mm}]    (F);
				\path (F) edge[-myto,line width=0.42pt,transform canvas={xshift=0mm}]      (B);
				
				
				\path (F) edge[-myhook,line width=0.42pt]    (E);
				\path (E) edge[-myto,line width=0.42pt]      (F);
				
				\path (A) edge[-myhook,line width=0.42pt,transform canvas={xshift=0mm}]    (E);
				\path (E) edge[-myto,line width=0.42pt,transform canvas={xshift=0mm}]      (A);

				\path (D) edge[-myhook,line width=0.42pt]    (C);
				\path (D) edge[-myto,white, line width=4pt]  (C);
				\path (C) edge[-myto,line width=0.42pt]      (D);

				\path (A) edge[-myhook,line width=0.42pt]    (C);
				\path (C) edge[-myto,line width=0.42pt]      (A);

				\path (B) edge[-myhook,line width=0.42pt]    (D);
				\path (D) edge[-myto,line width=0.42pt]      (B);
				
				\path (D) edge[-myhook,line width=0.42pt,transform canvas={xshift= 0mm,yshift= 0mm}]    (F);
				\path (F) edge[-myto,line width=0.42pt,transform canvas={xshift= 0mm,yshift= 0mm}]      (D);

				\path (C) edge[-myhook,line width=0.42pt,transform canvas={xshift= 0mm,yshift= 0mm}]    (E);
				\path (E) edge[-myto,line width=0.42pt,transform canvas={xshift= 0mm,yshift= 0mm}]      (C);


		\end{tikzpicture}

				

		\begin{tikzpicture}

			\matrix[matrix of math nodes,row sep=0.8cm,column sep=0.8cm]{
				                               & & |(A)| B(B_{n})   & & &                         & & |(B)| B(B_{n+1})     \\
				                               & &                  & & &                         & &                      \\								
				      |(C)| B(\tilde{A_{n-1}}) & &                  & & &  |(D)| B(\tilde{A_{n})} & &                      \\
  			                                    & &                  & & &                         & &                      \\								
				                               & & |(E)| B(A_{n-1}) & & &                         & & |(F)| B(A_{n})       \\
				};

				\path (B) edge[-myhook,line width=0.42pt]                                  (A);
				\path (A) edge[-myto,line width=0.42pt]                                    (B);
				
				\path (B) edge[-myonto,line width=0.42pt,transform canvas={xshift=0mm}]    (F);
				
				
				\path (F) edge[-myhook,line width=0.42pt]    (E);
				\path (E) edge[-myto,line width=0.42pt]      (F);
				
				\path (A) edge[-myonto,line width=0.42pt,transform canvas={xshift=0mm}]    (E);

				\path (D) edge[-myhook,line width=0.42pt]    (C);
				\path (D) edge[-myto,white, line width=4pt]  (C);
				\path (C) edge[-myto,line width=0.42pt]      (D);

				\path (A) edge[-myhook,line width=0.42pt]    (C);
				\path (C) edge[-myto,line width=0.42pt]      (A);

				\path (B) edge[-myhook,line width=0.42pt]    (D);
				\path (D) edge[-myto,line width=0.42pt]      (B);
				
				\path (D) edge[-myonto,line width=0.42pt,transform canvas={xshift= 0mm,yshift= 0mm}]    (F);

				\path (C) edge[-myonto,line width=0.42pt,transform canvas={xshift= 0mm,yshift= 0mm}]    (E);


		\end{tikzpicture}

			\end{center}
			

	     \subsection{Affine links and closures} \label{subsec1_4_5}
	     
	     $  \    $   
			\vspace{0.5cm}  
			
			In what follows we give some definitions and basic results in the theory of links, the classical well known results concerning invariants of links are to be mentioned briefly. The aim is to define the concept of  "affine" links, defining dual concepts and conventions to those in the classical theory. Most of theorems and results here are well explained in the literature, and we will not give details.\\ 
			
	   		Let $C_{i}$ be a circle in $\mathds{R}^{3}$, where $ 1 \leq i \leq n $. Let $C^{n}:=\cup_{i}C_{i}$ be the disjoint union of those $n$ circles. We call $C^{n}$ a rough link. Take the isotopy class of $C^{n}$, call it $C$, we call $C$ a circle link or simply a link. We consider here only the piecewise  linear links. If we orient the circles forming $C$, we say that it is an oriented link. Notice that orienting a circle is independent of orienting another one, since they do not intersect. Inverting the orientation of one of the circles gives a different oriented link.\\ 
				
			Roughly speaking, the problem of finding an invariant for the set of links in $\mathds{R}^{3}$ is to give names to links, in such way that any two links which have the same "shape" (i.e., we can arrive to one from the other by pulling and pushing the circles forming a link without cutting, adding or omitting any of the circles) have the same name.\\ 
				
			We recall the results of Alexander and Markov concerning braids and links.\\
					
			Consider a braid $b$, that is: a an element of $B(A_{n})$ for some $1 \leq n $. As we have $ n+1$ strands with $ n+1$ points at the top (the same at the bottom), a path from a point (say the $i$-th point) at the top to the $i$-th at the bottom makes the $i$-th strand turn into a deformation of a circle in $\mathds{R}^{3}$, repeating the same step with all the points (using non crossing paths) gives a union of disjoint deformed circles, hence a link. Thus we have defined a mapping from $\bigcup_{1 \leq i} B(A_{i}) $ into the set of links in  $\mathds{R}^{3}$. We call the image of $b$: the closure of $b$, denoted by~$\hat{b}$.	
				
			 \AddToShipoutPicture*{\POSITION{4cm}{13cm}{
			 \begin{tikzpicture}
			 \begin{scope}[xscale = 1.5]
\draw (0,0) -- +(0,3);
\draw (1,0) -- +(0,3);
\draw (3.5,0) -- +(0,3);
\filldraw[fill= white] (-0.5,0.5)-- (4, 0.5) -- (4, 2.5) -- (-0.5, 2.5) -- cycle; 
\draw (3.5,0).. controls +(0,-0.5) and +(0,-1.5) .. (4.5,1.5).. controls +(0,1.5) and +(0, 0.5) .. (3.5,3);
\draw (1,0).. controls +(0,-0.6) and +(0,-2.7) .. (5,1.5).. controls +(0,2.7) and 
+(0, 0.5) .. (1,3);
\draw (0,0).. controls +(0,-0.7) and +(0,-3) .. (5.5,1.5).. controls +(0,3) and 
+(0, 0.5) .. (0,3);
\end{scope}

			 \end{tikzpicture}
			 }}
			\vspace{4.5cm}
		  	\begin{figure}[h]
		  		\centering
		  		\parbox{3cm}{\textcolor{white}{...................................}}\parbox{5cm}{\caption{ }}
			\end{figure}		
				
				\vspace{1cm} 
				
			We are interested in a map in the opposite direction. \\
				
			\begin{theoreme}
				(Alexander) Suppose that $C$ is an oriented link, then there exists an integer $1 \leq n$ and $b$ in $B(A_{n})$ such that $\hat{b} = C$.\\ 
			\end{theoreme}
			
			In other terms, set $OK$ to be the set of all oriented links in $\mathds{R}^{3}$. Then the following map is surjective:  
			\begin{eqnarray}
				\bigcup_{1 \leq i} B(A_{i}) &\longrightarrow& OK  \nonumber\\
				b &\longmapsto& \hat{b} .\nonumber
			\end{eqnarray}
			
			Let us present the main idea of the proof. Take any link $C$. Suppose it is the union of $n$ deformed circles. We project it on  $\mathds{R}^{2}$ respecting the crossing points (the concept of positive and negative crossing makes it doable). We take any point in $\mathds{R}^{2}$, say $P$. We take an orientation of every circle (by arbitrary orientation of those circles we get all the possible orientations of $C$). The point $P$ defines negative and positive rounds, say that the negative, for example, can be shifted to the 'right' of $P$ the positive on the 'left', hence we arrive to some presentation of $C$ as the following 			
				\clearpage		
		  	\begin{figure}[h]
		  		\centering
				 \begin{tikzpicture}
				 \begin{scope}[xscale = 1.5]
\draw (0,0) -- +(0,3);
\draw (1,0) -- +(0,3);
\draw (3.5,0) -- +(0,3);
\filldraw[fill= white] (-0.5,0.5)-- (4, 0.5) -- (4, 2.5) -- (-0.5, 2.5) -- cycle; 
\draw (3.5,0).. controls +(0,-0.5) and +(0,-1.5) .. (4.5,1.5).. controls +(0,1.5) and +(0, 0.5) .. (3.5,3);
\draw (1,0).. controls +(0,-0.6) and +(0,-2.7) .. (5,1.5).. controls +(0,2.7) and 
+(0, 0.5) .. (1,3);
\draw (0,0).. controls +(0,-0.7) and +(0,-3) .. (5.5,1.5).. controls +(0,3) and 
+(0, 0.5) .. (0,3);

  \draw (4,1.5) -- (6.5, 1.5);
	\node at (7,1.5) {$\rho$};
	
\end{scope}

				  \end{tikzpicture}
				\caption{ }
			\end{figure}

			Then we can cut the circle on the axis $\rho$, getting a braid whose closure is $C$.\\
		
			 The next question is, obviously: when do two braids have the same closure? The answer is a theorem of Markov -- actually it was  	announced by Markov himself, finally proven by Birman.\\

			\begin{theoreme} \label{1_4_13}
				(Markov) Two braids have the same closure if and only if there  exists an integer $ 1\leq n$ such that: starting from one of the two braids we can arrive to the other by a finite number of transformations of the two following types:\\
				
				\begin{itemize}[label=$\bullet$, font=\normalsize, font=\color{black}, leftmargin=2cm, parsep=0cm, itemsep=0.25cm, topsep=0cm]
					 \item $ab \leftrightarrow ba  $ where $a$ and $b$ are in $B(A_{n})$,
					 \item $x \leftrightarrow x\sigma_{n}  $ where $x$  is  in $B(A_{n-1})$.\\
				\end{itemize} 
			\end{theoreme}

			After this theorem a very elegant answer would be to find a family of applications $t_{n+1}$ defined over $B(A_{n})$  such that for all $1 \leq n $: \\ 
			
			\begin{itemize}[label=$\bullet$, font=\normalsize, font=\color{black}, leftmargin=2cm, parsep=0cm, itemsep=0.25cm, topsep=0cm]				
				 \item $t_{n}(ab) = t_{n}(ba)$  for all $a$ and $b$   in $B(A_{n-1})$,
				 \item $t_{n+1}(x\sigma_{n}) = t_{n}(x) $ for any  $x$     in $B(A_{n-1})$.\\
			\end{itemize} 
			
			The answer given by Jones was exactly of this form. Here we attempt to define an "affine link" as a result of closing an affine braid, where we mean by affine braid an element of an $\tilde{A}$-type braid group. This task is not as evident as for the $A$-type braid group, for we have many geometrical presentations of  $\tilde{A}$-type braid groups. A choice must be made here, this is what we are about to do in the rest of this section. \\ 
			
			In order to simplify we call an oriented link simply a link in $\mathds{R}^{3}$ (in the literature $S^{3}$ is often used). \\ 
			
			We call a $B$-braid any element in a given $B$-type braid. Clearly any affine braid is a $B$-braid, which has a presentation as a cylindrical braid which could be closed at least in two ways. Here we view $B$-braids as we did above, braids with one fixed strand. Now we consider the following application:\\  
			\begin{eqnarray}
				_{n}I:B(B_{n}) &\longrightarrow& B(A_{n}) \nonumber\\
				\sigma_{i} &\mapsto& \sigma_{i+1}$ for $1\leq i\leq n-1 \nonumber\\
				t &\mapsto& \sigma^{2}_{1}. \nonumber
			\end{eqnarray} 
			
			\vspace{0.5cm} 
			
			It is geometrically presented as follows:\\  
			 
		  	\begin{figure}[h]
		  		\centering
				 \begin{tikzpicture}
				 \begin{scope}[xscale = 1.5, yscale = 1.75]
\draw[draw= white, double = black, line width = 0.3pt, double distance = 2pt] (0,2) -- (0,1);
\draw[draw= white, double = black, line width = 2pt, double distance = 0.4pt] (1, 0) .. controls +(0, 0.5) and +(0, -0.5) .. (-1, 1) .. controls +(0,0.5) and +(0,-0.5) .. (1,2); 
\draw[draw= white, double = black, line width = 3pt, double distance = 2pt] (0, 1) -- (0,0);
\node at (1, -0.5) {$1$};
\filldraw (1,0) circle (1.5pt);
\filldraw (1,2) circle (1.5pt);
\end{scope}
\path (4,2) edge[-myto,line width=0.42pt] (5,2);
\begin{scope}[xscale = 1.5, yscale = 1.75, xshift = 5cm]
\draw[draw= white, double = black, line width = 2pt, double distance = 0.4pt] (0, 1) .. controls +(0, -0.5) and +(0, +0.5) .. (1, 0);
\draw[draw= white, double = black, line width = 2pt, double distance = 0.4pt] (0,0) .. controls +(0, 0.5) and +(0, -0.5) .. ( 1,1);
\node at (0, -0.5) {$1$};
\filldraw (0,0) circle (1.5pt);
\node at (1, -0.5) {$2$};
\filldraw (1,0) circle (1.5pt);
\end{scope}
\begin{scope}[xscale = 1.5, yscale = 1.75, xshift = 5cm, yshift=1cm]
\draw[draw= white, double = black, line width = 2pt, double distance = 0.4pt] (0, 1) .. controls +(0, -0.5) and +(0, +0.5) .. (1, 0);
\draw[draw= white, double = black, line width = 2pt, double distance = 0.4pt] (0,0) .. controls +(0, 0.5) and +(0, -0.5) .. ( 1,1);
\filldraw (0,1) circle (1.5pt);
\filldraw (1,1) circle (1.5pt);


\end{scope}

				  \end{tikzpicture}
				\caption{$t\longmapsto\sigma_{1}^{2}$}
			\end{figure}
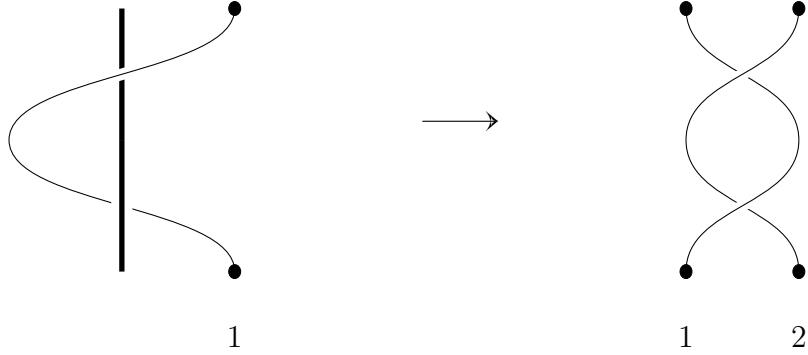		

			\vspace{0.5cm}
			 
			We see that
			\begin{eqnarray}
				_{n}I(t\sigma_{1}t\sigma_{1} ) &=& \sigma^{2}_{1} \sigma_{2}  \sigma^{2}_{1} \sigma_{2} =  \underbrace{\sigma_{1} \sigma_{2}\sigma_{1}}_{} \underbrace{\sigma_{2} \sigma_{1} \sigma_{2} }_{}  \nonumber\\
				&=& \sigma_{2} \sigma_{1}\underbrace{\sigma_{2} \sigma_{1} \sigma_{2} }_{}\sigma_{1}  = \sigma_{2}\sigma^{2}_{1} \sigma_{2} \sigma^{2}_{1}= _{n}I (\sigma_{1}t\sigma_{1}t). \nonumber\\\nonumber 
			\end{eqnarray}
			
			In other terms $_{n}I$ is a homomorphism. Moreover, it is a monomorphism (see \cite{Crisp_1999}) . The following  diagram, of injections, is  commutative:
  \clearpage
		 \begin{figure}[h]
		  		 \centering
				\begin{tikzpicture}

			\matrix[matrix of math nodes,row sep=1cm,column sep=1cm]{
			|(A)| B(B_{n})                 & & & & &   |(B)| B(A_{n})    \\
			                               & & & & &                     \\		
			                               & & & & &                     \\
			|(C)| B(\tilde A_{n-1})               & & & & &                     \\
				};
				
				\path (A) edge[-myhook,line width=0.42pt] node[above, xshift=-1mm, yshift=0mm, rotate=90] {\footnotesize $i_{n-1}$}  node[below, xshift=1mm, yshift=0mm, rotate=90] {  }  (C);
				\path (C) edge[-myto,line width=0.42pt]      (A);
				
				\path (B) edge[-myhook,line width=0.42pt]  node[above, xshift=-1mm, yshift=0mm] {\footnotesize $\ _n I $}  (A);
				\path (A) edge[-myto,line width=0.42pt]      (B);

				\path (B) edge[-myhook,line width=0.42pt] node[above, xshift=0.5mm, yshift=0mm, rotate=40] {\footnotesize $\bar{x_{n}} $} node[below, xshift=1.5mm, yshift=0mm, rotate=35] { }    (C);
				\path (C) edge[-myto,line width=0.42pt]      (B);

 		\end{tikzpicture}
		\caption{ }
	 			 \end{figure}

			We call  $\bar{x_{n}}$ the composition  
			$\ _n I \circ i_{n-1}$; it is  
			clearly an injection.\\

	  We can now subsume our conclusion about the closure of affine braids, which can be viewed as the "affine" version of theorem  \ref{1_4_13}. \\

			\begin{corollaire} \label{1_4_14}\textbf{Any $B$-braid, hence any affine 
			one, can be viewed as a braid in some $A$-type braid group. So, we can define the closure of an affine braid as the closure of its image under $\bar{x_{j}}$ 
			for some positive integer $j$. This injection means that any condition forcing any two affine braids to have the same closure is a consequence of the two Markov conditions. }
			\end{corollaire}
						
\medskip

			\begin{proposition} \label{1_4_15}
				Let $x$ be any affine braid   in $B(\tilde{A_{n-1}})$ for some $ 2 \leq n$. Then: \\
			     \begin{enumerate}[label=\arabic*), font=\normalsize, font=\color{black}, leftmargin=2cm, parsep=0cm, itemsep=0.25cm, topsep=0cm]
					\item  Given $y$ in $B(\tilde{A_{n-1}})$ such that $\phi_{n}(y)=x$ then $\hat{y} = \hat{x}$ (in other terms $\hat{}$ is invariant under the action of the Dynkin automorphism).
					\item  $\widehat{x a_{n+1}} = \hat{x}$.\\
				\end{enumerate}
			\end{proposition}
						
			\begin{demo}
			
				Suppose $\phi_{n}(y)=x$. That is equivalent to saying that $ D_{n} x D^{-1}_{n} = y $ in $B(\tilde{A_{n}})$. But by the first move   $ D_{n}  x D^{-1}_{n} \leftrightarrow x D^{-1}_{n} D_{n} = x $. Thus $x \leftrightarrow y $.\\
			
				On the other hand $ x a_{n+1} = x \sigma^{-1}_{n}a_{n}\sigma_{n} = x a_{n} \sigma_{n}a_{n}^{-1} \leftrightarrow a_{n}^{-1} x a_{n} \sigma_{n}$, by the first move. But $a_{n}^{-1} x a_{n} \sigma_{n} \leftrightarrow  a_{n}^{-1} x a_{n} $, by the second move. Hence we are reduced to $a_{n}^{-1} x a_{n} \leftrightarrow x a_{n}a_{n}^{-1}  = x $, by the first move, which means $ x a_{n+1} \leftrightarrow x$.
			
			\end{demo}
			
			Set $\widehat{B(A_{n})}$ to be $\left\{ \widehat{b}; b \in B(A_{n})  \right\}$.\\ 
			
			Now we reformulate our description of what we called "$B$-links", defined as the closures of $B$-braids. In $\mathds{R}^{3}$ the $B$-links are those links in which there is an oriented unknotted  fixed circle. Now we can talk about $\widehat{B(A_{n})}$ defined above without ambiguity (so as for $B\widehat{(\tilde{A_{n}})}$). It is clear geometrically, that considering the fixed circle as a circle among the others gives the way in which  $\widehat{B(B_{n})}$ is contained in $\widehat{B(A_{n})}$, while links in which there is no strings around the fixed circle, gives the inclusion of  $\widehat{B(A_{n-1})}$ in $\widehat{B(A_{n})}$. It is well known that $B$-links  represent the links in a solid torus, the string which makes the round around the hole of the torus represents $\hat{t}$ or $\hat{t^{-1}}$, which depends of course on  the orientation of the string. \\
							
			In the same spirit we see that affine links are $B$-links in which the number of positive rounds equals the number of negative rounds, of course around the fixed circle. Links which do not make rounds are counted here, actually they describe $\widehat{B(\tilde{A_{n}})}$ containing $\widehat{B(A_{n})}$.\\
            \begin{corollaire}
			
			Suppose that $L_{1}$ and $L_{2}$ are affine (oriented) links. Suppose that $l_{1}$ and $l_{2}$ are two affine braids such that $\widehat{l_{1}} = L_{1}$ and $\widehat{l_{2}} = L_{2}$. Then $L_{1}$ and $L_{2}$ are isotopic if and only if $l_{1}$ and $l_{2}$ are equivalent in  the sense of Markov, when being viewed in $\cup_{1 \leq i}B(A_{i})$ , that is: if and only if  $\bar{x_{j}}(l_{1})$ and $\bar{x_{j}}(l_{2})$ are equivalent for some $j$.
			
			\end{corollaire}
			
			\begin{demo}
			
			See theorem 5.2 in \cite{Geck_Lambropoulou_1997} with proposition \ref{1_4_15} .

			\end{demo}

		  	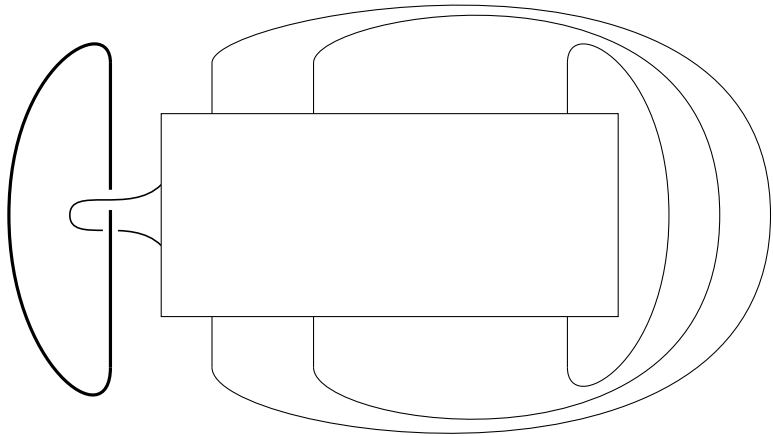
\begin{figure}[h]
		  		\centering
				  \begin{tikzpicture}
				  \begin{scope}[scale = 1.35]
\draw[draw= white, double = black, line width = 0.1cm] (-0.5, 1.2) .. controls +(-0.3, 0.3) and +(0, -0.3) .. (-1.4, 1.5) .. controls +(0,0.3) and +(-0.3,-0.3) .. (-0.5,1.8); 
\draw [white, line width = 0.2cm] (-1,0)-- +(0,1.5);
\draw[very thick] (-1,0)-- +(0,1.55);
\draw[very thick] (-1,1.75)-- (-1,3);
\draw (0,0) -- +(0,3);
\draw (1,0) -- +(0,3);
\draw (3.5,0) -- +(0,3);
\filldraw[fill= white] (-0.5,0.5)-- (4, 0.5) -- (4, 2.5) -- (-0.5, 2.5) -- cycle; 
\draw (3.5,0).. controls +(0,-0.5) and +(0,-1.5) .. (4.5,1.5).. controls +(0,1.5) and +(0, 0.5) .. (3.5,3);
\draw (1,0).. controls +(0,-0.6) and +(0,-2.7) .. (5,1.5).. controls +(0,2.7) and 
+(0, 0.5) .. (1,3);
\draw (0,0).. controls +(0,-0.7) and +(0,-3) .. (5.5,1.5).. controls +(0,3) and 
+(0, 0.5) .. (0,3);
\draw[very thick] (-1,0).. controls +(0,-0.7) and +(0,-1.5) .. (-2,1.5).. controls +(0,1.5) and  +(0, 0.5) .. (-1,3);
\end{scope}

				    \end{tikzpicture}
				\caption{Affine closure}
			\end{figure}

	\vspace{2cm}

\renewcommand{\refname}{REFERENCES}

\end{document}